\documentclass[10pt,letterpaper]{amsart}
\usepackage[top=22mm, bottom=22mm, left=18mm, right=18mm]{geometry}
\usepackage{newpxtext,newpxmath}
\usepackage{bm}
\usepackage{charter}
\usepackage[OT2,OT1]{fontenc}
\usepackage[latin9]{inputenc}
\usepackage{units}

\usepackage{amstext}
\usepackage{amsthm}

\usepackage{cancel}

\makeatletter

\pdfpageheight\paperheight
\pdfpagewidth\paperwidth

\DeclareTextSymbolDefault{\textquotedbl}{T1}
\providecommand{\tabularnewline}{\\}

\numberwithin{equation}{section}
\numberwithin{figure}{section}
\theoremstyle{plain}
\newtheorem{thm}{\protect\theoremname}[section]
\theoremstyle{definition}
\newtheorem{defn}[thm]{\protect\definitionname}
\theoremstyle{plain}
\newtheorem{prop}[thm]{\protect\propositionname}
\theoremstyle{plain}
\newtheorem{lem}[thm]{\protect\lemmaname}
\theoremstyle{remark}
\newtheorem{rem}[thm]{\protect\remarkname}
\theoremstyle{plain}
\newtheorem{cor}[thm]{\protect\corollaryname}

\usepackage{graphicx}
\usepackage{amscd}
\usepackage{enumerate}
\usepackage{url}
\usepackage{xfrac}

\newcommand{\cyr}{
\renewcommand\rmdefault{wncyr} \renewcommand\sfdefault{wncyss} \renewcommand\encodingdefault{OT2} \normalfont
\selectfont
}
\DeclareTextFontCommand{\textcyr}{\cyr}    

   \def\@settitle
{
\begin{center} \baselineskip14\p@\relax \LARGE\@title \end{center}
}

\usepackage[normalem]{ulem}
\usepackage[mathscr]{eucal}
\numberwithin{equation}{section}

\usepackage{verbatim}
\usepackage{amscd}
\usepackage[mathscr]{eucal}
\usepackage[all]{xy}
\usepackage{enumerate}
\usepackage{color}
\usepackage{colortbl}

\input xyimport.tex

\title{{\bf Constructing and characterizing prime $\mathbb{Q}$-Fano threefolds of genus one and with six $\sfrac{1}{2}(1,1,1)$-singularities via key varieties}}

\author{Hiromichi Takagi}

\address{Department of Mathematics, Gakushuin University, 
Mejiro, Toshima-ku, Tokyo 171-8588, Japan}
\email{hiromici@math.gakushuin.ac.jp}

\newcommand{\sE}{\mathcal{E}}

\newcommand{\sN}{\mathcal{N}}

\newcommand{\sO}{\mathcal{O}}

\newcommand{\sQ}{\mathcal{Q}}

\newcommand{\sU}{\mathcal{U}}

\newcommand{\mA}{\mathbb{A}}
\newcommand{\mC}{\mathbb{C}}
\newcommand{\mF}{\mathbb{F}}

\newcommand{\mN}{\mathbb{N}}
\newcommand{\mP}{\mathbb{P}}
\newcommand{\mQ}{\mathbb{Q}}

\newcommand{\Bs}{\mathrm{Bs}\,}

\newcommand{\Sing}{\mathrm{Sing}\,}

\newcommand{\rank}{\mathrm{rank}\,}

\theoremstyle{definition}

\newtheorem{con}{Convention}[section]
\numberwithin{equation}{section}

 \newcounter{myparagraph}[subsection]

\makeatother

\providecommand{\corollaryname}{Corollary}
\providecommand{\definitionname}{Definition}
\providecommand{\lemmaname}{Lemma}
\providecommand{\propositionname}{Proposition}
\providecommand{\remarkname}{Remark}
\providecommand{\theoremname}{Theorem}

\begin{document}
\maketitle 
\begin{abstract}
We consider the classification problem of prime $\mQ$-Fano 3-folds
with at most $\nicefrac{1}{2}(1,1,1)$-singularities, which was initiated
in \cite{Taka2}. We construct two distinct classes of such 3-folds
with genus one and six $\nicefrac{1}{2}(1,1,1)$-singularities, each
equipped with a prescribed Sarkisov link. Our method involves constructing
certain higher-dimensional $\mQ$-Fano varieties $\Sigma$, referred
to as key varieties, by extending the Sarkisov links to higher dimensions.
We prove that each such 3-fold $X$ arises as a linear section of
the corresponding key variety $\Sigma$, and conversely, any general
linear section of $\Sigma$ yields such an $X$. Various geometric
properties of the key varieties $\Sigma$ are also investigated and
clarified.
\end{abstract}

\maketitle
\markboth{$\mQ$-Fano 3-fold of genus $1$ with six $\nicefrac{1}{2}$-singularities}{Hiromichi
Takagi} {\small{\tableofcontents{}}}{\small\par}

2020\textit{ Mathematics subject classification}: 14J45, 14E30, 14N05


\section{\textbf{Introduction\label{sec:Introduction}}}

\subsection{Background for classification of $\mQ$-Fano 3-folds\label{subsec:Background-for-classification}}

A complex projective variety is called \textit{a $\mathbb{Q}$-Fano
variety} if it is a normal variety with only terminal singularities
and the anticanonical divisor is ample.

The minimal model theory provides a framework for the birational classification
of projective varieties by adopting those with numerically well-behaved
canonical divisors as selected models. $\mathbb{Q}$-Fano varieties
constitute one such class of models, making their biregular classification
a crucial aspect of the birational classification of projective varieties.

The classifications of $\mathbb{Q}$-Fano varieties in dimensions
one and two were settled by the early 20th century around the time
when the classical minimal model theory due to the Italian school
was completed. In dimension three, although the minimal model theory
was established toward the end of the last century, the classification
of three-dimensional $\mathbb{Q}$-Fano varieties $(\ensuremath{\mathbb{Q}}$-Fano
3-folds) remains incomplete. In the nonsingular case, however, the
classification of Fano 3-folds has been extensively studied from both
classical and modern perspectives throughout the last century, culminating
in a comprehensive and mature classification (\cite{F}, \cite{I},
\cite{MM}, \cite{Mu2}). The fascination with smooth Fano 3-folds
is widely recognized among geometers.

Following the success of the classification of nonsingular Fano 3-folds,
we aim to complete the classification of singular $\mathbb{Q}$-Fano
3-folds while also clarifying how certain rich geometric structures
arise in them.

\subsection{Prime $\mQ$-Fano $3$-folds with only $\nicefrac{1}{2}(1,1,1)$-singularities\label{subsec:Prime--Fano--folds 1/2}}

\subsubsection{\textbf{Background and previous research}}

Motivated by this, since the publication of \cite{Taka2}, we have
been investigating the classification of $\mathbb{Q}$-Fano 3-folds
with only $\nicefrac{1}{2}(1,1,1)$-singularities (abbreviated as
$\nicefrac{1}{2}$-singularities) and this has led to substantial
results. Since the $\nicefrac{1}{2}$-singularity is the simplest
among non-Gorenstein 3-fold terminal singularities, the classification
of $\mathbb{Q}$-Fano $3$-folds with only such singularities naturally
follows that of nonsingular ones. Moreover, $\mQ$-Fano 3-folds with
$\nicefrac{1}{2}$-singularities seem to provide a glimpse into the
overall picture of general non-Gorenstein $\mQ$-Fano 3-folds.

To be more precise, {[}ibid.{]} provides a detailed study of the classification
of $\mathbb{Q}$-Fano 3-folds $X$ with only $\nicefrac{1}{2}$-singularities,
genus $\geq2$, and whose anticanonical divisor generates the divisor
class group up to numerical equivalence. A $\mathbb{Q}$-Fano $3$-fold
satisfying the last condition on the anticanonical divisor is called
\textit{a prime $\mathbb{Q}$-Fano $3$-fold.} \textit{The genus $g(X)$
}of a $\mQ$-Fano $3$-fold $X$ is defined to be $h^{0}(-K_{X})-2$.
This classification method extends Takeuchi\textquoteright s approach
for nonsingular prime Fano 3-folds \cite{Take}: selecting any $\nicefrac{1}{2}$-singularity
of $X$ and performing the blow-up at that point, we can construct
a so-called Sarkisov link, which we explain in Subsection \ref{subsec:Sarkisov-link}.
The classification process is then divided into the following two
steps:

\vspace{3pt}

\noindent \textit{Step 1.} Listing the possible invariants associated
with the various objects appearing in the Sarkisov link for $X$.

\noindent \textit{Step 2.} Constructing the Sarkisov link with those
invariants and proving the existence of $X$.

\vspace{3pt}

In {[}ibid.{]}, the first step of the classification was carried out
with high accuracy, and for many of the resulting possibilities, it
was demonstrated that the Sarkisov link with the corresponding invariants
could indeed be constructed, thereby completing the second step of
the classification in many cases. 

Under this situation, we would like to step into the classification
of prime $\mQ$-Fano $3$-folds of genus at most $1$. Fortunately,
in {[}ibid.{]} we have shown that the above method also applies when
the condition that genus $\geq2$ is weakened to the one that genus
$\geq-1$ and anticanonical volume $\geq1$, allowing the first step
of the classification to be carried out. The details of this proof
will be published separately in \cite{Taka7}. 

\subsubsection{\textbf{The aim of this paper}}

In this paper, we complete the second step of the classification for
two specific classes of prime $\mQ$-Fano 3-folds introduced in Subsection
\ref{subsec:The-case-of} (No.$\,$12960 in \cite{GRDB}) by explicitly
constructing the Sarkisov links for them. Our approach is to extend
the Sarkisov links to higher dimensions, thereby obtaining certain
higher-dimensional $\mQ$-Fano varieties $\Sigma$, which we call
\textit{key varieties}. We then demonstrate that each 3-fold $X$
under consideration can be realized as a linear section of its associated
key variety $\Sigma$, and conversely, that a general linear section
of $\Sigma$ recovers such an $X$. This establishes a precise correspondence
between these 3-folds and their key varieties, providing a characterization
of the $3$-folds $X$. This method originates from Mukai\textquoteright s
studies on nonsingular prime Fano 3-folds \cite{Mu2}. The extension
from $X$ to $\Sigma$ is achieved through a precise analysis of the
invariants of the objects appearing in the Sarkisov link. The characterizations
of $\mQ$-Fano $3$-folds as linear sections of their key varieties
have already been obtained in \cite{Taka4} in several cases classified
in \cite{Taka2}. This type of characterization of $\mQ$-Fano $3$-folds
is expected to play a significant role for constructing and describing
their moduli spaces (cf.~\cite{Mu1}).

The paper \cite{Taka7} includes some examples for which the second
step of the classification has been also completed, but many examples
have not reached that stage. Moreover, as seen in the examples treated
in this paper, there are also cases where the construction is not
straightforward, making it appropriate to investigate them in separate
papers. Therefore, it may be seen that the classification in the cases
of genus at most one has just begun.

\subsection{Sarkisov link\label{subsec:Sarkisov-link}}

Let $X$ be a prime $\mQ$-Fano $3$-fold $X$ with $g(X)\geq-1$
and $(-K_{X})^{3}\geq1$, and $\widehat{f}\colon\widehat{X}\to X$
the blow-up at an arbitrarily chosen $\nicefrac{1}{2}$-singularity
$x$. By \cite[Part I, Thm.4.1 and Prop.4.2]{Taka2}, $-K_{\widehat{X}}$
is nef and big. Therefore the anticanonical model $\widehat{g}\colon\widehat{X}\to\overline{X}$
is defined. Although $\widehat{g}$ can be a crepant divisorial contraction,
we do not consider this case. Hereafter, we suppose that\textit{ $\widehat{g}$
is a flopping contraction.} Then, by playing the so-called $2$-ray
game (note that the Picard number of $\widehat{X}$ is $2$), we obtain
a uniquely determined birational map $\widehat{X}\dashrightarrow\widetilde{X}$
isomorphic in codimension $1$ such that $\widetilde{X}$ has a non-small
extremal contraction $\widetilde{f}\colon\widetilde{X}\to Y$. The
birational map $\widehat{X}\dashrightarrow\widetilde{X}$ is a composition
of flops and flips in this order in general, but again for simplicity,
we assume that \textit{$\widehat{X}\dashrightarrow\widetilde{X}$
is one flop.} Hence $\widetilde{X}$ has the flopped contraction $\widetilde{g}\colon\widetilde{X}\to\overline{X}$
to $\widehat{g}$. The\textit{ Sarkisov link} starting from $\widehat{f}$
is the following diagram summarizing the above construction:

\begin{equation}\label{eq:3-foldSarkisov} \xymatrix{& \widehat{X}\ar@{-->}[rr]\ar[dl]_{\widehat{f}}\ar[dr]^{\widehat{g}} &  &\widetilde{X}\ar[dr]^{\widetilde{f}}\ar[dl]_{\widetilde{g}}\\
X& &\overline{X} &  & Y.}
\end{equation}The variety $\overline{X}$ is usually called the \textit{midpoint}. 

Even if we do not restrict $\widehat{f}$ in this manner, there are
various situations where such a diagram can be constructed for some
birational morphism to $X$. The advantage of imposing the above restriction
on $\widehat{f}$, however, is that $\widehat{X}$ and $\widetilde{X}$
still possess only $\nicefrac{1}{2}$-singularities (this holds even
if we allow flips for $\widehat{X}\dashrightarrow\widetilde{X}$).
Then the extremal contractions appearing in $\widetilde{f}$, have
a $3$-fold with only $\nicefrac{1}{2}$-singularities as their source,
which is favorable because their classification has already been completed.
By substituting the numerical data extracted from this classification
into the equations obtained from the numerical restrictions of the
diagram (\ref{eq:3-foldSarkisov}), we can form Diophantine equations
and figure out a finite set of solutions. This is Step 1 as explained
in Subsection \ref{subsec:Prime--Fano--folds 1/2}. As stated in there,
this step turns out to be highly precise if $g(X)\geq2$, hence many
of these solutions can indeed be realized as Sarkisov diagrams of
some prime $\mathbb{Q}$-Fano 3-folds (Step 2).

\subsection{The two classes\label{subsec:The-case-of}}

Now we introduce the two classes of prime $\mQ$-Fano $3$-folds which
we investigate in this paper.
\begin{defn}
\label{def:dep2}We denote by ${\rm dP}_{2}$ a double cover of $\mP^{3}$
along a quartic surface $B$ with only $5$ Du Val singularities of
type $A_{1}$ or $A_{2}$ (for brevity, we call such singularities
\textit{$A_{\leq2}$-singularities}). In other words, a ${\rm dP}_{2}$
is a del Pezzo $3$-fold of degree $2$ and with $5$ singularities
each of which is an ordinary double point or a singularity analytically
isomorphic to $\{xy+z^{2}+w^{3}=0\}\subset\mC^{4}$. We denote by
$\mathsf{q}_{1},\dots,\mathsf{q}_{5}$ the $5$ singularities of $B$
and by $\mathsf{\widetilde{q}}_{1},\dots,\widetilde{\mathsf{q}}_{5}$
the corresponding $5$ singularities of the ${\rm dP}_{2}$. By \cite[Thm.4.3]{HP},
\cite[Thm.3]{Ch} and \cite[Thm.4.20 (1)]{KOPP}, a ${\rm dP}_{2}$
is factorial (this condition implies that it has Picard number $1$).
We denote by $\sO_{{\rm dP}_{2}}(1)$ the pull-back of $\sO_{\mP^{3}}(1)$
by the double cover.
\end{defn}

\begin{defn}
\label{def:TwoTypes}Let $X$ be a prime $\mQ$-Fano $3$-fold of
genus one and with six $\nicefrac{1}{2}$-singularities. Let $x$
be a $\nicefrac{1}{2}$-singularity of $X$, and $\widehat{f}\colon\widehat{X}\to X$
the blow-up at $x$ and $\widehat{E}$ the $\widehat{f}$-exceptional
divisor. We construct the Sarkisov link (\ref{eq:3-foldSarkisov})
starting from $\widehat{f}$. We denote by $\widetilde{L}$ the $\widetilde{f}$-pull-back
of the ample generator of ${\rm Pic}\,Y$, and by $\widehat{E}_{\widetilde{X}}\subset\widetilde{X}$
the strict transform of $\widehat{E}$. Assume that the equality 
\begin{equation}
\widetilde{L}=2(-K_{\widetilde{X}})-\widehat{E}_{\widetilde{X}}\label{eq:L=00003D2(-K)-E}
\end{equation}
 holds. Under this condition, we define the following two classes
of $\mQ$-Fano $3$-folds:

\vspace{3pt}

\noindent \textbf{Type $R$:} Assume that $Y\simeq\mP^{3}$ and $\widetilde{f}\colon\widetilde{X}\to Y$
is a divisorial contraction whose center $C\subset\mP^{3}$ is a curve
of degree $12$ and $p_{g}(C)=7$ and with \textbf{$5$ }non-planar
triple points $\mathsf{q}_{1},\dots,\mathsf{q}_{5}$ each of which
can be resolved by the blow-up at the point. 

In this case, we say that \textit{$X$ is of Type $R$ with respect
to $x$.} Note that $X$ is rational in this case (the name Type $R$
reflects this fact).

\vspace{3pt}

\noindent \textbf{Type $I\!R$:} Assume that $Y$ is a ${\rm dP}_{2}$
and $\widetilde{f}\colon\widetilde{X}\to Y$ is a divisorial contraction
whose center $C\subset Y$ is an elliptic curve through the five singularities
$\widetilde{\mathsf{q}}_{1},\dots,\mathsf{\widetilde{q}}_{5}$ of
$Y$ and of degree $4$ with respect to $\sO_{{\rm dP}_{2}}(1)$.

In this case, we say that \textit{$X$ is of Type $I\!R$ with respect
to $x$}. Note that $X$ is irrational in this case if $Y$ has only
ordinary double points by \cite[Thm.1.2]{CPS} (the name Type $I\!R$
reflects this fact). 

\vspace{3pt}

For the Sarkisov link (\ref{eq:3-foldSarkisov}) as above of a Type
$R$ or Type $I\!R$ $\mQ$-Fano $3$-fold, we denote by $\widetilde{E}$
the $\widetilde{f}$-exceptional divisor. When there is no risk of
confusion, the part \textquotedbl with respect to a $\nicefrac{1}{2}$-singularity\textquotedbl{}
after Type $R$ or Type $I\!R$ is omitted.
\end{defn}

\subsection{Main result}
\begin{defn}
Let $\Sigma$ be a projective variety in a weighted projective space.
We say that a subvariety $X\subset\Sigma$ is a \textit{linear section
}of $\Sigma$ if $X$ is defined in $\Sigma$ by weighted homogeneous
hypersurfaces $H_{1},\dots,H_{c}$ of weight one, where $c$ is the
codimension of $X$ in $\Sigma$.
\end{defn}

\begin{defn}
We say that five points in $\mP^{3}$ are \textit{in  general position}
if no four are in a plane nor no three are in a line. 
\end{defn}

\begin{thm}
\label{thm:main}

Assume that the $5$ points $\mathsf{q}_{1},\dots,\mathsf{q}_{5}\in\mP^{3}$
are in general position.

\begin{enumerate}[$(1)$]

\item There exists a $10$-dimensional terminal $\mQ$-factorial
$\mQ$-Fano variety $\mathscr{R}$ in $\mP(1^{10},2^{5})$ with an
action of the symmetric group $\mathfrak{S}_{6}$ satisfying $-K_{\mathscr{R}}=\sO_{\mathscr{R}}(8)$
such that the following $(1\text{{\rm -}}1)$ and $(1\text{-}2)$
hold:

\begin{enumerate}[$({1}\text{-}1)$] 

\item Any Type $R$ $\mQ$-Fano $3$-fold $X$ is a linear section
of $\mathscr{R}$ of codimension $7$. 

\item Conversely, any linear section of $\mathscr{R}$ of codimension
$7$ with only six $\nicefrac{1}{2}$-singularities is a Type $R$
$\mQ$-Fano $3$-fold with respect to any $\nicefrac{1}{2}$-singularity.

\end{enumerate}

\item For any ${\rm dP}_{2}$, there exists a $4$-dimensional terminal
$\mQ$-factorial $\mQ$-Fano variety $\mathscr{IR}$ (depending on
${\rm dP}_{2}$) in $\mP(1^{4},2^{5})$ satisfying $-K_{\mathscr{IR}}=\sO_{\mathscr{IR}}(2)$
such that the following $(2\text{{\rm -}}1)$ and $(2\text{-}2)$
hold: 

\begin{enumerate}[$({2}\text{-}1)$] 

\item Any Type $I\!R$ $\mQ$-Fano $3$-fold $X$ with $Y={\rm dP}_{2}$
is a linear section of $\mathscr{IR}$ of codimension $1$. 

\item Conversely, any linear section of $\mathscr{I}\mathscr{R}$
of codimension $1$ with only six $\nicefrac{1}{2}$-singularities
is a Type $I\!R$ $\mQ$-Fano $3$-fold with respect to\textbf{ }one
specified $\nicefrac{1}{2}$-singularity of $\mathscr{IR}$.

\end{enumerate}

\end{enumerate}
\end{thm}

In Theorem \ref{thm:main}, there are actually two cases for Type
$I\!R$ $\mQ$-Fano $3$-folds. We distinguish and treat the two cases
separately as General and Special Cases, under a common framework,
in Section \ref{sec:TypeIRGen}. We hope to investigate in the future
how this distinction will emerge in the description of the moduli
space of Type $I\!R$ prime $\mQ$-Fano $3$-folds. In the case of
Type $I\!R$, we will also consider the problem of whether the result
(2-2) holds for \textit{any} $\nicefrac{1}{2}$-singularity of $\mathscr{IR}$.

A type $R$ $\mQ$-Fano 3-fold $X$ obtained in Theorem \ref{thm:main}
turns out to be anticanonically embedded of codimension 4 in the sense
of \cite{Taka6} (cf.$\,$Proposition \ref{prop:Gorcodim4}). While
many such examples have been constructed in \cite{CD}, \cite{Taka6},
such an $X$ as in Theorem \ref{thm:main} is obtained for the first
time in this paper. As for the Sarkisov link for such an $X$, it
can be seen that $X\dashrightarrow\overline{X}$ is a so-called Type
I projection (cf.$\,$Subsection \ref{subsec:Constructing--via unproj}).
It follows from \cite{Ca} that it is neither of Tom nor Jerry Type
since the Sarkisov link for $X$ constructed in this paper does not
appear in the list in the big table of {[}ibid.{]} (see also {[}ibid.,
Rem.4.2{]}). This property also highlights a new aspect of such an
$X$, not found in previous examples. As for a Type $I\!R$ $\mQ$-Fano
3-folds, similar properties are expected but we cannot conclude them
since the explicit equations of a Type $I\!R$ $\mQ$-Fano 3-fold
are not available in this paper.

For $\mQ$-Fano $3$-folds in Theorem \ref{thm:main}, a particularly
noteworthy observation is that, in both cases, the final step of the
Sarkisov link involves five points in general position in a three-dimensional
projective space, which is closely related to the \textit{Segre cubic.}
The construction of the aforementioned varieties $\Sigma$ is made
possible by revealing this underlying classical projective geometry
of the Sarkisov link (see Section\ref{sec:Preliminaries}). 

\subsection{A future research plan}

In the companion paper \cite{Taka8}, we prove a result similar to
Theorem \ref{thm:main}, replacing the assumption that the five points
$\mathsf{q}_{1},\dots,\mathsf{q}_{5}\in\mP^{3}$ are in general position
with a weaker condition: namely, that they are in \textit{quasi-general
position}, meaning that exactly four of them lie in a plane, and no
three lie on a line.

\subsection{Structure of the paper}

In Section \ref{sec:Preliminaries}, we elucidate the projective geometry
underlying the final map of the Sarkisov link, $\widetilde{f}\colon\widetilde{X}\to Y$,
for $\mQ$-Fano 3-folds of Type $R$ and Type $I\!R$. 

The midpoint $\overline{X}$ of the Sarkisov link for $X$ of Type
$R$ or Type $I\!R$, based on its numerical data, is suggested to
have a relation with a weighted Grassmannian, and this turns out to
be the case. With this in mind, in Section \ref{sec:ConeG(2,5)},
we conduct a somewhat general investigation of weighted Grassmannians.

Sections \ref{sec:TypeR} and \ref{sec:TypeIRGen} are devoted to
the constructions of the key varieties $\mathscr{R}$ and $\mathscr{IR}$
for $\mQ$-Fano 3-folds of Type $R$ and Type $I\!R$, respectively.
We prove Theorem \ref{thm:main} in Sections \ref{subsec:Proof-of-Theorem R}
and \ref{sec:Proof-of-Theorem(2)} for Type $R$ and Type $I\!R$
$\mQ$-Fano $3$-folds respectively. For a more detailed explanation
of the content of each section, the reader is referred to its beginning.

\vspace{10pt}

\noindent\textbf{ Notation and Conventions:}

\begin{itemize}

\item

$A(t):={\footnotesize \left(\begin{array}{ccccc}
0 & t_{12} & t_{13} & t_{14} & t_{15}\\
 & 0 & t_{23} & t_{24} & t_{25}\\
 &  & 0 & t_{34} & t_{35}\\
 &  &  & 0 & t_{45}\\
 &  &  &  & 0
\end{array}\right)}:$ the $5\times5$ skew-symmetric matrix with $t_{ij}$ as $(i,j)$-entries
$(i<j)$.

\item ${\rm Pf}_{ijkl}(t)$: the $4\times4$ Pfaffian of the skew-symmetric
matrix obtained by keeping the $i$, $j$, $k$, and $l$ rows and
columns of $A(t)$. We also call ${\rm Pf}_{1234}(t)$, ${\rm Pf}_{1235}(t)$,
${\rm Pf}_{1245}(t)$, ${\rm Pf}_{1345}(t)$, and ${\rm Pf}_{2345}(t)$
\textit{the five Pl\"ucker relations of $A(t)$.}

\item $\mP(\sE)$: the projectivization of the vector bundle $\sE$
on a variety. Note that we do not follow the Grothendieck notation;
in some references, $\mP(\sE)$ is written as $\mP_{*}(\sE)$.

\item \textit{A $\nicefrac{1}{2}$-singularity}: the singularity
analytically isomorphic to the one appearing at the origin of the
quotient of $\mC^{n}$ with coordinate $x_{1},\dots,x_{n}$ by the
involution $(x_{1},\dots,x_{n})\mapsto(-x_{1},\dots,-x_{n})$.

\item \textit{A $c({\rm G}(2,5))$-singularity}: the singularity
analytically isomorphic to the one appearing at the vertex of the
cone over the Grassmannian ${\rm G}(2,5)$.

\item $\mP(x_{1},\dots,x_{n})$ and $\mA(x_{1},\dots,x_{n})$: the
projective space and the affine space with coordinates $x_{1},\dots,x_{n}$,
respectively.

\item Let $V$ be a vector space and $\bm{v}$ an element of $V$.
We denote by $[\bm{v}]$ the point of $\mP(V)$ corresponding to $\bm{v}$.
We also apply this convention for a weighted projectivization of $V$.

\item \textit{A flopping contraction}: a projective morphism $f\colon A\to B$
with only connected fibers satisfying (1) $A$ is $\mQ$-factorial
and has only terminal singularities, (2) the $f$-exceptional locus
has codimension at least 2, (3) any $f$-exceptional curve is numerically
trivial for $K_{A}$, and (4) the relative Picard number of $f$ is
$1$. 

\item\textit{ Flopped contractions}: if $f\colon A\to B$ and $f'\colon A'\to B$
are flopping contractions, and $A$ and $A'$ are not isomorphic over
$B$, we say $f$ is the \textit{flopped contraction }for $f'$ and
vice versa. The induced birational maps $A\dashrightarrow A'$ and
$A'\dashrightarrow A$ are called the flops.

\item \textit{Divisor and sheaf convention}: We denote a Weil divisor
on a normal variety by the same notation as the corresponding divisorial
sheaf. Although the notation feels somewhat awkward, we think that
it is better than unnecessarily increasing the number of symbols.

\end{itemize}

All computations in this paper are elementary; those that could not
be carried out by hand were performed using built-in commands in \textit{Mathematica}
\cite{W}.

\vspace{3pt}

\noindent \textbf{Acknowledgment}: This work was supported by a Grant-in-Aid
for Scientific Research (C), JSPS, Japan, Grant Number 25K06923.

\section{\textbf{Projective geometry underlying the Sarkisov links\label{sec:Preliminaries}}}

The main purpose of this section is to reveal the projective geometry
underlying the morphism $\widetilde{f}\colon\widetilde{X}\to Y$ in
the Sarkisov diagram (\ref{eq:3-foldSarkisov}) for Type $R$ and
$I\!R$ $\mQ$-Fano $3$-folds $X$. In Subsections \ref{subsec:Segre}
and \ref{subsec:Birational-geometry-of dp2}, we exhibit basic results
with respect to the Segre cubic ${\rm S}_{3}$ and its double cover
$\mathcal{S}$, respectively. Subsection \ref{subsec:CS} is one of
cores of this paper, where we describe the birational model of the
curve $C$ appearing in the Sarkisov diagram (\ref{eq:3-foldSarkisov})
on ${\rm S}_{3}$ in the case of Type $R$ and on $\mathcal{S}$ in
the case of Type $I\!R$. To describe them, we introduce the projective
bundle $\mathsf{R}$ over $\mathsf{S}_{3}$ in the case of Type $R$
and the projective bundle $\mathsf{IR}$ over $\mathcal{S}$ in the
case of Type $I\!R$, which serve as essential objects for constructing
the key variety of $X$.

\subsection{Birational geometry of the Segre cubic ${\rm S}_{3}$\label{subsec:Segre}}
\begin{defn}
\label{def:Segre}Let $\mathsf{q}_{1},\dots,\mathsf{q}_{5}$ be $5$
points of $\mP^{3}$ in general position. It is easy to check that
$\mathsf{q}_{1},\dots,\mathsf{q}_{5}$ are unique up to projective
equivalence, and the space of quadratic forms vanishing at them is
$5$-dimensional. Moreover the image of the rational map $\mP^{3}\dashrightarrow\mP^{4}$
defined by $5$ linearly independent quadratic forms vanishing at
them is a cubic $3$-fold. The cubic is classically called the Segre
cubic primal (cf.~\cite{D}), which we call simply \textit{the Segre
cubic} and denote it by ${\rm S}_{3}$. 

In the case of Type $R$, we use the specified equation of ${\rm S}_{3}$
as in (\ref{eq:S3}), while in the case of Type $I\!R$, we do not
specify an equation of ${\rm S_{3}}$. We denote by 
\[
\mu_{S}\colon\mP^{3}\dashrightarrow{\rm S}_{3}
\]
 the induced rational map.
\end{defn}

\begin{prop}
\label{prop:Segre} Assume that the $5$ points $\mathsf{q}_{1},\dots,\mathsf{q}_{5}\in\mP^{3}$
are in general position. The following holds:

\begin{enumerate}[$(1)$]

\item The rational map $\mu_{{\rm S}}\colon\mP^{3}\dashrightarrow{\rm S}_{3}$
is birational, and  its indeterminacy locus is resolved by the  blow-up ${\rm bl_{\mP^{3}}}\colon{\rm Bl}\,\mP^{3}\to\mP^{3}$
at the $5$ points $\mathsf{q}_{1},\dots,\mathsf{q}_{5}$. 

\item Let $\mu'_{{\rm S}}\colon{\rm Bl}\,\mP^{3}\to{\rm S}_{3}$
be the induced morphism. Let $\ell_{ij}$ be the lines through $\mathsf{q}_{i}$
and $\mathsf{q}_{j}$, and $\ell'_{ij}$ their strict transforms on
${\rm Bl}\,\mP^{3}\,(1\leq i<j\leq5)$. The $\mu'_{{\rm S}}$-exceptional
locus coincides with $\cup_{1\leq i<j\leq5}\ell_{ij}'$ and $\Sing{\rm S}_{3}$
consists of the $10$ points $\mathsf{q}_{ij}:=\mu'_{{\rm S}}(\ell'_{ij})$,
which are ordinary double points. 

\item Let $\mathsf{P}_{i}'$ be the exceptional divisor of ${\rm Bl}\,\mP^{3}\to\mP^{3}$
over $\mathsf{q}_{i}\,(1\leq i<j\leq5)$. The images of $\mathsf{P}_{i}'$
on ${\rm S}_{3}$ are planes, which we denote by $\mathsf{P}_{i}$.

\end{enumerate}
\end{prop}

\begin{proof}
Everything follows from computations based on explicit equations of
$\mathsf{q}_{1},\dots,\mathsf{q}_{5}$ (see Subsection \ref{subsec:LR and S3}
for an example of explicit equations). We only mention that $\mu_{{\rm S}}\colon\mP^{3}\dashrightarrow{\rm S}_{3}$
factors through the blow-up ${\rm bl_{\mP^{3}}}\colon{\rm Bl}\,\mP^{3}\to\mP^{3}$
at the $5$ points $\mathsf{q}_{1},\dots,\mathsf{q}_{5}$ since they
are defined scheme theoretically by $5$ linearly independent quadratic
forms vanishing at them.
\end{proof}

\subsection{Birational geometry of a ${\rm dP}_{2}$\label{subsec:Birational-geometry-of dp2}}

As for a ${\rm dP}_{2}$ as in Definition \ref{def:dep2}, we subsequently
assume that \textit{the $5$ singularities $\mathsf{q}_{1},\dots,\mathsf{q}_{5}$
of $B$ are in general position. }

We fix a ${\rm dP}_{2}$ in this subsection. Let $\rho_{\mP}\colon{\rm dP}_{2}\to\mP^{3}$
and $\rho_{{\rm Bl}\,\mP}\colon{\rm BldP}_{2}\to{\rm Bl}\,\mP^{3}$
be the double covers branched along $B$ and along the strict transform
$B'$ of $B$ respectively. Since $B$ has only $A_{\leq2}$-singularities,
$B'$ is smooth, and hence so is ${\rm BldP}_{2}.$ We can locally
check that the Stein factorization of the composite ${\rm BldP}_{2}\overset{\rho_{{\rm Bl}\,\mP}}{\to}{\rm Bl}\,\mP^{3}\overset{{\rm bl}_{\mP^{3}}}{\to}\mP^{3}$
is ${\rm BldP}_{2}\overset{{\rm bl}_{{\rm dP}_{2}}}{\to}{\rm dP}_{2}\overset{\rho_{\mP}}{\to}\mP^{3}$,
where ${\rm bl}_{{\rm dP}_{2}}\colon{\rm BldP}_{2}\to{\rm dP}_{2}$
is the blow-up at $\widetilde{\mathsf{q}}_{1},\dots,\widetilde{\mathsf{q}}_{5}$.
Let ${\rm BldP}_{2}\overset{\nu_{\mathcal{S}}'}{\to}\mathcal{S}\overset{\rho_{{\rm S}}}{\to}{\rm S}_{3}$
be the Stein factorization of the composite ${\rm BldP}_{2}\overset{2:1}{\to}{\rm Bl}\,\mP^{3}\overset{\mu_{{\rm S}}'}{\to}{\rm S}_{3}$.
Again we can locally check that $\rho_{{\rm S}}\colon\mathcal{S}\to{\rm S}_{3}$
is the double cover branched along the $\mu_{{\rm S}}$-birational
image $B_{{\rm S}}$ of $B$. In summary, we obtain the following
diagram, where all the vertical arrows are double covers:

\begin{equation}\label{eq:doublecover} \xymatrix{& \rm{Bl dP}_2\ar[dl]_{\rm{bl}_{\rm{dP}_2}}\ar[dr]^{\nu'_{\mathcal{S}}} \ar[d]_{\rho_{\rm{Bl}\,\mP}}& \\ \rm{dP}_2\ar[d]_{\rho_{\mP}}  & \rm{Bl} \mP^3\ar[dl]_{\rm{bl}_{\mP^3}}\ar[dr]^{\mu'_{\rm{S}}}  & \mathcal{S}\ar[d]^{\rho_{\rm{S}}}\\
\mP^3\ar@{-->}[rr]_{\mu_{\rm{S}}} & & \rm{S}_3.}
\end{equation}
\begin{prop}
\label{prop:Sing calS}The singular locus of $\mathcal{S}$ coincides
with $\cup_{1\leq i<j\leq5}\rho_{{\rm S}}^{-1}(\mathsf{q}_{ij})$.
\end{prop}

\begin{proof}
The assertion follows since ${\rm {BldP}_{2}}$ is smooth and the
image of the $\nu'_{\mathcal{S}}$-exceptional locus is $\cup_{1\leq i<j\leq5}\rho_{{\rm S}}^{-1}(\mathsf{q}_{ij})$. 
\end{proof}
\begin{prop}
\label{prop:23CI}The surface $B_{{\rm S}}$ is the complete intersection
of the Segre cubic ${\rm S}_{3}$ and a unique quadric $3$-fold,
which we denote by $Q_{B}$. In particular, $B_{{\rm S}}$ is also
a $K3$ surface with only Du Val singularities.
\end{prop}

\begin{proof}
Since $\mu_{{\rm S}}$ is defined by the linear system of quadrics
through $\mathsf{q}_{1},\dots,\mathsf{q}_{5}$, and $B$ is a quartic
surface with double points at $\mathsf{q}_{1},\dots,\mathsf{q}_{5}$,
it holds that the surface $B_{{\rm S}}$ is a quadric section of ${\rm S}_{3}$.
Thus the first assertion follows. The second assertion follows since
$K_{B_{{\rm S}}}\sim0$ and $B_{{\rm S}}$ is birational to the $K3$
surface $B$.
\end{proof}
The following proposition leads to an essential case division of the
key varieties of Type $I\!R$ $\mQ$-Fano $3$-folds (we refer to
Sections \ref{sec:TypeIRGen} for the treatment of the respective
cases).
\begin{prop}
\label{prop:rkQB} Let $Q_{B}$ be the quadric $3$-fold as in Proposition
\ref{prop:23CI}. The rank of $Q_{B}$ is $4$ or $5$. If the rank
of $Q_{B}$ is $4$, then the singular point of $Q_{B}$ (the vertex
of $Q_{B}$) is not contained in ${\rm S}_{3}$.
\end{prop}

\begin{proof}
Since $B_{{\rm S}}=Q_{B}\cap{\rm S}_{3}$ is an irreducible surface,
we have $\rank Q_{B}\geq3$. Thus we have only to exclude the case
where $\rank Q_{B}=3$. Assume for contradiction that $\rank Q_{B}=3$.
Then $Q_{B}$ is singular along a line $l$, and hence $B_{{\rm S}}={\rm S}_{3}\cap Q_{B}$
is singular at the points in ${\rm S}_{3}\cap l$. By Proposition
\ref{prop:Sing calS}, this implies that ${\rm S}_{3}$ is singular
at the points in ${\rm S}_{3}\cap l$. This contradicts Proposition
\ref{prop:23CI} since then both $Q_{B}$ and ${\rm S}_{3}$ are singular
at the points in ${\rm S}_{3}\cap l$ and then $B_{{\rm S}}$ has
non-hypersurface singularities at them.

The preceding discussion also shows that when $\rank Q_{B}=4$, the
singular point of $Q_{B}$ is not contained in ${\rm S}_{3}$.
\end{proof}
We denote by $\mathcal{Q}_{B}$ the double cover of $\mP^{4}$ branched
along $Q_{B}$. Note that $\mathcal{Q}_{B}$ is a quadric 4-fold in
$\mP^{5}$ and $\rank\mathcal{Q}_{B}=5$ and $6$ when $\rank Q_{B}=$
$4$ and $5$ respectively. We regard $\mathcal{S}$ as a cubic section
of $\mathcal{Q}_{B}$. Indeed, we may write 
\begin{equation}
\mathcal{S}=\mathcal{Q}_{B}\cap c({\rm S}_{3}),\label{eq:S=00003DQS}
\end{equation}
where $c({\rm S}_{3})$ is the cone over ${\rm S}_{3}$ in $\mP^{5}$. 

Let $\mathsf{Q}_{i}'$ be the $\rho_{{\rm Bl}\,\mP}$-inverse image
of $\mathsf{P}_{i}'\,(1\leq i\leq5)$, and $\mathsf{Q}_{i}$ the image
of $\mathsf{Q}_{i}'$ on $\mathcal{Q}_{B}$. We note that $\mathsf{Q}_{i}'$
is the ${\rm bl}_{{\rm dP}_{2}}$-exceptional divisor over $\widetilde{\mathsf{q}}_{i}$,
and hence is $\mP^{1}\times\mP^{1}$ (resp.~$\mP(1^{2},2)$) if $\mathsf{\mathsf{q}}_{i}$
is an $A_{1}$ (resp.~$A_{2}$)-singularity of $B$. We also note
that $\mathsf{Q}_{i}$ is a quadric surface isomorphic to $\mathsf{Q}_{i}'$,
which is the double cover of $\mathsf{P}_{i}$ and is a linear section
of the quadric $4$-fold $\mathcal{Q}_{B}$ of codimension $2$.

\subsection{Birational model of $C\subset Y$ via the Segre cubic\label{subsec:CS}}

\subsubsection{\textbf{Type $R$\label{subsec:TypeR}}}
\begin{lem}
\label{lem:not on a quadric}For a Type $R$ $\mQ$-Fano $3$-fold
$X$, the curve $C$ as in Definition \ref{def:TwoTypes} is not in
a quadric surface.
\end{lem}

\begin{proof}
Since $-K_{\widetilde{X}}=\widetilde{f}^{*}(-K_{\mP^{3}})-\widetilde{E}=4\widetilde{L}-\widetilde{E}$,
we have 
\begin{equation}
h^{0}(\sO_{\mP^{3}}(4)-C)=h^{0}(-K_{\widetilde{X}})=h^{0}(-K_{\widehat{X}})=h^{0}(-K_{X})=3.\label{eq:4L-C}
\end{equation}
If $C$ were contained in a quadric surface, then it would hold that
\[
3=h^{0}(\sO_{\mP^{3}}(4)-C)=h^{0}(\sO_{\mP^{3}}(2)+\sO_{\mP^{3}}(2)-C)\geq h^{0}(\sO_{\mP^{3}}(2))=10,
\]
 a contradiction. 
\end{proof}
The following proposition reveals the projective geometry underlying
the morphism $\widetilde{f}\colon\widetilde{X}\to Y=\mP^{3}$ as in
(\ref{eq:3-foldSarkisov}).

\begin{prop}
\label{prop:KeyLemma} \begin{enumerate}[$(A)$]

\item For a Type $R$ $\mQ$-Fano $3$-fold $X$, let $C_{{\rm S}}\subset{\rm S}_{3}$
be the birational image of $C$. The following hold:

\begin{enumerate}[$(1)$]

\item The curve $C_{{\rm S}}$ is a smooth curve of genus $7$ and
of degree $9$.

\item There exists a cubic scroll $T\subset\mP^{4}$ such that $C_{{\rm S}}={\rm S}_{3}\cap T$.

\item Both ${\rm S}_{3}$ and $T$ are smooth along $C_{{\rm S}}$.

\item The strict transform of $C_{{\rm S}}$ on ${\rm Bl}\,\mP^{3}$
and $\ell'_{ij}\,(1\leq i<j\leq5)$ are disjoint.

\end{enumerate}

\item Conversely, let $C_{{\rm S}}$ be the smooth curve which is
the intersection between ${\rm S}_{3}$ and a cubic scroll $T$, and
$C$ the birational image of ${\rm C}_{{\rm S}}$ by $\mu_{{\rm S}}\colon\mP^{3}\dashrightarrow{\rm S}_{3}$.
It holds that $C$ is a curve of degree $12$ and is smooth outside
$\mathsf{q}_{1},\dots,\mathsf{q}_{5}$ and has non-planar triple points
at them. Moreover, by the blow-up of $\mP^{3}$ at them, the triple
points of $C$ are resolved.

\end{enumerate}
\end{prop}

\begin{proof}
We show the assertion (A). 

\vspace{3pt}

\noindent (1) and (2). Since $\sO_{\mP^{3}}(2)-\mathsf{q}_{1}-\cdots-\mathsf{q}{}_{5}$
is the pull-back of $\sO_{{\rm S}_{3}}(1)$, and $C$ has triple points
at the five points $\mathsf{q}_{1},\dots,\mathsf{q}_{5}$, it holds
that $\sO_{{\rm S_{3}}}(1)\cdot C_{{\rm S}}=2\deg C-3\times5=9$,
namely, $\deg C_{{\rm S}}=9$.

Note that the strict transform of $C$ on the blow-up ${\rm Bl}\,\mP^{3}$
of $\mP^{3}$ at $\mathsf{q}_{1},\dots,\mathsf{q}_{5}$ is a smooth
curve by the assumption of Type $R$. By (\ref{eq:4L-C}), we can
take three quartic surfaces $G_{1},G_{2},G_{3}$ containing $C$ such
that their defining equations are linearly independent. Since $C$
has nonplanar triple points at the five points $\mathsf{q}_{1},\dots,\mathsf{q}_{5}$,
$G_{i}\,(1\leq i\leq3)$ are singular at them. Since $\sO_{\mP^{3}}(4)-2\mathsf{q}_{1}-\cdots-2\mathsf{q}{}_{5}$
is the pull-back of $\sO_{{\rm S}_{3}}(2)$, we see that $C_{{\rm S}}$
is contained in the $3$ quadrics $\widetilde{G}_{i}\,(1\leq i\leq3$)
corresponding to $G_{i}$. Let us see that $C_{{\rm S}}$ is not contained
in a hyperplane of $\mP^{4}$, where $\mP^{4}$ is the ambient projective
space of ${\rm S_{3}}$. Indeed, if $C_{{\rm S}}$ were contained
in a hyperplane $\mathsf{L}$, then $C$ would be contained in the
quadric surface which is the strict transform of $\mathsf{L}\cap{\rm S}_{3}$,
a contradiction to Lemma \ref{lem:not on a quadric}. Thus the quadrics
$\widetilde{G}_{i}\,(1\leq i\leq3)$ are irreducible. If $\widetilde{G}_{1}\cap\widetilde{G}_{2}\cap\widetilde{G}_{3}$
were a curve, then it would hold that $\deg C_{S}\leq8$. This contradicts
$\deg C_{S}=9$ which we have obtained above. Thus $\widetilde{G}_{1}\cap\widetilde{G}_{2}\cap\widetilde{G}_{3}$
is $2$-dimensional. If $\widetilde{G}_{1}\cap\widetilde{G}_{2}$
were irreducible, then $\widetilde{G}_{1}\cap\widetilde{G}_{2}\cap\widetilde{G}_{3}$
must be a curve, a contradiction. Thus $\widetilde{G}_{1}\cap\widetilde{G}_{2}$
is a reducible surface. Since $C_{{\rm S}}$ is not contained in a
$3$-plane, $C_{{\rm S}}$ is not contained in a surface of degree
$\leq2$. Therefore $\widetilde{G}_{1}\cap\widetilde{G}_{2}$ is the
union of an irreducible non-degenerate cubic surface $T$ in $\mP^{4}$
and a plane, and $C_{{\rm S}}$ is contained in $T$. It is known
classically that $T$ is a scroll, i.e., $T\simeq\mP(\sO_{\mP^{1}}(-1)\oplus\sO_{\mP^{1}}(-2))$
or the cone over a twisted cubic curve. Since $C_{{\rm S}}$ is also
contained in the cubic ${\rm S}_{3}$, we see that $C_{{\rm S}}\subset T\cap{\rm S}_{3}$.
Since both $C_{{\rm S}}$ and $T\cap{\rm S}_{3}$ are curves of degree
9, we conclude that $C_{{\rm S}}=T\cap{\rm S}_{3}$. By the assumption,
we have $g(C)=7$. By intersection theory of curves on $T$, we see
that the arithmetic genus of $C_{{\rm S}}$ is $7$. Therefore $C_{{\rm S}}$
is a smooth curve of genus $7$ and degree 9.

\vspace{3pt}

\noindent (3). Since $C_{{\rm S}}$ is a Cartier divisor on $T$
defined by ${\rm S}_{3}$, $T$ is smooth along $C_{{\rm S}}$. This
implies that $T$ is a local complete intersection near $C_{{\rm S}}$.
Then, since $C_{{\rm S}}$ is cut out from ${\rm S}_{3}$ by $T$,
${\rm S}_{3}$ is smooth along $C_{{\rm S}}$. 

\vspace{3pt}

\noindent (4). (3) implies (4) by Proposition \ref{prop:Segre} (2).

\vspace{3pt}

The assertion (B) follows by reversing the proof of (A). Therefore,
we omit the detailed proof, but note the following point. Since $C_{{\rm S}}={\rm S}_{3}\cap T$
and $\mathsf{P}_{i}\subset{\rm S}_{3}$, we have $C_{{\rm S}}\cap\mathsf{P}_{i}=T\cap\mathsf{P}_{i}$.
Then, since $C_{\mathsf{S}}\not\subset\mathsf{P}_{i}$, we see that
$C_{{\rm S}}\cap\mathsf{P}_{i}$ is a $0$-dimensional noncolinear
subscheme of length 3.
\end{proof}
The following proposition shows that, starting from a general situation
of Proposition \ref{prop:KeyLemma} (B) (see the assumption of (A)
(b) of the proposition below), we can conversely construct the morphism
$\widetilde{f}\colon\widetilde{X}\to\mP^{3}$:
\begin{prop}
\label{prop:GenExOf f'} 

\begin{enumerate}[$(A)$]

\item

There exists a cubic scroll $T\subset\mP^{4}$ satisfying

\begin{enumerate}[$(a)$]

\item $C_{{\rm S}}:=T\cap{\rm S}_{3}$ is a smooth curve of genus
$7$ and degree 9 (this implies that $T$ and ${\rm S}_{3}$ are smooth
near $C_{{\rm S}}$ by the proof of Proposition \ref{prop:KeyLemma}
(A) (3)), and

\item $T\cap\mathsf{P}_{i}$ consists of three distinct points for
each $i=1,\dots,5$, which are not colinear since $T$ has no $3$-secant
lines.

\end{enumerate}

\item Under the situation of $(A)$, let $C$ be the birational image
of $C_{{\rm S}}$ by $\mu_{{\rm S}}\colon\mP^{3}\dashrightarrow{\rm S}_{3}$
(see Proposition \ref{prop:KeyLemma} (B) for the properties of $C$).
It holds that there exists a $3$-fold $\widetilde{X}$ with only
five $\nicefrac{1}{2}$-singularities, and an extremal divisorial
contraction $\widetilde{f}\colon\widetilde{X}\to\mP^{3}$ (in particular,
$\rho(\widetilde{X})=2$) contracting a prime divisor $\widetilde{E}$
to $C$ and satisfying the following properties:

\begin{enumerate}[$(1)$]

\item The $\widetilde{f}$-images of the five $\nicefrac{1}{2}$-singularities
of $\widetilde{X}$ are $\mathsf{q}_{1},\dots,\mathsf{q}_{5}\in\mP^{3}$.

\item $\widetilde{f}$ coincides with the blow-up along $C$ outside
the fibers over $\mathsf{q}_{1},\dots,\mathsf{q}_{5}$.

\item The fiber over $\mathsf{q}_{i}$ consists of three copies of
$\mP^{1}$ of degree $\nicefrac{1}{2}$ with respect to $-K_{\widetilde{X}}$.

\item It holds that $(-K_{\widetilde{X}})^{3}=5/2$ and $(-K_{\widetilde{E}})^{2}=-138$.

\end{enumerate}

\end{enumerate}

\end{prop}

\begin{proof}
The assertion (A) follows by \cite[2 Thm. (ii)]{Kl} since $\mP^{4}$
is homogeneous under the action of ${\rm PGL}_{5}$ (note also that
the remark in the proof of Proposition \ref{prop:KeyLemma} (B)). 

Let us show the assertion (B) by constructing a birational model $\widetilde{X}$
of ${\rm S}_{3}$ in the following steps (cf.~\cite{H}): 

\vspace{3pt}

\noindent \textit{Step 1.} Take the small resolution $\mu_{{\rm S}}'\colon{\rm Bl}\,\mP^{3}\to{\rm S}_{3}$
as in Proposition \ref{prop:Segre}. Let $\mathsf{P}_{i}',C'_{{\rm S}}\subset{\rm Bl}\,\mP^{3}$
be the strict transforms of $\mathsf{P}_{i}\,(1\leq i\leq5),\,C_{{\rm S}}$
respectively. Recall that $\mathsf{P}_{i}'$ is a $\mu_{{\rm S}}$-exceptional
divisor and $\mathsf{\ensuremath{P}}_{i}'\simeq\mP^{2}$ by Proposition
\ref{prop:Segre}.

\vspace{3pt}

\noindent \textit{Step 2.} Take the blow-up ${\rm (Bl}\,\mP^{3})'$
of ${\rm Bl}\,\mP^{3}$ along $C_{{\rm S}}'$. Let $\widetilde{E}'$
be the exceptional divisor over $C_{{\rm S}}'$ and $\mathsf{P}_{i}''\subset{\rm (Bl}\,\mP^{3})'$
the strict transforms of $\mathsf{P}_{i}'\,(1\leq i\leq5).$ Note
that $\mathsf{P}_{i}''\to\mathsf{P}_{i}'$ is the blow-up along $\mathsf{P}_{i}'\cap C_{{\rm S}}'$
which consists of three non-colinear points.

\vspace{3pt}

\noindent \textit{Step 3.} Flop the strict transforms in $\mathsf{P}_{i}''$
of the three lines in $\mathsf{P}_{i}'$ each of which intersects
$C'_{S}$ at two points. Denote the flop by $({\rm Bl}\,\mP^{3})'\dashrightarrow({\rm Bl}\,\mP^{3})^{+}$,
and by $\mathsf{P}_{i}^{+}\subset({\rm Bl}\,\mP^{3})^{+}$ the strict
transforms of $\mathsf{P_{i}''}\,(1\leq i\leq5)$. Then it holds that
$\mathsf{P}_{i}^{+}\simeq\mP^{2}$ and $\mathsf{P}_{i}^{+}|_{\mathsf{P_{i}^{+}}}=\sO_{\mP^{2}}(-2)$
(we remark that $\mathsf{P}_{i}'\dashrightarrow\mathsf{P}_{i}^{+}$
is the Cremona transformation).

\vspace{3pt}

\noindent \textit{Step 4.} Take the birational morphism $({\rm Bl}\,\mP^{3})^{+}\to Y'$
contracting $\mathsf{P}_{i}^{+}\,(1\leq i\leq5)$ to $\nicefrac{1}{2}$-singularities.

\vspace{3pt}

By the construction, a birational morphism $\widetilde{f}\colon\widetilde{X}\to\mP^{3}$
is induced. It holds that $\widetilde{f}$ contracts the strict transform
$\widetilde{E}$ of $\widetilde{E}'$ to the curve $C$. We can check
the properties (1)--(4) of $\widetilde{f}$ in a elementary way by
the construction. We remark that the formula $(-K_{\widetilde{E}})^{2}=-138$
in (4) also follows from \cite[Prop.7.1]{Taka3} (interpreting $E=\widetilde{E}$,
$g(\overline{C})=7$, $m=0$ and $m'=5$).
\end{proof}
\begin{rem}
A crutial point of this explicit construction is that we can verify
$\widetilde{X}$ has Picard number 2.
\end{rem}

To construct $\widetilde{f}\colon\widetilde{X}\to\mP^{3}$, including
non-general situations of Proposition \ref{prop:KeyLemma} (B), and
further to construct a Type $R$ $\mQ$-Fano 3-fold $X$ by reversing
the Sarkisov diagram, it is appropriate to first extend the blow-up
of ${\rm S}_{3}$ along $C_{{\rm S}}$ to a higher dimensional situation,
and reverse the higher dimensional Sarkisov diagram. 

The following proposition indicate how to extend the blow-up of ${\rm S}_{3}$
along $C_{{\rm S}}$ to a higher dimensional situation:
\begin{prop}
\label{prop:The-cubic-scroll zero} Let $V$ be a $5$-dimensional
vector space and $V^{*}$ its dual. We denote the coordinates of $V^{\oplus2}$
by $\mathsf{x}_{j},\mathsf{y}_{j}\,(1\leq j\leq5)$, where $\mathsf{x}_{j}$
and $\mathsf{y}_{j}$ are the coordinates of the first and the second
factors of $V^{\oplus2}$ respectively. Let $T$ be a cubic scroll
in $\mP(V^{*})\simeq\mP^{4}$. There exist 
\[
\mathsf{H}_{i}:=\{\sum_{j=1}^{5}a_{ij}\mathsf{x}_{j}+\sum_{j=1}^{5}b_{ij}\mathsf{y}_{j}=0\}\subset\mP(V^{*})\times\mP(V^{\oplus2})\,(1\leq i\leq7)
\]
 such that the natural morphism $\mP(\Omega_{\mP(V^{*})}(1)^{\oplus2})\cap\bigcap_{i=1}^{7}\mathsf{H}_{i}\to\mP(V^{*})$
is the blow-up of $\mP(V^{*})$ along $T$.

In particular, if $C_{{\rm S}}={\rm S}_{3}\cap T$ is a smooth curve,
then $\mP(\Omega_{\mP(V^{*})}(1)^{\oplus2}|_{{\rm S}_{3}})\cap\bigcap_{i=1}^{7}\mathsf{H}_{i}\to{\rm S}_{3}$
is the blow-up of ${\rm S}_{3}$ along $C_{{\rm S}}$.
\end{prop}

\begin{rem}
To consider $T\subset\mP(V^{*})$ (rather than $T\subset\mP(V)$)
is compatible with the subsequent constructions. 
\end{rem}

\begin{proof}
We may explicitly find $\bigcap_{i=1}^{7}\mathsf{H}_{i}$ and check
that $\mP(\Omega_{\mP(V^{*})}(1)^{\oplus2})\cap\bigcap_{i=1}^{7}\mathsf{H}_{i}\to\mP(V^{*})$
is the blow-up of $\mP(V^{*})$ along $T$ as follows:

\vspace{3pt}

\noindent \textit{Case 1: $T\simeq\mP(\sO_{\mP^{1}}(-1)\oplus\sO_{\mP^{1}}(-2))$.}
By the ${\rm GL}(V)$-action on $\mP(V^{*})$, we may assume that
\[
T=\left\{ \rank\left(\begin{array}{ccc}
\mathsf{z}_{1} & \mathsf{z}_{2} & \mathsf{z}_{3}\\
\mathsf{z}_{2} & \mathsf{z}_{4} & \mathsf{z}_{5}
\end{array}\right)\leq1\right\} .
\]
Then we can take 
\[
\bigcap_{i=1}^{7}\mathsf{H}_{i}=\{\mathsf{x}_{1}=0,\mathsf{x}_{2}+\mathsf{y}_{1}=0,\mathsf{x}_{3}=0,\mathsf{x}_{4}+\mathsf{y}_{2}=0,\mathsf{x}_{5}+\mathsf{y}_{3}=0,\mathsf{y}_{4}=0,\mathsf{y}_{5}=0\}.
\]

\vspace{3pt}

\noindent \textit{Case 2: $T$ is the cone over a twisted cubic curve.
}By the ${\rm GL}(V)$-action on $\mP(V^{*})$, we may assume that
\[
T=\left\{ \rank\left(\begin{array}{ccc}
\mathsf{z}_{1} & \mathsf{z}_{2} & \mathsf{z}_{3}\\
\mathsf{z}_{2} & \mathsf{z}_{3} & \mathsf{z}_{4}
\end{array}\right)\leq1\right\} .
\]
Then we can take 

\[
\bigcap_{i=1}^{7}\mathsf{H}_{i}=\{\mathsf{x}_{1}=0,\mathsf{x}_{2}+\mathsf{y}_{1}=0,\mathsf{x}_{3}+\mathsf{y}_{2}=0,\mathsf{x}_{4}+\mathsf{y}_{3}=0,\mathsf{x}_{5}=0,\mathsf{y}_{4}=0,\mathsf{y}_{5}=0\}.
\]
\end{proof}
Since the projective bundle $\mP(\Omega_{\mP(V^{*})}(1)^{\oplus2}|_{{\rm S}_{3}})$
plays important roles subsequently, let us give it a name:

\begin{equation}
\mathsf{R}:=\mP(\Omega_{\mP(V^{*})}(1)^{\oplus2}|_{{\rm S}_{3}}).\label{eq:mathsfR}
\end{equation}

\subsubsection{\textbf{Type $I\!R$\label{subsec:TypeIR}}}

The approach in this subsection is similar to that in Subsection \ref{subsec:TypeR}.
\begin{lem}
\label{lem:elliptic quartic}

For a Type $I\!R$ $\mQ$-Fano $3$-fold $X$, let $C_{\mP}$ be the
image of $C$ as in Definition \ref{def:TwoTypes} by the double covering
$\rho_{\mP}\colon Y={\rm dP}_{2}\to\mP^{3}$ as in (\ref{eq:doublecover}).
The induced morphism $C\to C_{\mathbb{P}}$ is an isomorphism and
$C_{\mathbb{P}}$ is a quartic elliptic curve and is the complete
intersection of two quadric surfaces.
\end{lem}

\begin{proof}
Suppose that $C\to C_{\mP}$ is a double covering. Then $C_{\mP}$
is a conic, and hence, by the Hurwitz formula, the branch locus of
$C\to C_{\mP}$ consists of 4 points. The branch locus of $C\to C_{\mP}$
contains, however, at least the five points $\mathsf{q}_{1},\dots,\mathsf{q_{5}}$,
a contradiction. Therefore $C\to C_{\mP}$ is birational. Since $-K_{\widetilde{X}}=\widetilde{f}{}^{*}(-K_{Y})-\widetilde{E}=2\widetilde{L}-\widetilde{E}$,
we have 
\[
h^{0}(\sO_{Y}(2)-C)=h^{0}(-K_{\widetilde{X}})=h^{0}(-K_{\widehat{X}})=h^{0}(-K_{X})=3.
\]
If $C_{\mP}$ were contained in a plane, then 
\[
3=h^{0}(\sO_{Y}(1)+\sO_{Y}(1)-C)\geq h^{0}(\sO_{Y}(1))=4,
\]
a contradiction. Since $h^{0}(\sO_{Y}(2)-C)=3$ and $h^{0}(\sO_{\mP^{3}}(2))=h^{0}(\sO_{Y}(2))-1$,
we see that $h^{0}(\sO_{\mP^{3}}(2)-C_{\mP})\geq2.$ Let $Q_{1},Q_{2}$
be two quadrics containing $C_{\mP}$ such that the equations of $Q_{1},Q_{2}$
are linearly independent. Since $C_{\mP}$ is not contained in a plane,
$Q_{1}$ and $Q_{2}$ are irreducible. Therefore $Q_{1}\cap Q_{2}$
is 1-dimensional. Since $C_{\mP}\subset Q_{1}\cap Q_{2}$ and $\deg C_{\mP}=\deg Q_{1}\cap Q_{2}=4$,
we must have $C_{\mP}=Q_{1}\cap Q_{2}$. Since the arithmetic genus
of $Q_{1}\cap Q_{2}$ is 1, and $C\to C_{\mP}$ is birational, this
is an isomorphism. 
\end{proof}
\begin{prop}
\label{prop:CS cubic}\begin{enumerate}[$(A)$]

\item For a Type $I\!R$ $\mQ$-Fano $3$-fold, let $C_{\mathsf{\mathcal{S}}}$
be the birational image of $C$ as in Definition \ref{def:TwoTypes}
by the birational map $Y={\rm dP}_{2}\dashrightarrow\mathcal{S}$
as in (\ref{eq:doublecover}). It holds that $C_{\mathsf{\mathcal{S}}}$
is a cubic elliptic curve and is the section of the cubic cone $c({\rm S}_{3})$
by a plane in the quadric $4$-fold $\mathcal{Q}_{B}$ as constructed
in Subsection \ref{subsec:Birational-geometry-of dp2}. Moreover,
$C_{\mathcal{S}}$ does not contain a singular point of $\mathcal{S}$,
and when $\mathsf{Q}_{i}$ is singular, $C_{\mathcal{S}}$ does not
contain the singular point of $\mathsf{Q}_{i}$.

\item Conversely, for any ${\rm d}P_{2}$ as in Definition \ref{def:dep2},
let $C_{\mathcal{S}}$ be the cubic elliptic curve which is the section
of $c({\rm S}_{3})$ by a plane $\Pi_{C}$ in the quadric $4$-fold
$\mathcal{Q}_{B}$ (we define $\mathcal{S}$ and $\mathcal{Q}_{B}$
as in Subsection \ref{subsec:Birational-geometry-of dp2}). Let $C$
be the birational image of ${\rm C}_{\mathcal{S}}$ by ${\rm dP}_{2}\dashrightarrow\mathcal{S}$
as in (\ref{eq:doublecover}). It holds that $C$ is a quartic elliptic
curve through $\mathsf{\widetilde{q}}_{1},\dots,\widetilde{\mathsf{q}}_{5}$.

\end{enumerate}
\end{prop}

\begin{proof}
We show the assertion (A). Let $C_{{\rm S}}\subset{\rm S}_{3}$ be
the birational image of $C_{\mP}$ by $\mu_{{\rm S}}\colon\mP^{3}\dashrightarrow{\rm S}_{3}$.
Since $\sO_{\mP^{3}}(4)-2\mathsf{q}_{1}-\cdots-2\mathsf{q}{}_{5}$
is the pull-back of $\sO_{{\rm S}_{3}}(2)$, and $C$ contains the
five points $\mathsf{q}_{1},\dots,\mathsf{q}_{5}$, it holds that
$\sO_{{\rm S}_{3}}(1)\cdot C_{{\rm S}}=2\deg C_{\mP}-5=3$. Therefore,
since $C_{\mP}$ is an elliptic curve by Lemma \ref{lem:elliptic quartic},
$C_{{\rm S}}$ is a cubic elliptic curve. Since $C\to C_{\mP}$ is
an isomorphism by Lemma \ref{lem:elliptic quartic}, so is $C_{\mathcal{S}}\to C_{{\rm S}}$.
Therefore $C_{\mathcal{S}}$ is also a cubic elliptic curve. Note
that the plane containing $C_{\mathcal{S}}$ is contained in $\mathcal{Q}_{B}$
since otherwise the intersection between the plane and $\mathcal{Q}_{B}$
is a conic and cannot contain $C_{\mathcal{S}}$. Therefore, by (\ref{eq:S=00003DQS}),
$C_{\mathcal{S}}$ is the intersection between the plane and the cubic
$c({\rm S}_{3})$. 

Since the smooth curve $C_{{\rm S}}$ is cut out by a plane from ${\rm S}_{3}$,
it is a complete intersection, thus it does not pass through the singular
points $\mathsf{q}_{ij}\,(1\leq i<j\leq5)$ of ${\rm S}_{3}$ and
then $C_{\mathcal{S}}$ does not contain a point of $\rho_{{\rm S}}^{-1}(\mathsf{q}_{ij})$.
Therefore $C_{\mathcal{S}}$ is disjoint from ${\rm Sing\,\mathcal{S}}$
by Proposition \ref{prop:Sing calS}. Since $C$ is smooth, the strict
transform of $C$ on ${\rm BldP}_{2}$ intersects $\mathsf{Q}_{i}'$
at one smooth point simply (note that $\mathsf{Q}_{i}'$ is a ${\rm bl}_{{\rm dP}_{2}}$-exceptional
divisor). Therefore $C_{\mathcal{S}}$ does not contain the singular
point of $\mathsf{Q}_{i}$ when $\mathsf{Q}_{i}$ is singular. Hence
the last assertion of (A) follows. 

The assertion (B) follows by reversing the proof of (A). Therefore,
we omit the detailed proof, but note the following point. Since $C_{{\rm S}}=c({\rm S}_{3})\cap\Pi_{C}$
and $\mathsf{Q}_{i}\subset c({\rm S}_{3})$, we have $C_{{\rm S}}\cap\mathsf{Q}_{i}=\Pi_{C}\cap\mathsf{Q}_{i}$.
Note that $\Pi_{C}\subset\mathcal{Q}_{B}$ and $\mathsf{Q}_{i}$ is
a linear section of $\mathcal{Q}_{B}$. Therefore, since $C_{\mathsf{S}}\not\subset\mathsf{Q}_{i}$,
we see that $C_{{\rm S}}\cap\mathsf{Q}_{i}$ consists of a smooth
point of $\mathsf{Q}_{i}$.
\end{proof}
As can be seen from the following proposition, in the case of Type
$I\!R$, the morphism $\widetilde{f}\colon\widetilde{X}\to{\rm dP}_{2}$
can be conversely constructed from any situation of Proposition \ref{prop:CS cubic}
(B).
\begin{prop}
\label{prop:GenExOf f'-IR} For any ${\rm dP}_{2}$ as in Definition
\ref{def:dep2}, we define $\mathcal{S}$ and $\mathcal{Q}_{B}$ as
in Subsection \ref{subsec:Birational-geometry-of dp2}. The following
assertions hold:

\begin{enumerate}[$(A)$]

\item

There exists a plane $\Pi_{C}$ in $\mathcal{Q}_{B}$ satisfying $C_{{\rm \mathcal{S}}}:=\mathcal{S}\cap\Pi_{C}=c({\rm S}_{3})\cap\Pi_{C}$
is a cubic elliptic curve.

\item Under the situation of $(A)$, let $C$ be the birational image
of $C_{{\rm \mathcal{S}}}$ by ${\rm dP}_{2}\dashrightarrow\mathcal{S}$
(see Proposition \ref{prop:CS cubic} (B) for the properties of $C$).
There exists a $3$-fold $\widetilde{X}$ with only five $\nicefrac{1}{2}$-singularities,
and an extremal divisorial contraction $\widetilde{f}\colon\widetilde{X}\to{\rm dP}_{2}$
(in particular, $\rho(\widetilde{X})=2$) contracting a prime divisor
$\widetilde{E}$ to $C$ and satisfying the following properties:

\begin{enumerate}[$(1)$]

\item The images of the five $\nicefrac{1}{2}$-singularities of
$\widetilde{X}$ by the composite $\widetilde{X}\to{\rm dP}_{2}\to\mP^{3}$
are $\mathsf{q}_{1},\dots,\mathsf{q}_{5}\in\mP^{3}$.

\item $\widetilde{f}$ coincides with the blow-up along $C$ outside
the fibers over $\mathsf{\widetilde{q}}_{1},\dots,\mathsf{\widetilde{q}}_{5}$.

\item An irreducible component of the fiber over $\mathsf{\widetilde{q}}_{i}$
is $\mP^{1}$ of degree $\nicefrac{1}{2}$ with respect to $-K_{\widetilde{X}}$.
If $\mathsf{q}_{i}$ is an $A_{1}$ (resp. $A_{2}$)-singularity of
$B$, then the fiber over $\mathsf{\widetilde{q}}_{i}$ is the union
of two copies of $\mP^{1}$ (resp. a double fiber).

\item It holds that $(-K_{\widetilde{X}})^{3}=5/2$ and $(-K_{\widetilde{E}})^{2}=-10$.

\end{enumerate} 

\end{enumerate} 
\end{prop}

\begin{proof}
Note that, when $\rank\mathcal{Q}_{B}=6$, $\mathcal{Q}_{B}$ is homogeneous
under the action of ${\rm SL}_{4}$, so the assertion (A) follows
by applying \cite[2 Thm. (ii)]{Kl} on $\mathcal{Q}_{B}$. When $\rank\mathcal{Q}_{B}=5$,
$\mathcal{Q}_{B}\setminus\{v\}$ is homogeneous under the action of
${\rm Sp}_{4}\times\mC^{*}$, so the assertion (A) follows similarly. 

Let us show the assertion (B) by constructing a birational model $\widetilde{X}$
of $\mathcal{S}$ in the steps below. Note that, by the comment in
the proof of Proposition \ref{prop:CS cubic} (B), $C_{{\rm S}}\cap\mathsf{Q}_{i}$
consists of a smooth point of $\mathsf{Q}_{i}$.

\vspace{3pt}

\noindent \textit{Step 1.} We take the birational model ${\rm BldP}_{2}$
of $\mathcal{S}$ as in the diagram (\ref{eq:doublecover}). Let $C'_{\mathcal{S}}\subset{\rm BldP}_{2}$
be the strict transform of $C_{\mathcal{S}}$. Recall that $\mathsf{Q}_{i}'$
is a ${\rm {bl}_{{\rm {dP}_{2}}}}$-exceptional divisor and $\mathsf{\ensuremath{Q}}_{i}'\simeq\mathsf{Q}_{i}$.

\vspace{3pt}

\noindent \textit{Step 2.} We take the blow-up ${\rm (BldP}_{2})'$
of ${\rm BldP}_{2}$ along $C_{\mathcal{S}}'$. Let $\widetilde{E}'$
be the exceptional divisor over $C_{\mathcal{S}}'$ and $\mathsf{Q}_{i}^{(2)}\subset{\rm (BldP}_{2})'$
the strict transforms of $\mathsf{Q}_{i}'\,(1\leq i\leq5),$ where
$\mathsf{Q}_{i}^{(2)}\to\mathsf{Q}_{i}'$ is the blow-up at the one
point $\mathsf{Q}_{i}'\cap C_{\mathcal{S}}'$ in the smooth locus
of $\mathsf{Q}_{i}'$.

\vspace{3pt}

Since ${\rm dP}_{2}$ is factorial as we noted in Definition \ref{def:dep2},
the relative Picard number of ${\rm bl}_{{\rm dP}_{2}}$ is $1$ near
each $\mathsf{\widetilde{\mathsf{q}}}_{i}$, hence the relative Picard
number of ${\rm (BldP}_{2})'\to{\rm dP}_{2}$ is $2$ near each $\widetilde{\mathsf{q}}_{i}$. 

\vspace{3pt}

\noindent \textit{Step 3.} We play the 2-ray game over ${\rm dP}_{2}$
near $\widetilde{\mathsf{q}}_{i}$ starting from ${\rm (BldP}_{2})'$.
We can check that the (well-known) birational morphism $\mathsf{Q}_{i}^{(2)}\to\mP^{2}$$\,(1\leq i\leq5)$
extend to a numerically $K$-trivial small contraction ${\rm (BldP}_{2})'\to{\rm (BldP}_{2})^{(2)}$
over ${\rm dP}_{2}$ which is a flopping contraction near each $\mathsf{Q}_{i}^{(2)}$.
By \cite{Ko}, there exists the birational map ${\rm (BldP}_{2})'\dashrightarrow{\rm (BldP}_{2})^{+}$
which coincides with the flop for ${\rm (BldP}_{2})'\to{\rm (BldP}_{2})^{(2)}$
near each $\mathsf{Q}_{i}^{(2)}$. We denote by $\mathsf{Q}_{i}^{+}\subset({\rm BldP_{2}})^{+}$
the strict transforms of $\mathsf{Q}_{i}^{(2)}$. Let us check $\mathsf{Q}_{i}^{+}\simeq\mP^{2}$
by constructing ${\rm (BldP}_{2})'\dashrightarrow{\rm (BldP}_{2})^{+}$
explicitly near each $\mathsf{Q}_{i}^{(2)}$. 

\textit{Case $\mathsf{Q}_{i}'\simeq\mP^{1}\times\mP^{1}$}: The flopping
curves on $\mathsf{Q}_{i}^{(2)}$ are the strict transforms $\delta_{1}$
and $\delta_{2}$ of two rulings in $\mathsf{Q}_{i}'$ through $\mathsf{Q}_{i}'\cap C_{\mathcal{S}}'$.
The flop is of Atiyah type, and then, by the standard construction
of the flop, we see that the flop induces the contraction $\mathsf{Q}_{i}^{(2)}\to\mP^{2}$
of $\delta_{1}$ and $\delta_{2}$.

\textit{Case $\mathsf{Q}_{i}'\simeq\mP(1^{2},2)$}: Let $\delta\subset\mathsf{Q}_{i}^{(2)}$
be the strict transform of the ruling of $\mathsf{Q}_{i}'$ through
one point $\mathsf{Q}_{i}'\cap C_{\mathcal{S}}'$. We can check that
$\delta$ is a flopping curve on ${\rm (BldP}_{2})'$. We construct
the flop of $\delta$ as follows (this is nothing but the two-tiered
Pagoda construction of the flop of a $(0,-2)$-curve given in \cite{R1}):
Let $\beta\colon{\rm (BldP}_{2})''\to{\rm (BldP}_{2})'$ be the blow-up
along $\delta$ and $\mathsf{Q}_{i}^{(3)}\subset{\rm (BldP}_{2})''$
the strict transform of $\mathsf{Q}_{i}^{(2)}$. Let $\delta'\subset\mathsf{Q}_{i}^{(3)}$
be the strict transform of $\delta$ by $\beta|_{\mathsf{Q}_{i}^{(3)}}\colon\mathsf{Q}_{i}^{(3)}\to\mathsf{Q}_{i}^{(2)}$
and $\gamma$ the $\beta|_{\mathsf{Q}_{i}^{(3)}}$-exceptional curve
(note that the image of $\gamma$ on $\mathsf{Q}_{i}^{(2)}$ is the
unique singular point of $\mathsf{Q}_{i}^{(2)}$). By a local computation,
we see that $(\beta|_{\mathsf{Q}_{i}^{(3)}})^{*}\delta=\delta'+\gamma$
and this coincides with the intersection between the $\beta$-exceptional
divisor and $\mathsf{Q}_{i}^{(3)}$. By this, we see that the normal
bundle of $\delta'$ in ${\rm (BldP}_{2})''$ is $\sO_{\mP^{1}}(-1)^{\oplus2}$,
hence, $\delta'$ is a flopping curve of Atiyah type. Constructing
the flop of $\delta'$ in the standard way, we see that the flop of
$\delta'$ induces the contraction $\mathsf{Q}_{i}^{(3)}\to\mathsf{Q}_{i}^{(4)}$
of $\delta'$. The strict transform of $\gamma$ on the flopped 3-fold
for $\delta'$ coincides with the intersection between the strict
transform of the $\beta$-exceptional divisor and $\mathsf{Q}_{i}^{(4)}$.
Finally we can contract the strict transform of the $\beta$-exceptional
divisor and the flop of $\delta$ is completed. Then the image of
$\mathsf{Q}_{i}^{(4)}$ on the flopped 3-fold is isomorphic to $\mP^{2}$.

\vspace{3pt}

We can check explicitly that $\mathsf{Q}_{i}^{+}|_{\mathsf{Q_{i}^{+}}}=\sO_{\mP^{2}}(-2)$
in the both cases.

\noindent \textit{Step 4.} We take the birational morphism $({\rm BldP_{2}})^{+}\to Y'$
contracting $\mathsf{Q}_{i}^{+}\,(1\leq i\leq5)$ to $\nicefrac{1}{2}$-singularities.

By the construction, a birational morphism $\widetilde{f}\colon\widetilde{X}\to{\rm dP}_{2}$
is induced. It holds that $\widetilde{f}$ contracts the strict transform
$\widetilde{E}$ of $\widetilde{E}'$ to $C$. We can check the properties
(1)--(4) of $\widetilde{f}$ in a standard way by the construction.
We remark that the formula $(-K_{\widetilde{E}})^{2}=-10$ in (4)
also follows from \cite[Prop.7.1]{Taka3} (interpreting $E=\widetilde{E}$,
$g(\overline{C})=1$, $m=5$ and $m'=0$).
\end{proof}
We finally present the proposition indicating how to extend $\widetilde{f}\colon Y\to{\rm dP}_{2}$
to a higher dimensional situation, accompanied with the case division
indispensable to construct the key variety $\mathscr{IR}$ for a Type
$I\!R$ $\mQ$-Fano $3$-fold $X$. 

Let $\Pi_{C}$ be a plane in the quadric $4$-fold $\mathcal{Q}_{B}$
such that $C_{\mathcal{S}}:=\mathcal{S}\cap\Pi_{C}$ is a cubic elliptic
curve.

\vspace{3pt}

\noindent \textbf{Case division.}

\vspace{3pt}

\noindent\textbf{General Case: }Assume that $\rank\mathcal{Q}_{B}=6$.
Then we may identify $\mathcal{Q}_{B}$ with the Grassmanian ${\rm G}(2,W)$
where $W$ is a $4$-dimensional vector space. For later convenience,
\textit{we choose this identification such that the plane $\Pi_{C}$
parameterizes lines in a fixed plane of $\mP(W)$ }(this is possible
by replacing $W$ with $W^{*}$ if necessary). Then it is well-known
that $\Pi_{C}$ is the zero scheme of a global section of the universal
subbundle $\mathcal{U}$ of rank $2$ on $\mathrm{G}(2,W)$. Therefore
\begin{equation}
C_{\mathcal{S}}\ \text{is the}\ 0\text{-scheme of a global section of the vector bundle}\ \sU|_{\mathcal{S}}.\label{eq:0schemev1}
\end{equation}
\noindent\textbf{Special Case: }Assume that $\rank\mathcal{Q}_{B}=5$.
Then $\mathcal{Q}_{B}$ is the cone over a smooth quadric 3-fold,
which we denote by $\mathcal{Q}_{B}'$. We denote by $v$ the vertex
of $\mathcal{Q}_{B}$. We may identify $\mathcal{Q}_{B}'$ with a
hyperplane section of the Grassmanian ${\rm G}(2,W)$ where $W$ is
a $4$-dimensional vector space. Let $p\colon\mathcal{Q}_{B}\setminus\{v\}\to\mathcal{Q}_{B}'$
be the natural projection. Note that $v\in\Pi_{C}$ since any plane
in a rank 5 quadric $4$-fold contains its vertex. It is also well-known
that $\Pi_{C}\setminus\{v\}$ is the zero scheme of a global section
of 
\begin{equation}
\sU':=p^{*}(\sU|_{\mathcal{Q}_{B}'}).\label{eq:U'}
\end{equation}
Therefore we may consider 
\begin{equation}
C_{\mathcal{S}}\ \text{is the}\ 0\text{-scheme of a global section of the vector bundle}\ \sU'|_{\mathcal{S}}\label{eq:0schemev2}
\end{equation}
(note that $\mathcal{S}$ does not contain $v$ by Proposition \ref{prop:rkQB}). 

\vspace{3pt}
\begin{prop}
\label{prop:YS}Let $\Pi_{C}$ be a plane in the quadric $4$-fold
$\mathcal{Q}_{B}$ such that $C_{\mathcal{S}}:=\mathcal{S}\cap\Pi_{C}$
is a cubic elliptic curve (note that $C_{\mathcal{S}}$ is disjoint
from ${\rm Sing}\,\mathcal{S}$ by Proposition \ref{prop:CS cubic}
(B)). Let $\mathsf{X}$ be the tautological divisor of \textbf{$\mP(\sU|_{\mathcal{S}})$}
when $\rank\mathcal{Q}_{B}=6$ (resp.\textbf{ $\mP(\sU'|_{\mathcal{S}})$}
when $\rank\mathcal{Q}_{B}=5$) corresponding to a global section
of \textbf{$\sU|_{\mathcal{S}}$} (resp.\textbf{ $\sU'|_{\mathcal{S}}$})
defining $C_{\mathcal{S}}$ as in (\ref{eq:0schemev1}) (resp.~(\ref{eq:0schemev2})).
Then the naturally induced morphism $\mathsf{X}\to\mathcal{S}$ is
the blow-up of $\mathcal{S}$ along $C_{\mathcal{S}}$.
\end{prop}

\begin{proof}
The assertion follows from the following general result.
\end{proof}
\begin{lem}
\label{lem:rk2 blup}Let $A$ be a smooth variety and $\sE$ a rank
$2$ vector bundle on $A$. Let $s$ be a global section of $\sE$
and $D$ the tautological divisor of $\mP(\sE)$ corresponding to
$s$. Assume that the zero locus $Z$ of $s$ is of codimension $2$.
It holds that the natural morphism $D\to A$ is the blow-up along
$Z$.
\end{lem}

\begin{proof}
We can show the assertion by taking a local trivialization of $\sE$.
\end{proof}
Let us give the projective bundles $\mP(\sU|_{\mathcal{S}})$ and
\textbf{$\mP(\sU'|_{\mathcal{S}})$} a unified name since they play
important roles subsequently: 

\begin{equation}
\mathsf{IR}:=\mP(\sU|_{\mathcal{S}})\,\,\text{{(resp.\,\ensuremath{\mP(\sU'|_{\mathcal{S}})}) if \ensuremath{\rank\mathcal{Q}_{B}=6} (resp.\,\ensuremath{\rank\mathcal{Q}_{B}=5})}}.\label{eq:mathsfIR}
\end{equation}

\subsection{Comments on how to proceed with the rest of the paper}

For Type $R$ (resp.~Type $I\!R$), Proposition \ref{prop:The-cubic-scroll zero}
(resp.~Proposition \ref{prop:YS}) indicates that the projective
bundle $\text{\ensuremath{\mathsf{R}}}$ (resp.~$\mathsf{IR}$) has
a birational model over $Y$ extending $\widetilde{f}\colon\widetilde{X}\to Y$
as in the Sarkisov diagram (\ref{eq:3-foldSarkisov}). Such a model
$\widetilde{\mathscr{R}}$ (resp.~$\mathscr{\widetilde{IR}}$) over
$Y$ actually exists (see Subsection \ref{subsec:The-quasi--bundle R}
and Proposition \ref{prop:Rtilde and P(Omega2)} for Type $R$, and
Subsection \ref{subsec:The-quasi--bundle-IR} and Proposition \ref{prop:IRtilde and P(U)S}
for Type $I\!R$). The projective bundle $\text{\ensuremath{\mathsf{R}}}$
(resp.~$\mathsf{IR}$), however, is not sufficient to define $\widetilde{\mathscr{R}}$
(resp.~$\mathscr{\widetilde{IR}}$). Fortunately, further information
can be extracted from the morphism $\widetilde{g}\colon\widetilde{X}\to\overline{X}$;
the numerical data of $\overline{X}$ indicates that it can be extended
to a suitable weighted Grassmannian. Indeed, the equations of $\widetilde{\mathscr{R}}$
(resp.~$\mathscr{\widetilde{IR}}$) are obtained by modifying a combination
of those of the weighted Grassmannian and those of the projective
bundle $\text{\ensuremath{\mathsf{R}}}$ (resp.~$\mathsf{IR}$).
Having this in mind, we exhibit a general study of weighted Grassmannians
in Section \ref{sec:ConeG(2,5)}, and specialize them for Type $R$
in Subsection \ref{subsec:The-cone GR}, and for Type $I\!R$ in Subsection
\ref{subsec:Birational-geometry-of IR Goverline}.

\section{\textbf{Cones over ${\rm G}(2,5)$} \textbf{in weighted projective
spaces} \label{sec:ConeG(2,5)}\protect 
}
\subsection{Basics}

In this section, we present constructions common to Types $R$ and
$I\!R$ related with cones over ${\rm G}(2,5)$. 

Let $V$ be a $5$-dimensional vector space and $W$ an $n$-dimensional
vector space for some $n\in\mN$. We define the weighted projective
space $\mP(1^{n},2^{10})$ considering the vector space $W\oplus\wedge^{2}V$
and putting weight 1 on the coordinates $\overline{x}_{1},\dots,\overline{x}_{n}$
of $W$, and weight 2 on the coordinates $\overline{p}_{ij}\,(1\leq i<j\leq5)$
of $\wedge^{2}V$. We also denote this by $\mP(W\oplus\wedge^{2}V)$,
the subspace $\{\overline{p}_{ij}=0\,(1\leq i<j\leq5)\}$ by $\mP(W)$,
and the subspace $\{\overline{\ensuremath{x}}_{1}=\cdots=\overline{x}_{n}=0\}$
by $\mP(\wedge^{2}V).$ Let 
\[
\overline{{\rm G}}\subset\mP(1^{n},2^{10})
\]
be the (weighted) cone over ${\rm G}(2,V)\simeq{\rm G}(2,5)$ defined
by the five Pl\"ucker relations of the anti-symmetric matrix $A(\overline{p})=(\overline{p}_{ij})$
(cf.~Notation and Conventions in the end of Section \ref{sec:Introduction}).
The following fact is obvious but we write down it to keep in mind:
\begin{prop}
\label{prop:singoverlineG}The variety $\overline{{\rm G}}$ has $\nicefrac{1}{2}$-singularities
along ${\rm G}(2,V)=\overline{{\rm G}}\cap\mP(\wedge^{2}V)$, and
$c({\rm G}(2,5))$-singularities along $\mP(W)$. 
\end{prop}

Let $S_{\overline{{\rm G}}}$ be the polynomial ring associated to
$\mP(1^{n},2^{10})$ and $R_{\overline{{\rm G}}}$ the quotient ring
of $S_{\overline{{\rm G}}}$ by the homogeneous ideal of $\overline{{\rm G}}$.
By \cite{BE}, we have the following $S_{\overline{{\rm G}}}$-free
resolution of $R_{\overline{{\rm G}}}$: 
\[
0\leftarrow R_{\overline{{\rm G}}}\leftarrow S_{\overline{{\rm G}}}\leftarrow S_{\overline{{\rm G}}}(-4)^{\oplus5}\leftarrow S_{\overline{{\rm G}}}(-6)^{\oplus5}\leftarrow S_{\overline{{\rm G}}}(-10)\leftarrow0.
\]
 From this we deduce 
\begin{equation}
-K_{\overline{{\rm G}}}=\sO_{\overline{{\rm G}}}(n+10),\label{eq:KG}
\end{equation}
and 
\begin{equation}
\deg\overline{{\rm G}}=\nicefrac{5}{2^{7}}.\label{eq:degG}
\end{equation}

Let $q_{ij}(\overline{x})\,(1\leq i<j\leq5)$ be quadratic forms with
respect to $\overline{x}:=(\overline{x}_{1},\dots,\overline{x}_{n})$
satisfying the five Pl\"ucker relations. These $q_{ij}(\overline{x})$
will be chosen suitably for Types $R$ and $I\!R$. Then $\overline{{\rm G}}$
contains 
\[
\Pi:=\{[\overline{x},q(\overline{x})]\mid\overline{x}\in V\}\simeq\mP(W)\simeq\mP^{n-1},
\]
where we set $q(\overline{x})=\left(q_{ij}(\overline{x})\,(1\leq i<j\leq5)\right)$.
We denote by $q\colon\mP(W)\dashrightarrow{\rm G}(2,V)$ the quadratic
map defined by them.

Now we perform the coordinate change of $\mP(1^{n},2^{10})$ by 
\[
\overline{r}_{ij}:=\overline{p}_{ij}-q_{ij}(\overline{x})\,(1\leq i<j\leq5)
\]
with $\overline{x}$ being unchanged. By the new coordinates {[}$\overline{x}$,
$\overline{r}]$ with $\overline{r}:=(\overline{r}_{ij})\,(1\leq i<j\leq5)$,
$\overline{{\rm G}}$ is defined by the five Pl\"ucker relations
of the skew-symmetric matrix

{\footnotesize
\[
\left(\begin{array}{ccccc}
0 & \overline{r}_{12}+q_{12}(\overline{x}) & \overline{r}_{13}+q_{13}(\overline{x}) & \overline{r}_{14}+q_{14}(\overline{x}) & \overline{r}_{15}+q_{15}(\overline{x})\\
 & 0 & \overline{r}_{23}+q_{23}(\overline{x}) & \overline{r}_{24}+q_{24}(\overline{x}) & \overline{r}_{25}+q_{25}(\overline{x})\\
 &  & 0 & \overline{r}_{34}+q_{34}(\overline{x}) & \overline{r}_{35}+q_{35}(\overline{x})\\
 &  &  & 0 & \overline{r}_{45}+q_{45}(\overline{x})\\
 &  &  &  & 0
\end{array}\right),
\]
}which are expanded as 

\begin{align}
 & (\overline{r}_{ij}q_{kl}(\overline{x})+\overline{r}_{kl}q_{ij}(\overline{x}))-(\overline{r}_{ik}q_{jl}(\overline{x})+\overline{r}_{jl}q_{ik}(\overline{x}))+(\overline{r}_{il}q_{jk}(\overline{x})+\overline{r}_{jk}q_{il}(\overline{x}))+{\rm Pf}_{ijkl}(\overline{r})=0\label{eq:neweq}\\
 & \text{for}\ (i,j,k,l)=(1,2,3,4),(1,2,3,5),(1,2,4,5),(1,3,4,5),(2,3,4,5)\nonumber 
\end{align}
 (recall that we are assuming $q_{ij}(\overline{x})$ satisfy the
five Pl\"ucker relations). We also note that 
\[
\Pi=\{\overline{r}=0\},
\]
and this is disjoint from $\nicefrac{1}{2}$-singularities, and intersects
$\mP(W)$ along $\{\overline{r}=0,q_{ij}(\overline{x})=0\,(1\leq i<j\leq5)\}$. 

It turns out that the equation (\ref{eq:neweq}) has the following
meaning: 
\begin{prop}
\label{prop:bundleE}Let $\sE$ be the vector bundle on ${\rm {\rm G}(2,V)}$
of rank $7$ defined by 
\[
0\to\sE\to\wedge^{2}V\otimes\sO_{{\rm G}(2,V)}\to\wedge^{2}\sQ\to0,
\]
where $\sQ$ is the universal quotient bundle of rank $3$ on ${\rm G}(2,V)$.
\end{prop}

\begin{lem}
\begin{enumerate}

\item It holds that $\mP(\sE)$ is defined in $\mP(\wedge^{2}V)\times{\rm {\rm {\rm G}(2,V)}}$
by the equations
\[
(r_{ij}q_{kl}+r_{kl}q_{ij})-(r_{ik}q_{jl}+r_{jl}q_{ik})+(r_{il}q_{jk}+r_{jk}q_{il})=0
\]
for\textup{ $(i,j,k,l)=(1,2,3,4),(1,2,3,5),(1,2,4,5),(1,3,4,5),(2,3,4,5),$
where $r_{ij},q_{ij}\,(1\leq i<j\leq5)$ are the Pl\"ucker coordinates
of the first and second factors of }$\mP(\wedge^{2}V)\times{\rm {\rm G}(2,V)}$
respectively.

\item Let $\ell$ be a line in $\mP(V)$. The fiber of the first
projection $\mP(\sE)\to\mP(\wedge^{2}V)$ over $\ell\in{\rm G}(2,V)$
parameterizes lines in $\mP(V)$ intersecting $\ell$.

\end{enumerate}
\end{lem}

\begin{proof}
The assertion (1) follows by a direct computation on the charts of
${\rm G}(2,V)$. Using the ${\rm GL}(V)$-action on $\mP(\sE)$, the
assertion (2) can be checked explicitly for a special line $\ell$. 
\end{proof}
This lemma is effectively employed in the subsequent arguments.

\subsection{The rational map $\mu\colon\mP(\wedge^{2}V)\protect\dashrightarrow\mP(V^{*})$\label{subsec:PwedgeVtoPV*}}

For a nonzero element $\widetilde{r}:=(\widetilde{r}_{ij})\in\wedge^{2}V$,
we say the point $[\widetilde{r}]\in\mP(\wedge^{2}V)$ is of rank
$k$ if the corresponding skew-symmetric matrix $A(\widetilde{r})$
is of rank $k$ (cf.~Notation and Conventions in the end of Section
\ref{sec:Introduction}). Since $\dim V=5$, $\rank\widetilde{r}$
is at most $4$. Assume that $\widetilde{r}$ is of rank $4$. Then
there exist linearly independent 4 vectors $\bm{v}_{1},\bm{v}_{2},\bm{v}_{3},\bm{v}_{4}\in V$
such that $\bm{v}_{1}\wedge\bm{v}_{2}+\bm{v}_{3}\wedge\bm{v}_{4}=\widetilde{r}$.
Since $\widetilde{r}\wedge\widetilde{r}=2(\bm{v}_{1}\wedge\bm{v}_{2}\wedge\bm{v}_{3}\wedge\bm{v}_{4})$,
$V_{4}:=\langle\bm{v}_{1},\bm{v}_{2},\bm{v}_{3},\bm{v}_{4}\rangle$
coincides with the kernel of the linear map $V\ni t\mapsto t\wedge\widetilde{r}\wedge\widetilde{r}\in\wedge^{5}V\simeq\mC$.
In particular, $V_{4}$ is uniquely determined from $\widetilde{r}$. 

Subsequently, we need the rational map 
\[
\mu\colon\mP(\wedge^{2}V)\dashrightarrow\mP(V^{*})
\]
defined by $[\widetilde{r}]\mapsto[V_{4}\subset V]$ on the rank $4$
locus. By the consideration as above, we see that the map $\mu$ coincides
with the map $\mu'\colon\mP(\wedge^{2}V)\dashrightarrow\mP(\wedge^{4}V)\simeq\mP(V^{*})$
defined as the composite
\[
[\widetilde{r}]\mapsto[\widetilde{r}\wedge\widetilde{r}]\mapsto[V\ni t\mapsto t\wedge\widetilde{r}\wedge\widetilde{r}\in\wedge^{5}V\simeq\mC].
\]
Let $\mathsf{z}_{i}\,(1\leq i\leq5)$ be the dual coordinates of $V^{*}$.
Using the definition of $\mu'$, we may verify that the map $\mu=\mu'$
is defined by

\begin{equation}
[\widetilde{r}]\mapsto[\mathsf{z}]=[{\rm Pf}_{2345}(\widetilde{r}),-{\rm Pf}_{1345}(\widetilde{r}),{\rm Pf}_{1245}(\widetilde{r}),-{\rm Pf}_{1235}(\widetilde{r}),{\rm Pf}_{1234}(\widetilde{r})].\label{eq:zPf}
\end{equation}

\subsection{The rational map $\rho\colon\overline{{\rm G}}\protect\dashrightarrow\mP(\wedge^{2}V)$\label{subsec:GtoPwedgeV}}

In subsequent discussions, it is convenient to clarify the meaning
of the map 
\[
\rho\colon\overline{{\rm G}}\dashrightarrow\mP(\wedge^{2}V)\quad[\overline{x},\overline{r}]\mapsto[\overline{r}].
\]
Note that, for a general point $[\overline{x},\overline{r}]\in\overline{{\rm G}}$,
the rank of $\overline{r}=\overline{p}-q(\overline{x})$ is $4$.
In this case, $\overline{r}$ corresponds to $\bm{v}_{1}\wedge\bm{v}_{2}+\bm{v}_{3}\wedge\bm{v}_{4}$
for some linearly independent vectors $\bm{v_{1}},\dots,\bm{v}_{4}\in V$.
Setting $V_{4}:=\langle\bm{v}_{1},\dots,\bm{v}_{4}\rangle$, we have
$\rho([\overline{x},\overline{r}])=[\overline{r}]\in\mP(\wedge^{2}V_{4})$.
Conversely, when we start from a point $[\widetilde{r}]\in\mP(\wedge^{2}V_{4})$
of rank $4$ for some $4$-dimensional vector space $V_{4}\subset V$,
the $\rho$-fiber over the point $[\widetilde{r}]$ consists of points
$[\overline{x},\overline{p}-q(\overline{x})]\in\overline{{\rm G}}$
such that $\overline{p}-q(\overline{x})=s\widetilde{r}$ for some
$s\not=0$. Since both $\overline{p}$ and $q(\overline{x})$ are
of rank $2$, we see that $[\overline{p}],[q(\overline{x})]\in{\rm G}(2,V_{4})$.
Since ${\rm G}(2,V_{4})\subset\mP(\wedge^{2}V_{4})$ is a quadratic
hypersurface, $\overline{p}$ is uniquely determined by the condition
$\overline{p}=q(\overline{x})+s\widetilde{r}$ once we fix $\overline{x}$
and $\widetilde{r}$. Therefore the $\rho$-fiber over the point $\widetilde{r}$
is birational to $q^{-1}({\rm G}(2,V_{4}))\subset\mP(W)$.

Subsequently we will determine $q^{-1}({\rm G}(2,V_{4}))$ in the
case of Types $R$ and $I\!R$, respectively.

\section{\textbf{The key variety $\mathscr{R}$ for Type $R$ $\mQ$-Fano
$3$-folds \label{sec:TypeR}}}

In Subsection \ref{subsec:The-cone GR}, we specialize the construction
of Section \ref{sec:ConeG(2,5)} to the setting of Type $R$ $\mQ$-Fano
3-folds. In Subsection \ref{subsec:Constructing-the-keyR}, we construct
the key variety $\mathscr{R}$ for a Type $R$ $\mQ$-Fano $3$-fold
$X$, featuring the following $\mathfrak{S}_{5}$-equivariant Sarkisov
link:

\begin{equation}\label{eq:Sarkisov} \xymatrix{& \widehat{\mathscr{R}}\ar@{-->}[rr]\ar[dl]_{\widehat{f}_{\mathscr{R}}}\ar[dr]^{\widehat{g}_{\mathscr{R}}} &  &\widetilde{\mathscr{R}}\ar[dr]^{\widetilde{f}_{\mathscr{R}}}\ar[dl]_{\widetilde{g}_{\mathscr{R}}}\\
\mathscr{R}& &\overline{\mathscr{R}}&  & L_\mathscr{R}.}
\end{equation}The construction of $\mathscr{R}$ consists of several steps and proceeds
with the aid of minimal model theory. We start by constructing explicitly
the midpoint $\overline{\mathscr{R}}$ of this Sarkisov link as a
linear section of $\overline{{\rm G}}_{\mathscr{R}}$ and confirm
that $\overline{\mathscr{R}}$ carries an $\mathfrak{S}_{5}$-action
(Subsections \ref{subsec:The-midpoint} and \ref{subsec:S5-action}).
Next, we identify a $3$-space $L_{\mathscr{R}}\subset\mP(\wedge^{2}V)$
compatible with $\overline{\mathscr{R}}$, and construct the quasi-$\mP^{7}$-bundle
$\widetilde{f}_{\mathscr{R}}\colon\widetilde{\mathscr{R}}\to L_{\mathscr{R}}$
(Subsections \ref{subsec:LR and S3} and \ref{subsec:The-quasi--bundle R}).
We study basic properties of this morphism in Subsection \ref{subsec:Basic-properties-of}.
In Subsection \ref{subsec:RTildePOmega}, we establish the existence
of a birational map \textbf{$\widetilde{{\rm \mathscr{R}}}\dashrightarrow\mathsf{R}=\mP(\Omega_{\mP(V^{*})}(1)^{\oplus2}|_{{\rm S}_{3}})$}.
We then construct the birational morphism $\widetilde{g}_{\mathscr{R}}\colon\widetilde{\mathscr{R}}\to\overline{\mathscr{R}}$,
which is shown to be a flopping contraction (Subsection \ref{subsec:The-flopping-contraction}).
After introducing necessary notation in Subsection \ref{subsec:Glossary-of-notation},
we construct the key variety $\mathscr{R}$ from $\overline{\mathscr{R}}$
via unprojection (Subsection \ref{subsec:Constructing--via unproj}).
We verify $\mathscr{R}$ admits the $\mathfrak{S}_{6}$-action extending
birationally the $\mathfrak{S}_{5}$-action on $\overline{\mathscr{R}}$
(Subsection \ref{subsec:S6}). Finally, we complete the Sarkisov diagram
(\ref{eq:Sarkisov}) decomposing the unprojection process $\overline{\mathscr{R}}\dashrightarrow\mathscr{R}$
into the composite of the inverse of the flopped contraction $\widehat{g}_{\mathscr{R}}\colon\widehat{\mathscr{R}}\to\overline{\mathscr{R}}$
for $\widetilde{g}_{\mathscr{R}}$ and the blow-up $\widehat{f}_{\mathscr{R}}\colon\widehat{\mathscr{R}}\to\mathscr{R}$
at the $\nicefrac{1}{2}$-singularity produced by the unprojection
(Subsections \ref{subsec:The-blow-up-} and \ref{subsec:The-flopped-contraction}). 

The construction is carried out while explicitly determining the defining
equations of the birational models of $\mathscr{R}$ at each intermediate
step. One reason for adopting this concrete approach is that the flop
$\widetilde{\mathscr{R}}\dashrightarrow\widehat{\mathscr{R}}$ in
the construction is complex and could not be handled without explicit
equations. One benefit of this concrete method is that it allowed
us to clarify that $\mathscr{R}$ admits an $\mathfrak{S}_{6}$-action.
Furthermore, although not discussed in this paper in detail, we are
also able to elucidate the relationship between the key variety $\mathscr{R}$
and a $\mathbb{P}^{2}\times\mathbb{P}^{2}$-fibration (Remark \ref{rem:P2P2Fib}).

\subsection{The cone $\overline{{\rm G}}_{\mathscr{R}}$ over ${\rm G}(2,5)$\label{subsec:The-cone GR}}

\subsubsection{\textbf{The setting}}

In the case of Type $R$, we consider the cone $\overline{{\rm G}}$
over ${\rm G}(2,5)$ as in Section \ref{sec:ConeG(2,5)} under the
situation $W=V\oplus V$. We denote $\overline{{\rm G}}$ by $\overline{{\rm G}}_{\mathscr{R}}$.
We denote by $\overline{x}=(\overline{x}_{1},\dots,\overline{x}_{5})$
and $\overline{y}=(\overline{y}_{1},\dots,\overline{y}_{5})$ the
coordinates of the first and the second factors of $W$ respectively.
We also denote by $V_{x}$ and $V_{y}$ the first and the second factors
of $W$ respectively. We set

\begin{equation}
q_{ij}(x,y):=\left|\begin{array}{cc}
x_{i} & x_{j}\\
y_{i} & y_{j}
\end{array}\right|\,(1\leq i<j\leq5).\label{eq:qij}
\end{equation}
As the quadrics $q_{ij}(\overline{x})$ in Section \ref{sec:ConeG(2,5)},
we take these $q_{ij}(\overline{x},\overline{y})$. Note that these
satisfy the five Pl\"ucker relations.

We consider the standard action of ${\rm {\rm GL}}(V)$ on $V$. Then
the actions of ${\rm {\rm GL}}(V)$ on $\mP(\wedge^{2}V)$, $\mP(V^{*}),$
$\overline{{\rm G}}_{\mathscr{R}}$, etc, are naturally induced.

\subsubsection{\textbf{A rational map $\overline{{\rm G}}_{\mathscr{R}}\protect\dashrightarrow\mP(\Omega_{\mP(V^{*})}(1)^{\oplus2})$}}

Note that the projective bundle $\mP(\Omega_{\mP(V^{*})}(1)^{\oplus2})$
is defined in $\mP(V^{\oplus2})\times\mP(V^{*})$ by 
\begin{equation}
\sum_{i=1}^{5}\mathsf{x}_{i}\mathsf{z}_{i}=0,\sum_{i=1}^{5}\mathsf{y}_{i}\mathsf{z}_{i}=0,\label{eq:P(Ogema)}
\end{equation}
where $\mathsf{x}_{i}$, $\mathsf{\ensuremath{y}}_{i}\,(1\leq i\leq5)$
are the coordinates of the first and the second factors of $V^{\oplus2}$
respectively, and $\mathsf{z}_{i}\,(1\leq i\leq5)$ the dual coordinates
of $V^{*}$. Since $q^{-1}({\rm G}(2,V_{4}))=\mP(V_{4}^{\oplus2})$
for a $4$-dimensional vector space $V_{4}\subset V$, the consideration
as in Subsection \ref{subsec:GtoPwedgeV} suggests that the rational
map $\rho\colon\overline{G}_{\mathscr{R}}\dashrightarrow\mP(V^{*})$
induces a rational map \textbf{$\overline{{\rm G}}_{\mathscr{R}}\dashrightarrow\mP(\Omega_{\mP(V^{*})}(1)^{\oplus2})$}.
We can check the following by a direct calculation:
\begin{lem}
\label{lem:TwoRel}The points of $\overline{{\rm G}}_{\mathscr{R}}$
satisfy the following two relations: 

\begin{equation}
\begin{cases}
\overline{x}_{1}{\rm Pf}_{2345}(\overline{r})-\overline{x}_{2}{\rm Pf}_{1345}(\overline{r})+\overline{x}_{3}{\rm Pf}_{1245}(\overline{r})-\overline{x}_{4}{\rm Pf}_{1235}(\overline{r})+\overline{x}_{5}{\rm Pf}_{1234}(\overline{r})=0,\\
\overline{y}_{1}{\rm Pf}_{2345}(\overline{r})-\overline{y}_{2}{\rm Pf}_{1345}(\overline{r})+\overline{y}_{3}{\rm Pf}_{1245}(\overline{r})-\overline{y}_{4}{\rm Pf}_{1235}(\overline{r})+\overline{y}_{5}{\rm Pf}_{1234}(\overline{r})=0.
\end{cases}\label{eq:GR}
\end{equation}
\end{lem}

Therefore the composite $\overline{{\rm G}}_{\mathscr{R}}\overset{\rho}{\dashrightarrow}\mP(\wedge^{2}V)\overset{\mu}{\dashrightarrow}\mP(V^{*})$
(cf.~Subsections \ref{subsec:PwedgeVtoPV*} and \ref{subsec:GtoPwedgeV})
certainly induces the rational map \textbf{$\overline{{\rm G}}_{\mathscr{R}}\dashrightarrow\mP(\Omega_{\mP(V^{*})}(1)^{\oplus2})$}
defined by 
\[
[\overline{x},\overline{y},\overline{r}]\mapsto[\overline{x},\overline{y}]\times[{\rm Pf}_{2345}(\overline{r}),-{\rm Pf}_{1345}(\overline{r}),{\rm Pf}_{1245}(\overline{r}),-{\rm Pf}_{1235}(\overline{r}),{\rm Pf}_{1234}(\overline{r})].
\]

\subsection{Constructing the key variety $\mathscr{R}$\label{subsec:Constructing-the-keyR}}

\subsubsection{\textbf{The midpoint $\overline{\mathscr{R}}$\label{subsec:The-midpoint}}}
\begin{defn}
We define 
\begin{equation}
\overline{\mathscr{R}}:=\overline{{\rm G}}_{\mathscr{R}}\cap\{\overline{r}_{12}=\overline{r}_{14}=\overline{r}_{45}=0,\,\overline{r}_{13}=\overline{r}_{34}+\overline{r}_{35},\,\overline{r}_{23}=-\overline{r}_{34}-\overline{r}_{35},\,\overline{r}_{25}=-\overline{r}_{35}\}.\label{eq:linear overlineR}
\end{equation}
\end{defn}

The scheme $\overline{\mathscr{R}}$ is chosen such that it has an
$\mathfrak{S}_{5}$-action (see Subsection \ref{subsec:S5-action}).

We consider $\overline{\mathscr{R}}$ as a subvariety of $\mP(1^{10},2^{4})$
with the coordinates 
\begin{equation}
\overline{x},\overline{y},\overline{r}_{15},\overline{r}_{24},\overline{r}_{34},\overline{r}_{35}.\label{eq:coord overline}
\end{equation}
 Let us write down the equations of $\overline{\mathscr{R}}$ explicitly,
which are obtained from the equations of $\overline{{\rm G}}_{\mathscr{R}}$.
Setting 
\begin{equation}
r_{12}=r_{14}=r_{45}=0,\,r_{13}=r_{34}+r_{35},\,r_{23}=-r_{34}-r_{35},\,r_{25}=-r_{35}\label{eq:S5r}
\end{equation}
for the 5 Pl\"ucker quadrics ${\rm Pf}_{1234}(r),\dots,{\rm Pf}_{2345}(r)$,
we obtain

\begin{align}
f_{1}(r) & :=-r_{24}(r_{34}+r_{35}),\,f_{2}(r):=-(r_{15}-r_{35})(r_{34}+r_{35}),\label{eq:f1--f5}\\
f_{3} & (r):=r_{15}r_{24},\,f_{4}(r):=r_{15}r_{34},f_{5}(r):=-(r_{24}+r_{34})r_{35},\nonumber 
\end{align}
where we set $r=(r_{15},r_{24},r_{34},r_{35})$ for (\ref{eq:f1--f5}).
Setting $q_{ij}:=q_{ij}(x,y)$, we define the following five polynomials:
\begin{equation}
\begin{cases}
RF_{1}(x,y,r):=-q_{13}r_{24}+\left(q_{12}-q_{14}-q_{24}\right)r_{34}-\left(q_{14}+q_{24}\right)r_{35}+f_{1}(r)=0,\\
RF_{2}(x,y,r):=q_{23}r_{15}-\left(q_{15}+q_{25}\right)r_{34}+\left(q_{12}+q_{13}-q_{15}-q_{25}\right)r_{35}+f_{2}(r)=0,\\
RF_{3}(x,y,r):=q_{24}r_{15}+q_{15}r_{24}+q_{14}r_{35}+f_{3}(r)=0,\\
RF_{4}(x,y,r):=q_{34}r_{15}+\left(q_{15}+q_{45}\right)r_{34}-\left(q_{14}-q_{45}\right)r_{35}+f_{4}(r)=0,\\
RF_{5}(x,y,r):=-q_{35}r_{24}+\left(q_{25}-q_{45}\right)r_{34}-\left(q_{24}+q_{34}+q_{45}\right)r_{35}+f_{5}(r)=0,
\end{cases}\label{eq:R}
\end{equation}
which are the five Pl\"ucker relations for the skew-symmetric matrix

{\footnotesize
\[
\left(\begin{array}{ccccc}
0 & q_{12} & r_{34}+r_{35}+q_{13} & q_{14} & r_{15}+q_{15}\\
 & 0 & -r_{34}-r_{35}+q_{23} & r_{24}+q_{24} & -r_{35}+q_{25}\\
 &  & 0 & r_{34}+q_{34} & r_{35}+q_{35}\\
 &  &  & 0 & q_{45}\\
 &  &  &  & 0
\end{array}\right).
\]
}Therefore we have the following:
\begin{prop}
It holds that 
\begin{equation}
\overline{\mathscr{R}}=\{RF_{1}(\overline{x},\overline{y},\overline{r})=RF_{2}(\overline{x},\overline{y},\overline{r})=RF_{3}(\overline{x},\overline{y},\overline{r})=RF_{4}(\overline{x},\overline{y},\overline{r})=RF_{5}(\overline{x},\overline{y},\overline{r})=0\}\label{eq:MidEq}
\end{equation}
 in $\mP(1^{10},2^{4}).$
\end{prop}

By Lemma \ref{lem:TwoRel}, we have immediately the following:
\begin{lem}
\label{lem:TwoRelonRvar}The points of $\overline{{\rm \mathscr{R}}}$
satisfy the following two relations: 

\begin{equation}
\begin{cases}
\overline{x}_{1}f_{5}(\overline{r})-\overline{x}_{2}f_{4}(\overline{r})+\overline{x}_{3}f_{3}(\overline{r})-\overline{x}_{4}f_{2}(\overline{r})+\overline{x}_{5}f_{1}(\overline{r})=0,\\
\overline{y}_{1}f_{5}(\overline{r})-\overline{y}_{2}f_{4}(\overline{r})+\overline{y}_{3}f_{3}(\overline{r})-\overline{y}_{4}f_{2}(\overline{r})+\overline{y}_{5}f_{1}(\overline{r})=0.
\end{cases}\label{eq:Rvar2rel}
\end{equation}
\end{lem}

\subsubsection{\textbf{$\mathfrak{S}_{5}$-action on $\overline{\mathscr{R}}$\label{subsec:S5-action}}}

We set $\sigma_{i}:=(i\,i+1)\in\mathfrak{S}_{5}\,(1\leq i\leq4)$.
By associating the following matrices to $\sigma_{i}$, the $4$-dimensional
linear space with coordinates $\overline{r}_{15},\overline{r}_{24},\overline{r}_{34},\overline{r}_{35}$
becomes an irreducible $\mathfrak{S}_{5}$-representation which is
isomorphic to the tensor product between the standard and the sign
representations:

\begin{align}
\sigma_{1}:{\footnotesize \left(\begin{array}{cccc}
-1 & 0 & 0 & 0\\
0 & 1 & 0 & 0\\
0 & -1 & -1 & 0\\
0 & 1 & 0 & -1
\end{array}\right)},\, & \sigma_{2}:{\footnotesize \left(\begin{array}{cccc}
0 & 0 & 0 & -1\\
-1 & -1 & 0 & 1\\
1 & 0 & -1 & -1\\
-1 & 0 & 0 & 0
\end{array}\right)},\label{eq:S5 4dim}\\
\sigma_{3}:{\footnotesize \left(\begin{array}{cccc}
0 & 0 & -1 & -1\\
0 & -1 & 0 & 0\\
-1 & 0 & 0 & 1\\
0 & 0 & 0 & -1
\end{array}\right)},\, & \sigma_{4}:{\footnotesize \left(\begin{array}{cccc}
-1 & 0 & 1 & 1\\
0 & -1 & 0 & 0\\
0 & 0 & 0 & 1\\
0 & 0 & 1 & 0
\end{array}\right)}.\nonumber 
\end{align}

Moreover, by associating the following matrices to $\sigma_{i}$,
$V_{x}$ and $V_{y}$ become $5$-dimensional irreducible representations:
\begin{align}
\sigma_{1}:{\footnotesize \left(\begin{array}{ccccc}
1 & 0 & 0 & 0 & 0\\
0 & 1 & 0 & 0 & 0\\
0 & -1 & -1 & 0 & 0\\
0 & 0 & 0 & 1 & 0\\
-1 & 0 & 0 & -1 & -1
\end{array}\right)},\, & \sigma_{2}:{\footnotesize \left(\begin{array}{ccccc}
0 & 0 & -1 & 0 & 0\\
1 & 1 & 1 & 0 & 0\\
-1 & 0 & 0 & 0 & 0\\
0 & 1 & 1 & -1 & -1\\
0 & 0 & 0 & 0 & 1
\end{array}\right)},\label{eq:standard S5}\\
\sigma_{3}:{\footnotesize \left(\begin{array}{ccccc}
1 & 0 & 0 & 0 & 0\\
-1 & 0 & 0 & -1 & 0\\
0 & 0 & 0 & 0 & -1\\
-1 & -1 & 0 & 0 & 0\\
0 & 0 & -1 & 0 & 0
\end{array}\right)},\, & \sigma_{4}:{\footnotesize \left(\begin{array}{ccccc}
1 & 0 & 0 & 0 & 0\\
0 & 1 & 0 & 0 & 0\\
0 & 0 & 1 & 0 & 0\\
0 & 1 & 0 & -1 & 0\\
-1 & 0 & 1 & 0 & -1
\end{array}\right).}\nonumber 
\end{align}
These $\mathfrak{S}_{5}$-representations induce an action of $\mathfrak{S}_{5}$
on $\mP(1^{10},2^{4})$. We may directly check the following: 
\begin{prop}
The $\mathfrak{S}_{5}$-action on $\mP(1^{10},2^{4})$ induces that
on $\overline{\mathscr{R}}$.
\end{prop}

\subsubsection{\textbf{The $3$-space $L_{\mathscr{R}}$ in $\mP(\wedge^{2}V)$
and the Segre cubic ${\rm S}_{3}$\label{subsec:LR and S3}}}

Let us start by a slightly general consideration.

Let $L_{\mathscr{R}}$ be a $3$-space of $\mP(\wedge^{2}V)$ such
that ${\rm G}(2,V)\cap L_{\mathscr{R}}$ consists of 5 points, which
we denote by $\mathsf{q}_{1},\dots,\mathsf{q}_{5}$. Since $\deg{\rm G}(2,V)=5$,
this certainly holds for a general $L_{\mathscr{R}}$. 
\begin{lem}
The $5$ points $\mathsf{q}_{1},\dots,\mathsf{q}_{5}$ are in  general
position.
\end{lem}

\begin{proof}
The set of $\mathsf{q}_{1},\dots,\mathsf{q}_{5}$ is Gorenstein since
so is ${\rm G}(2,V)$. Therefore, by \cite[Ex.2.1]{Mi}, the 5 points
are in  general position.
\end{proof}
Since $\mu\colon\mP(\wedge^{2}V)\dashrightarrow\mP(V^{*})$ as in
Subsection \ref{subsec:PwedgeVtoPV*} is defined by the 5 Pl\"ucker
quadrics defining ${\rm G}(2,V)$, its restriction $L_{\mathscr{R}}\dashrightarrow\mP(V^{*})$
is defined by the 5 quadrics defining the 5 points. Therefore this
situation is compatible with the one in Subsection \ref{subsec:Segre},
so we denote by ${\rm S}_{3}$ the image of $L_{\mathscr{R}}\dashrightarrow\mP(V^{*})$
and by $\mu_{{\rm S}}\colon L_{\mathscr{R}}\dashrightarrow{\rm S}_{3}$
the induced map.

\vspace{5pt}

Now we specify $L_{\mathscr{R}}$ as follows:
\[
L_{\mathscr{R}}:=\{\widetilde{r}_{12}=\widetilde{r}_{14}=\widetilde{r}_{45}=0,\,\widetilde{r}_{13}=\widetilde{r}_{34}+\widetilde{r}_{35},\,\widetilde{r}_{23}=-\widetilde{r}_{34}-\widetilde{r}_{35},\,\widetilde{r}_{25}=-\widetilde{r}_{35}\},
\]
which certainly cut out 5 points from ${\rm G}(2,V)$, and is compatible
with (\ref{eq:S5r}). We identify $L_{\mathscr{R}}$ with $\mP^{3}$
choosing $\widetilde{r}_{15},\widetilde{r}_{24},\widetilde{r}_{34},\widetilde{r}_{35}$
as coordinates of $L_{\mathscr{R}}$. Then the five points in ${\rm G}(2,V)\cap L_{\mathscr{R}}$
are the following:
\begin{equation}
\begin{cases}
\mathsf{q_{1}:=}\text{the}\ \ensuremath{\widetilde{r}_{15}}\text{-point},\mathsf{q_{2}:=}\text{the}\ \ensuremath{\widetilde{r}_{24}}\text{-point},\mathsf{q_{3}:=}\text{the}\ \ensuremath{\widetilde{r}_{34}}\text{-point},\\
\text{\ensuremath{\mathsf{q_{4}:=}}the point with \ensuremath{\{\widetilde{r}_{24}=\widetilde{r}_{34}=0,\widetilde{r}_{15}=\widetilde{r}_{35}\}},}\\
\text{\ensuremath{\mathsf{q_{5}:=}}the point with \ensuremath{\{\widetilde{r}_{15}=0,\widetilde{r}_{24}=-\widetilde{r}_{34}=\widetilde{r}_{35}\}}.}
\end{cases}\label{eq:points qi}
\end{equation}

We may check that the matrices (\ref{eq:S5 4dim}) induce an $\mathfrak{S}_{5}$-action
on $L_{\mathscr{R}}\simeq\mP^{3}$ with coordinates $\widetilde{r}_{15},\widetilde{r}_{24},\widetilde{r}_{34},\widetilde{r}_{35}$
symmetrically acting on the five points.

The restriction $L_{\mathscr{R}}\dashrightarrow\mP(V^{*})$ of $\mu\colon\mP(\wedge^{2}V)\dashrightarrow\mP(V^{*})$
is defined by 
\begin{equation}
[\widetilde{r}_{15},\widetilde{r}_{24},\widetilde{r}_{34},\widetilde{r}_{35}]\mapsto[f_{5}(\widetilde{r}),-f_{4}(\widetilde{r}),f_{3}(\widetilde{r}),-f_{2}(\widetilde{r}),f_{1}(\widetilde{r})]\label{eq:mu coord}
\end{equation}
and the image ${\rm S}_{3}\subset\mP(V^{*})$ is defined by
\begin{equation}
-\mathsf{z}_{1}(\mathsf{z}_{3}\mathsf{z}_{4}-\mathsf{z}_{2}\mathsf{z}_{5})=(\mathsf{z}_{2}-\mathsf{z}_{3})(\mathsf{z}_{2}+\mathsf{z}_{4})\mathsf{z}_{5}.\label{eq:S3}
\end{equation}

\subsubsection{\textbf{The quasi-$\mP^{7}$-bundle }$\widetilde{f}_{\mathscr{R}}\colon\widetilde{{\rm \mathscr{R}}}\to L_{\mathscr{R}}$\label{subsec:The-quasi--bundle R}}

By Lemma \ref{lem:TwoRelonRvar} and the equations (\ref{eq:R}) and
(\ref{eq:MidEq}) of $\overline{\mathscr{R}}$, we are led to define
a variety $\widetilde{{\rm \mathscr{R}}}$ over $L_{\mathscr{R}}$
as follows, which play a role connecting $\overline{\mathscr{R}}$
and $\mathsf{R}=\mP(\Omega_{\mP(V^{*})}(1)^{\oplus2}|_{{\rm S}_{3}})$.

Let $V_{\mathscr{R}}\subset\wedge^{2}V$ be the $4$-dimensional subspace
such that $\mP(V_{\mathscr{R}})=L_{\mathscr{R}}$. Following \cite[Chap.2]{R2}
with a modification as in \cite[p.16, and A.1]{BCZ}, we define the
toric variety $\mF_{\widetilde{\mathscr{R}}}$ as the quotient of
$(V_{\mathscr{R}}\setminus\{\bm{o}\})\times\left((V^{\oplus2}\oplus\mC)\setminus\{\bm{o}\}\right)$
by the action of $(\mC^{*})^{2}$ with the weights given as in Table
\ref{Table:FR}, where the entries of the second line are the coordinates
of $((V_{\mathscr{R}}\setminus\{\bm{o}\})\times\left((V^{\oplus2}\oplus\mC)\setminus\{\bm{o}\}\right)$,
those of the third and the fourth lines are the weights of the coordinates
by the action of the first and the second factors of $(\mC^{*})^{2}$
respectively.
\begin{table}[h]
\begin{tabular}{|c|c|c|c|c|}
\hline 
spaces & $V_{\mathscr{R}}\setminus\{\bm{o}\}$ & \multicolumn{3}{c|}{$(V^{\oplus2}\oplus\mC)\setminus\{\bm{o}\}$}\tabularnewline
\hline 
coordinates & $\widetilde{r}_{15},\widetilde{r}_{24},\widetilde{r}_{34},\widetilde{r}_{35}$ & $\widetilde{x}$ & $\widetilde{y}$ & $\widetilde{w}$\tabularnewline
\hline 
1st weights & $1$ & $0$ & $0$ & $-1$\tabularnewline
\hline 
2nd weights & $0$ & $1$ & $1$ & $2$\tabularnewline
\hline 
\end{tabular}

\caption{Toric variety $\mF_{\widetilde{\mathscr{R}}}$ }
\label{Table:FR}
\end{table}
 The quotient of $V_{\mathscr{R}}\setminus\{\bm{o}\}$ by the action
of the first factor of $(\mC^{*})^{2}$ is nothing but $L_{\mathscr{R}}=\mP(V_{\mathscr{R}})$.
Hence we have the natural projection $\mF_{\widetilde{\mathscr{R}}}\to L_{\mathscr{R}}$.
This is a $\mP(1^{10},2)$-bundle, i.e., the variety which is locally
the product of an open subset of $L_{\mathscr{R}}$ and $\mP(1^{10},2)$.
Indeed, by Table \ref{Table:FR}, we see that $\mF_{\widetilde{\mathscr{R}}}|_{\widetilde{r}_{ij}=1}\simeq\left((\wedge^{2}V)|_{\widetilde{r}_{ij}=1}\right)\times\mP(1^{10},2)$.
\begin{defn}
We set $\widetilde{q}_{ij}:=q_{ij}(\widetilde{x},\widetilde{y})\,(1\leq i<j\leq5)$.
In the toric variety $\mF_{\widetilde{\mathscr{R}}}$, we define $\widetilde{\mathscr{R}}$
by the following equations:
\begin{equation}
\begin{cases}
\widetilde{x}_{1}f_{5}(\widetilde{r})-\widetilde{x}_{2}f_{4}(\widetilde{r})+\widetilde{x}_{3}f_{3}(\widetilde{r})-\widetilde{x}_{4}f_{2}(\widetilde{r})+\widetilde{x}_{5}f_{1}(\widetilde{r})=0,\\
\widetilde{y}_{1}f_{5}(\widetilde{r})-\widetilde{y}_{2}f_{4}(\widetilde{r})+\widetilde{y}_{3}f_{3}(\widetilde{r})-\widetilde{y}_{4}f_{2}(\widetilde{r})+\widetilde{y}_{5}f_{1}(\widetilde{r})=0,\\
-\widetilde{q}_{13}\widetilde{r}_{24}+\left(\widetilde{q}_{12}-\widetilde{q}_{14}-\widetilde{q}_{24}\right)\widetilde{r}_{34}-\left(\widetilde{q}_{14}+\widetilde{q}_{24}\right)\widetilde{r}_{35}+\widetilde{w}f_{1}(\widetilde{r})=0,\\
\widetilde{q}_{23}\widetilde{r}_{15}-\left(\widetilde{q}_{15}+\widetilde{q}_{25}\right)\widetilde{r}_{34}+\left(\widetilde{q}_{12}+\widetilde{q}_{13}-\widetilde{q}_{15}-\widetilde{q}_{25}\right)\widetilde{r}_{35}+\widetilde{w}f_{2}(\widetilde{r})=0,\\
\widetilde{q}_{24}\widetilde{r}_{15}+\widetilde{q}_{15}\widetilde{r}_{24}+\widetilde{q}_{14}\widetilde{r}_{35}+\widetilde{w}f_{3}(\widetilde{r})=0,\\
\widetilde{q}_{34}\widetilde{r}_{15}+\left(\widetilde{q}_{15}+\widetilde{q}_{45}\right)\widetilde{r}_{34}-\left(\widetilde{q}_{14}-\widetilde{q}_{45}\right)\widetilde{r}_{35}+\widetilde{w}f_{4}(\widetilde{r})=0,\\
-\widetilde{q}_{35}\widetilde{r}_{24}+\left(\widetilde{q}_{25}-\widetilde{q}_{45}\right)\widetilde{r}_{34}-\left(\widetilde{q}_{24}+\widetilde{q}_{34}+\widetilde{q}_{45}\right)\widetilde{r}_{35}+\widetilde{w}f_{5}(\widetilde{r})=0,
\end{cases}\label{eq:Rtilde}
\end{equation}
where the first two equations (resp.~the remaining equations) originate
from Lemma \ref{lem:TwoRelonRvar} (resp.~the equations $RF_{1}(\overline{x},\overline{y},\overline{r})$--$RF_{5}(\overline{x},\overline{y},\overline{r})$
of $\overline{\mathscr{R}}$). Since the polynomials in (\ref{eq:Rtilde})
are quasi-homogeneous with respect to the $(\mC^{*})^{2}$-action,
$\widetilde{\mathscr{R}}$ is well-defined as a subvariety of $\mF_{\widetilde{\mathscr{R}}}$.
We denote by $\widetilde{f}{}_{\mathscr{R}}\colon\widetilde{\mathscr{R}}\to L_{\mathscr{R}}$
the natural projection.
\end{defn}

We will prove that there exists a flopping contraction $\widetilde{\mathscr{R}}\to\overline{\mathscr{R}}$
(Proposition \ref{prop:hR}).

\subsubsection{\textbf{The $\mathfrak{S}_{5}$-action on $\widetilde{\mathscr{R}}$}}

Using the $\mathfrak{S}_{5}$-representations as in Subsection \ref{subsec:S5-action},
we may define an $\mathfrak{S}_{5}$-action on $\widetilde{\mathscr{R}}$
by setting $w$ as an $\mathfrak{S}_{5}$-invariant (compare the equations
(\ref{eq:MidEq}) of $\overline{\mathscr{R}}$ and the equations (\ref{eq:Rtilde})
of $\widetilde{\mathscr{R}}$. we may also directly check that the
first two polynomials in (\ref{eq:Rtilde}) are $\mathfrak{S}_{5}$-invariants). 

\subsubsection{\textbf{Basic properties of $\widetilde{f}{}_{\mathscr{R}}\colon\widetilde{\mathscr{R}}\to L_{\mathscr{R}}$\label{subsec:Basic-properties-of}}}
\begin{defn}
In the weighted projective space $\mP(1^{10},2)$ with weight one
coordinates $x_{1},\dots,x_{10}$ and weight two coordinate $w$,
we define the $8$-dimensional subvariety $\mathsf{C}$ as follows:
\[
\mathsf{C}:=\left\{ \rank\left(\begin{array}{ccc}
x_{1} & x_{2} & x_{3}\\
x_{4} & x_{5} & x_{6}
\end{array}\right)\leq1\right\} \subset\mP(x_{1},\dots,x_{10},w),
\]
which is the weighted cone over $\mP^{1}\times\mP^{2}$.
\end{defn}

\begin{prop}
\label{prop:f'} The morphism $\widetilde{f}{}_{\mathscr{R}}$ is
a $\mP^{7}$-bundle outside of $\cup_{i=1}^{5}\widetilde{f}{}_{\mathscr{R}}^{-1}(\mathsf{q}_{i})$.
The $\widetilde{f}{}_{\mathscr{R}}$-fiber over $\mathsf{q}_{i}\,(1\leq i\leq5)$
is isomorphic to the variety $\mathsf{C}.$ In particular the Picard
number of $\widetilde{\mathscr{R}}$ is $2$.
\end{prop}

\begin{proof}
Note that the open subset $L_{\mathscr{R}}\setminus\{\mathsf{q}_{1},\dots,\mathsf{q}_{5}\}$
is the union of open subsets $\{\widetilde{f}_{i}\not=0\}\,(1\leq i\leq5)$.
If $\widetilde{f}_{1}\not=0$, then we can solve the first three equations
of (\ref{eq:Rtilde}) with respect to the three coordinates $\widetilde{x}_{5},\widetilde{y}_{5}$
and $\widetilde{w}$, and the other equalities of (\ref{eq:Rtilde})
hold by substituting the solutions into them. This implies that, over
the open subset $\{\widetilde{f}_{1}\not=0$\}, $\widetilde{\mathscr{R}}$
is isomorphic to $\mP^{7}\times\{\widetilde{f}_{1}\not=0\}$, where
the coordinates of $\mP^{7}$ are $\widetilde{x}_{i},\widetilde{y}_{i}\,(1\leq i\le4)$.
We can also show that $\widetilde{f}{}_{\mathscr{R}}$ is a $\mP^{7}$-bundle
over the other open subsets $\{\widetilde{f}_{i}\not=0$\} using the
$\mathfrak{S}_{5}$-action.

We determine the $\widetilde{f}_{\mathscr{R}}$-fiber over $\mathsf{q}_{i}\,(1\leq i\leq5)$.
By the $\text{\ensuremath{\mathfrak{S}_{5}}}$-action on $\widetilde{\mathscr{R}}$,
we have only to consider the $\widetilde{f}_{\mathscr{R}}$-fiber
over the $\widetilde{r}_{15}$-point. We can directly see that it
is 
\[
\left\{ \rank\left(\begin{array}{ccc}
\widetilde{x}_{2} & \widetilde{x}_{3} & \widetilde{x}_{4}\\
\widetilde{y}_{2} & \widetilde{y}_{3} & \widetilde{y}_{4}
\end{array}\right)\leq1\right\} \subset\mP(\widetilde{x},\widetilde{y},\widetilde{w})=\mP(1^{10},2),
\]
which is isomorphic to $\mathsf{C}$. 

Since there are no nontrivial morphism from $\mathsf{C}$, the relative
Picard number of $\widetilde{f}{}_{\mathscr{R}}$ is $1$, hence it
holds that the Picard number of $\widetilde{\mathscr{R}}$ is $2$. 
\end{proof}
\begin{prop}
\label{prop:QfacR}Singularities of $\widetilde{{\rm \mathscr{R}}}$
are contained in the fibers over $\mathsf{q}_{i}\,(1\leq i\leq5)$.
In each fiber, $\widetilde{{\rm \mathscr{R}}}$ has a $\nicefrac{1}{2}$-singularity
at the $w$-point, which we call $\widetilde{\mathsf{p}}_{i}$, and
$c({\rm G}(2,5))$-singularities along $\mathbb{P}^{3}$ (we write
down explicitly the five copies $\widetilde{\Gamma}_{i}$ of $\mP^{3}$
in the fibers over $\mathsf{q}_{i}\,(1\leq i\leq5)$ in the proof).
In particular, $\widetilde{\mathscr{R}}$ is $\mQ$-factorial and
terminal.
\end{prop}

\begin{proof}
By Proposition \ref{prop:f'}, $\widetilde{\mathscr{R}}$ is a $\mP^{7}$-bundle
outside $\cup_{i=1}^{5}\widetilde{f}{}_{\mathscr{R}}^{-1}(\mathsf{q}_{i})$,
thus is smooth there. 

By the $\mathfrak{S}_{5}$-action, we have only to determine the singularities
of $\widetilde{\mathscr{R}}$ on the $\widetilde{f}{}_{\mathscr{R}}$-fiber
over $\mathsf{q}_{3}=$the $\widetilde{r}_{34}$-point. Near the fiber,
setting $\widetilde{r}_{34}=1$, we may consider the following coordinate
change: 
\begin{align*}
\widetilde{x}'_{1} & :=\widetilde{x}_{1}+\widetilde{x}_{4}(1+\widetilde{r}_{35}),\widetilde{x}'_{2}:=\widetilde{x}_{2}-\widetilde{x}_{3}\widetilde{r}_{24}-\widetilde{x}_{4}(1+\widetilde{r}_{35}),\widetilde{x}'_{5}:=\widetilde{x}_{1}\widetilde{r}_{35}+\widetilde{x}_{5}(1+\widetilde{r}_{35}),\\
\widetilde{y}'_{1} & :=\widetilde{y}_{1}+\widetilde{y}_{4}(1+\widetilde{r}_{35}),\widetilde{y}'_{2}:=\widetilde{y}_{2}-\widetilde{y}_{3}\widetilde{r}_{24}-\widetilde{y}_{4}(1+\widetilde{r}_{35}),\widetilde{y}'_{5}:=\widetilde{y}_{1}\widetilde{r}_{35}+\widetilde{y}_{5}(1+\widetilde{r}_{35}),\\
\widetilde{w}' & :=(1+\widetilde{r}_{35})(\widetilde{w}+\widetilde{q}_{34}),
\end{align*}
which replace $\widetilde{x}_{1},\widetilde{x}_{2},\widetilde{x}_{5},\widetilde{y}_{1},\widetilde{y}_{2},\widetilde{y}_{5},w$
respectively. Then, using (\ref{eq:Rtilde}), we can check that $\widetilde{\mathscr{R}}$
is defined by the five Pl\"ucker relations for the skew-symmetric
matrix {\footnotesize
\begin{equation}
\left(\begin{array}{ccccc}
0 & \widetilde{w}' & \widetilde{x}_{1}' & \widetilde{x}_{2}' & \widetilde{x}_{5}'\\
 & 0 & \widetilde{y}_{1}' & \widetilde{y}_{2}' & \widetilde{y}_{5}'\\
 &  & 0 & \widetilde{r}_{24} & -\widetilde{r}_{15}\\
 &  &  & 0 & \widetilde{r}_{35}\\
 &  &  &  & 0
\end{array}\right)\label{eq:mat Rtilde}
\end{equation}
}in an open neighborhood of $\{\widetilde{r}_{15}=\widetilde{r}_{24}=\widetilde{r}_{35}=0\}$
in 
\[
\mP(\widetilde{x}_{1}',\widetilde{x}_{2}',\widetilde{x}_{3},\widetilde{x}_{4},\widetilde{x}_{5}',\widetilde{y}_{1}',\widetilde{y}_{2}',\widetilde{y}_{3},\widetilde{y}_{4},\widetilde{y}_{5}',\widetilde{w})\times\mA(\widetilde{r}_{15},\widetilde{r}_{24},\widetilde{r}_{35})\simeq\mP(1^{10},2)\times\mA^{3}.
\]
 Therefore, in the fiber over $\mathsf{q}_{3}$, $\widetilde{\mathscr{R}}$
has a $\nicefrac{1}{2}(1^{10})$-singularity at the $\widetilde{w}$-point,
and $c({\rm G}(2,5))$-singularities along the common zero locus of
the entries of the matrix (\ref{eq:mat Rtilde}), i.e.,
\[
\widetilde{\Gamma}_{3}:=\{\widetilde{x}_{1}+\widetilde{x}_{2}=\widetilde{x}_{1}+\widetilde{x}_{4}=\widetilde{x}_{5}=\widetilde{y}_{1}+\widetilde{y}_{2}=\widetilde{y}_{1}+\widetilde{y}_{4}=\widetilde{y}_{5}=\widetilde{w}\widetilde{r}_{34}+\widetilde{q}_{34}=0\}\simeq\mP^{3}.
\]
 By the $\mathfrak{S}_{5}$-action, we see that $\widetilde{\mathscr{R}}$
has $c({\rm G}(2,5))$-singularities along the following loci in the
$\widetilde{f}{}_{\mathscr{R}}$-fibers over $\mathsf{q}_{1},\dots,\mathsf{q}_{5}$
as follows:
\begin{align*}
\mathsf{q_{1}}:\widetilde{\Gamma}_{1} & :=\{\widetilde{x}_{2}=\widetilde{x}_{3}=\widetilde{x}_{4}=\widetilde{y}_{2}=\widetilde{y}_{3}=\widetilde{y}_{4}=\widetilde{w}\widetilde{r}_{15}+\widetilde{q}_{15}=0\},\\
\mathsf{q_{2}}:\widetilde{\Gamma}_{2} & :=\{\widetilde{x}_{1}=\widetilde{x}_{3}=\widetilde{x}_{5}=\widetilde{y}_{1}=\widetilde{y}_{3}=\widetilde{y}_{5}=\widetilde{w}\widetilde{r}_{24}+\widetilde{q}_{24}=0\},\\
\mathsf{q_{4}}:\widetilde{\Gamma}_{4} & :=\{\widetilde{x}_{4}=\widetilde{x}_{1}+\widetilde{x}_{2}=\widetilde{x}_{1}-\widetilde{x}_{3}+\widetilde{x}_{5}=\widetilde{y}_{4}=\widetilde{y}_{1}+\widetilde{y}_{2}=\widetilde{y}_{1}-\widetilde{y}_{3}+\widetilde{y}_{5}=\widetilde{w}\widetilde{r}_{35}+\widetilde{q}_{35}=0\},\\
\mathsf{q_{5}}:\widetilde{\Gamma}_{5} & :=\{\widetilde{x}_{1}=\widetilde{x}_{2}+\widetilde{x}_{3}=\widetilde{x}_{4}+\widetilde{x}_{5}=\widetilde{y}_{1}=\widetilde{y}_{2}+\widetilde{y}_{3}=\widetilde{y}_{4}+\widetilde{y}_{5}=\widetilde{w}\widetilde{r}_{24}+\widetilde{q}_{24}=0\}.
\end{align*}
To show that  $\widetilde{\mathscr{R}}$ is $\mQ$-factorial and
terminal, it suffices to show that so is the $c({\rm G}(2,5))$-singularity. Let $(V,o)$ be an analytic germ of the $c({\rm G}(2,5))$-singularity and $p\colon \widetilde{V}\to V$ be the blow-up of $V$ at $o$, which is a resolution of singularity. Let $E\simeq \rm{G}(2,5)$ be the $p$-exceptional divisor. Since $o$ is a cone singularity, we have $-E|_E=\sO_{\rm{G}(2,5)}(1)$. Moreover, since $-K_{\rm{G}(2,5)}=\sO_{\rm{G}(2,5)}(5)$, we have $-K_{\widetilde{V}}|_E=\sO_{\rm{G}(2,5)}(4)$. As the Picard number of $\rm{G}(2,5)$ is $1$, this implies that $p$ is a $K$-negative extremal divisorial contraction. Therefore $(V,o)$ is $\mQ$-factorial and
terminal by \cite[Prop.5-1-5]{KMM}.
\end{proof}

\subsubsection{\textbf{The birational map $\widetilde{{\rm \mathscr{R}}}\protect\dashrightarrow\mathsf{R}$\label{subsec:RTildePOmega}}}

By Table \ref{Table:FR} and the construction as in Subsection \ref{subsec:LR and S3},
we obtain the rational map $\mF_{\widetilde{\mathscr{R}}}\dashrightarrow{\rm S}_{3}\times\mP(V^{\oplus2})$
defined by 
\begin{equation}
[\widetilde{r};\widetilde{x},\widetilde{y},\widetilde{w}]\mapsto[f_{5}(\widetilde{r}),-f_{4}(\widetilde{r}),f_{3}(\widetilde{r}),-f_{2}(\widetilde{r}),f_{1}(\widetilde{r})]\times[\widetilde{x},\widetilde{y}].\label{eq:Rtilde to sfR}
\end{equation}
By (\ref{eq:P(Ogema)}) and the first 2 equations of (\ref{eq:Rtilde}),
this induces the rational map $\widetilde{\mathscr{R}}\dashrightarrow\mathsf{R}$,
where we recall $\mathsf{R}=\mP(\Omega_{\mP(V^{*})}(1)^{\oplus2}|_{{\rm S}_{3}}).$
Let $\mathsf{f}_{\mathscr{R}}\colon\mathsf{R}\to{\rm S}_{3}$ be the
natural projection. We can check the following square is commutative:

\begin{equation}\label{eq:sqGres}
\xymatrix{\widetilde{\mathscr{R}}\ar@{-->}[r]\ar[d]_{\widetilde{f}_{\mathscr{R}}} &\mathsf{R}\ar[d]^{\mathsf{f}_{\mathscr{R}}}\\
L_\mathscr{R}\ar@{-->}[r]& S_3.}
\end{equation}

By Proposition \ref{prop:Segre} and (\ref{eq:sqGres}), we obtain
the following:
\begin{prop}
\label{prop:Rtilde and P(Omega2)}The rational map $\widetilde{{\rm \mathscr{R}}}\dashrightarrow\mathsf{R}$
is an isomorphism outside $\widetilde{f}_{\mathscr{R}}^{-1}(\ell_{ij})\,(1\leq i<j\leq5)$.
The $\widetilde{f}_{\mathscr{R}}$-fiber over a point of $\ell_{ij}\setminus\{\mathsf{q}_{i},\mathsf{q}_{j}\}$
is mapped isomorphically to the $\mathsf{f}_{\mathscr{R}}$-fiber
over the singular point $\mathsf{q}_{ij}$.
\end{prop}

\begin{proof}
Note that $L_{\mathscr{R}}\setminus\{\mathsf{q}_{1},\dots,\mathsf{q}_{5}\}=\cup_{i=1}^{5}\{f_{i}(\widetilde{r})\not=0\}$.
We only consider the map $\widetilde{{\rm \mathscr{R}}}\dashrightarrow\mathsf{R}$
over $\{f_{1}(\widetilde{r})\not=0\}$ and $\{\mathsf{z}_{5}\not=0\}$
since we can treat it in the same way on the other charts by the $\mathfrak{S}_{5}$-action.
As in the proof of Proposition \ref{prop:f'}, over the open subset
$\{\widetilde{f}_{1}\not=0$\}, $\widetilde{\mathscr{R}}$ is isomorphic
to $\mP^{7}\times\{\widetilde{f}_{1}\not=0\}$, where the coordinates
of $\mP^{7}$ are $\widetilde{x}_{i},\widetilde{y}_{i}\,(1\leq i\le4)$.
Correspondingly, we can similarly show that, over $\{\mathsf{z}_{5}\not=0\}$,
$\mathsf{R}$ is isomorphic to $\mP^{7}\times\{\mathsf{z}_{5}\not=0\}$,
where the coordinates of $\mP^{7}$ are $\mathsf{x}_{i},\mathsf{y}_{i}\,(1\leq i\le4)$.
Therefore $\widetilde{{\rm \mathscr{R}}}\dashrightarrow\mathsf{R}$
over $\{\widetilde{f}_{1}\not=0\}$ and $\{\mathsf{z}_{5}\not=0\}$
is the direct product of the identity map of $\mP^{7}$ and the morphism
$\mu_{{\rm S}}\colon L_{\mathscr{R}}\dashrightarrow{\rm S}_{3}$.
This implies the assertion by noting the description of $\mu_{{\rm S}}$.
\end{proof}

\subsubsection{\textbf{The flopping contraction $\widetilde{g}_{\mathscr{R}}\colon\widetilde{\mathscr{R}}\to\overline{\mathscr{R}}$\label{subsec:The-flopping-contraction}}}

We define the morphism $\widetilde{g}_{\mathscr{R}}\colon\widetilde{\mathscr{R}}\to\overline{\mathscr{R}}$
by setting 

\begin{align}
\widetilde{\mathscr{R}} & \ni[\widetilde{r}_{15},\widetilde{r}_{24},\widetilde{r}_{34},\widetilde{r}_{35};\widetilde{x},\widetilde{y},\widetilde{w}]\mapsto[\widetilde{x},\widetilde{y},\widetilde{w}\widetilde{r}_{15},\widetilde{w}\widetilde{r}_{24},\widetilde{w}\widetilde{r}_{34},\widetilde{w}\widetilde{r}_{35}]\in\mP(1^{10},2^{4}),\label{eq:hR}
\end{align}
where we can check that the image of $\widetilde{g}_{\mathscr{R}}$
is certainly contained in $\overline{\mathscr{R}}$ by (\ref{eq:R})
and (\ref{eq:Rtilde}). We may also check that $\widetilde{g}_{\mathscr{R}}$
is $\mathfrak{S}_{5}$-equivariant. We will obtain properties of $\widetilde{g}_{\mathscr{R}}\colon\widetilde{\mathscr{R}}\to\overline{\mathscr{R}}$
and $\overline{\mathscr{R}}$ as follows:
\begin{lem}
\label{lem:isomgtildeR}The morphism $\widetilde{g}_{\mathscr{R}}$
is surjective and is an isomorphism outside $\widetilde{g}_{\mathscr{R}}^{-1}(\Pi)$.
\end{lem}

\begin{proof}
Let $[\overline{x},\overline{y},\overline{r}_{15},\overline{r}_{24},\overline{r}_{34},\overline{r}_{35}]$
be a point of $\overline{\mathscr{R}}$. The morphism $\widetilde{g}_{\mathscr{R}}$
is surjective since it holds that 
\[
\widetilde{g}_{\mathscr{R}}([\overline{r}_{15},\overline{r}_{24},\overline{r}_{34},\overline{r}_{35};\overline{x},\overline{y},1])=[\overline{x},\overline{y},\overline{r}_{15},\overline{r}_{24},\overline{r}_{34},\overline{r}_{35}],
\]
 and $[\overline{r}_{15},\overline{r}_{24},\overline{r}_{34},\overline{r}_{35};\overline{x},\overline{y},1]\in\widetilde{\mathscr{R}}$
by (\ref{eq:MidEq}), Lemma \ref{lem:TwoRelonRvar}, and (\ref{eq:Rtilde}).
Now we show that $\widetilde{g}_{\mathscr{R}}$ is injective over
$\widetilde{\mathscr{R}}\setminus\widetilde{g}_{\mathscr{R}}^{-1}(\Pi)$.
Let $[\overline{x},\overline{y},\overline{r}_{15},\overline{r}_{24},\overline{r}_{34},\overline{r}_{35}]$
be a point of $\overline{\mathscr{R}}\setminus\Pi$. Note that $(\overline{r}_{15},\overline{r}_{24},\overline{r}_{34},\overline{r}_{35})\not=\bm{o}$.
Assume that 
\[
\widetilde{g}_{\mathscr{R}}([\widetilde{r}_{15},\widetilde{r}_{24},\widetilde{r}_{34},\widetilde{r}_{35};\widetilde{x},\widetilde{y},\widetilde{w}])=[\overline{x},\overline{y},\overline{r}_{15},\overline{r}_{24},\overline{r}_{34},\overline{r}_{35}]
\]
 for a point $[\widetilde{r}_{15},\widetilde{r}_{24},\widetilde{r}_{34},\widetilde{r}_{35};\widetilde{x},\widetilde{y},\widetilde{w}]\in\widetilde{\mathscr{R}}$.
Since $(\overline{r}_{15},\overline{r}_{24},\overline{r}_{34},\overline{r}_{35})\not=\bm{o}$,
we have $\widetilde{w}\not=0$. By the $(\mC^{*})^{2}$-action, we
may assume that $w=1$. Therefore there exists $\alpha\in\mC^{*}$
such that $\widetilde{x}=\alpha\overline{x}$, $\widetilde{y}=\alpha\overline{y}$
and $(\widetilde{r}_{15},\widetilde{r}_{24},\widetilde{r}_{34},\widetilde{r}_{35})=\alpha^{2}(\overline{r}_{15},\overline{r}_{24},\overline{r}_{34},\overline{r}_{35})$.
By Table \ref{Table:FR}, we see that
\[
[\widetilde{r}_{15},\widetilde{r}_{24},\widetilde{r}_{34},\widetilde{r}_{35};\widetilde{x},\widetilde{y},1]=[\alpha^{-2}\widetilde{r}_{15},\alpha^{-2}\widetilde{r}_{24},\alpha^{-2}\widetilde{r}_{34},\alpha^{-2}\widetilde{r}_{35};\alpha^{-1}\widetilde{x},\alpha^{-1}\widetilde{y},1]=[\overline{r}_{15},\overline{r}_{24},\overline{r}_{34},\overline{r}_{35};\overline{x},\overline{y},1].
\]
 Therefore the injectivity of $\widetilde{g}_{\mathscr{R}}$ follows. 
\end{proof}
\begin{prop}
\label{prop:RvarGor}The scheme $\overline{\mathscr{R}}$ is a $\mQ$-Gorenstein
variety of codimension $3$ in $\mP(1^{10},2^{4})$.
\end{prop}

\begin{proof}
The scheme $\overline{\mathscr{R}}$ is a variety by Lemma \ref{lem:isomgtildeR}
since so is $\widetilde{\mathscr{R}}$ by Proposition \ref{prop:QfacR}.
Since $\dim\widetilde{\mathscr{R}}=10$, $\overline{\mathscr{R}}$
is of codimension $3$ in $\mP(1^{10},2^{4})$ by Lemma \ref{lem:isomgtildeR}.
Since $\overline{\mathscr{R}}$ is defined by the five Pl\"ucker
relations as in Subsection \ref{subsec:The-midpoint}, we see that
$\overline{\mathscr{R}}$ is $\mQ$-Gorenstein by \cite{BE}. Therefore
we have shown the assertion.
\end{proof}
Setting $q_{ij}:=q_{ij}(x,y)\,(1\leq i<j\leq5)$, we define 
\begin{equation}
\begin{cases}
 & S_{1}(x,y):=-q_{15}q_{24}+q_{45}q_{13}-q_{15}q_{14}-q_{15}q_{25}+q_{12}q_{15}+q_{13}q_{15}-q_{15}^{2},\\
 & S_{2}(x,y):=q_{24}q_{34}+q_{24}^{2}-q_{12}q_{24}+q_{14}q_{34}+q_{14}q_{24}+q_{24}q_{25}-q_{12}q_{34}+q_{15}q_{24},\\
 & S_{3}(x,y):=q_{15}q_{34}-q_{13}q_{34}-q_{13}q_{24}+q_{24}q_{35},\\
 & S_{4}(x,y):=-q_{15}q_{34}-q_{15}q_{23}-q_{24}q_{35}.
\end{cases}\label{eq:S}
\end{equation}

\begin{lem}
\label{lem:Mq}Setting $\overline{q}_{ij}:=q_{ij}(\overline{x},\overline{y})\,(1\leq i<j\leq5)$,
we define the matrix

{\footnotesize
\[
M_{q}:=\left(\begin{array}{cccc}
0 & -\overline{q}_{13} & \overline{q}_{12}-\overline{q}_{14}-\overline{q}_{24} & -\overline{q}_{14}-\overline{q}_{24}\\
\overline{q}_{23} & 0 & -\overline{q}_{15}-\overline{q}_{25} & \overline{q}_{12}+\overline{q}_{13}-\overline{q}_{15}-\overline{q}_{25}\\
\overline{q}_{24} & \overline{q}_{15} & 0 & \overline{q}_{14}\\
\overline{q}_{34} & 0 & \overline{q}_{15}+\overline{q}_{45} & -\overline{q}_{14}+\overline{q}_{45}\\
0 & -\overline{q}_{35} & \overline{q}_{25}-\overline{q}_{45} & -\overline{q}_{24}-\overline{q}_{34}-\overline{q}_{45}
\end{array}\right),
\]
}which comes from the coefficients of the linear parts of (\ref{eq:MidEq})
with respect to $r_{15}$, $r_{24}$, $r_{34}$ and $r_{35}$. The
following assertions hold:

\begin{enumerate}[$(1)$]

\item It holds that $\rank M_{q}\leq3$ at any point of $\Pi$.

\item The locus of points with $\rank M_{q}\leq2$ in $\Pi$ coincides
with 
\[
\overline{S}_{\mathscr{R}}:=\{S_{1}(\overline{x},\overline{y})=S_{2}(\overline{x},\overline{y})=S_{3}(\overline{x},\overline{y})=S_{4}(\overline{x},\overline{y})=0\}.
\]

\item The locus of points with $\rank M_{q}=1$ in $\Pi$ is the
union of the following $10$ loci:{\footnotesize
\begin{align}
 & \overline{\Delta}_{1}^{o}:=\{\overline{x}_{1}=\overline{x}_{2}=\overline{x}_{3}=\overline{y}_{1}=\overline{y}_{2}=\overline{y}_{3}=0,\overline{q}_{45}\not=0\}.\label{eq:10loci}\\
 & \overline{\Delta}_{2}^{o}:=\{\overline{x}_{1}=\overline{x}_{2}=\overline{x}_{4}=\overline{y}_{1}=\overline{y}_{2}=\overline{y}_{4}=0,\overline{q}_{35}\not=0\}.\nonumber \\
 & \overline{\Delta}_{3}^{o}:=\{\overline{x}_{1}=\overline{x}_{4}=\overline{x}_{5}=\overline{y}_{1}=\overline{y}_{4}=\overline{y}_{5}=0,\overline{q}_{23}\not=0\}.\nonumber \\
 & \overline{\Delta}_{4}^{o}:=\{\overline{x}_{3}=\overline{x}_{5}=\overline{x}_{2}-\overline{x}_{4}=\overline{y}_{3}=\overline{y}_{5}=\overline{y}_{2}-\overline{y}_{4}=0,\overline{q}_{12}\not=0\}.\nonumber \\
 & \overline{\Delta}_{5}^{o}:=\{\overline{x}_{1}=\overline{x}_{2}-\overline{x}_{4}=\overline{x}_{3}-\overline{x}_{5}=\overline{y}_{1}=\overline{y}_{2}-\overline{y}_{4}=\overline{y}_{3}-\overline{y}_{5}=0,\overline{q}_{23}\not=0\}.\nonumber \\
 & \overline{\Delta}_{6}^{o}:=\{\overline{x}_{3}+\overline{x}_{4}=\overline{x}_{2}-\overline{x}_{4}=\overline{x}_{1}+\overline{x}_{4}+\overline{x}_{5}=\overline{y}_{3}+\overline{y}_{4}=\overline{y}_{2}-\overline{y}_{4}=\overline{y}_{1}+\overline{y}_{4}+\overline{y}_{5}=0,\overline{q}_{12}\not=0\}.\nonumber \\
 & \overline{\Delta}_{7}^{o}:=\{\overline{x}_{3}=\overline{x}_{4}=\overline{x}_{1}+\overline{x}_{5}=\overline{y}_{3}=\overline{y}_{4}=\overline{y}_{1}+\overline{y}_{5}=0,\overline{q}_{12}\not=0\}.\nonumber \\
 & \overline{\Delta}_{8}^{o}:=\{\overline{x}_{4}=\overline{x}_{5}=\overline{x}_{2}+\overline{x}_{3}=\overline{y}_{4}=\overline{y}_{5}=\overline{y}_{2}+\overline{y}_{3}=0,\overline{q}_{12}\not=0\}.\nonumber \\
 & \overline{\Delta}_{9}^{o}:=\{\overline{x}_{3}=\overline{x}_{1}+\overline{x}_{2}=\overline{x}_{1}+\overline{x}_{4}+\overline{x}_{5}=\overline{y}_{3}=\overline{y}_{1}+\overline{y}_{2}=\overline{y}_{1}+\overline{y}_{4}+\overline{y}_{5}=0,\overline{q}_{12}\not=0\}.\nonumber \\
 & \overline{\Delta}_{10}^{o}:=\{\overline{x}_{5}=\overline{x}_{2}+\overline{x}_{3}=\overline{x}_{1}-\overline{x}_{3}=\overline{y}_{5}=\overline{y}_{2}+\overline{y}_{3}=\overline{y}_{1}-\overline{y}_{3}=0,\overline{q}_{12}\not=0\},\nonumber 
\end{align}
}each of which is the complement of a quadric surface in $\mP^{3}$. 

\item The locus of points with $\rank M_{q}=0$ is 
\[
\overline{\Delta}_{0}:=\{\overline{q}_{ij}=0\,(1\leq i<j\leq5)\}.
\]

\end{enumerate}
\end{lem}

\begin{proof}
All the assertions follows from direct computations.
\end{proof}
\begin{prop}
\label{prop:SingRpre} The singular locus of $\overline{\mathscr{R}}$
is the union of

\begin{enumerate}[$(a)$]

\item $\overline{S}_{\mathscr{R}}$, 

\item the images of $\widetilde{\Gamma}_{1},\dots,\widetilde{\Gamma}_{5}$
(cf.~the proof of Proposition \ref{prop:QfacR}), which are isomorphic
to $\mP^{3}$ and is denoted by $\overline{\Gamma}_{1},\dots,\overline{\Gamma}_{5}$,
and 

\item$\overline{\mathscr{R}}\cap\{\overline{x}=\overline{y}=0\}=\{f_{1}(\overline{r})=\cdots=f_{5}(\overline{r})=0\}$
consisting of $5$ points. 

\end{enumerate}

The variety $\overline{\mathscr{R}}$ is normal and has hypersurface
singularities along the locus of points of $\Pi$ with $\rank M_{q}=2$,
and $\nicefrac{1}{2}$-singularities at the five points in $\overline{\mathscr{R}}\cap\{\overline{x}=\overline{y}=0\}$.
\end{prop}

\begin{proof}
By Proposition \ref{prop:QfacR} and Lemma \ref{lem:isomgtildeR},
the singular locus of $\overline{\mathscr{R}}\setminus\Pi$ is determined
as in (b) and (c), and $\overline{\mathscr{R}}$ has $\nicefrac{1}{2}$-singularities
at the five points in $\overline{\mathscr{R}}\cap\{\overline{x}=\overline{y}=0\}$. 

Now we fix a point $\mathsf{p}:=(\overline{x},\overline{y})\in\Pi$.
Noting that $\overline{\mathscr{R}}$ is defined by $RF_{i}(\overline{x},\overline{y},\overline{r})=0\,(1\leq i\leq5)$
and the singularities of $\overline{\mathscr{R}}$ along $\Pi$ are
determined only by the linear parts of $RF_{i}(\overline{x},\overline{y},\overline{r})$
with respect to $\overline{r}_{15}$, $\overline{r}_{24}$, $\overline{r}_{34}$,
$\overline{r}_{35}$, we see that $\overline{\mathscr{R}}$ is singular
at $\mathsf{p}\in\Pi$ if and only if $\rank M_{q}\leq2$ (note also
that $\overline{\mathscr{R}}$ is codimension $3$ in $\mP(1^{10},2^{4})$
by Proposition \ref{prop:RvarGor}, and is Gorenstein near $\Pi$).
Therefore, by Lemma \ref{lem:Mq}, the singular locus of $\overline{\mathscr{R}}$
in $\Pi$ is $\overline{S}_{\mathscr{R}}$. We can check directly
that ${\rm codim}\,\overline{S}_{\mathscr{R}}=3$ in $\overline{\mathscr{R}}$,
hence $\overline{\mathscr{R}}$ is normal by Proposition \ref{prop:RvarGor}.
The condition that $\rank M_{q}=2$ at $\mathsf{p}$ implies that
$\overline{\mathscr{R}}$ has a hypersurface singularity at $\mathsf{p}$.
\end{proof}
\begin{prop}
\label{prop:hR}The morphism $\widetilde{g}_{\mathscr{R}}\colon\widetilde{\mathscr{R}}\to\overline{\mathscr{R}}$
is a flopping contraction such that the image of the exceptional locus
is $\overline{S}_{\mathscr{R}}$ in $\Pi$. The $\widetilde{g}_{\mathscr{R}}$-fibers
over general points of $\overline{S}_{\mathscr{R}}$ are isomorphic
to $\mP^{1}$ (in fact, we can determine all the $\widetilde{g}_{\mathscr{R}}$-fibers.
We refer to the proof below for this fact). In particular, $\overline{\mathscr{R}}$
has only terminal singularities.
\end{prop}

\begin{proof}
We fix a point $\mathsf{p}:=[\overline{x},\overline{y}]\in\Pi$. Note
that $\widetilde{g}_{\mathscr{R}}^{-1}(\Pi)=\{w=0\}$. Then, by (\ref{eq:Rtilde}),
the $\widetilde{g}_{\mathscr{R}}$-fiber over $\mathsf{p}$ is defined
by

\begin{equation}
\begin{cases}
\overline{x}_{1}f_{5}(\widetilde{r})-\overline{x}_{2}f_{4}(\widetilde{r})+\overline{x}_{3}f_{3}(\widetilde{r})-\overline{x}_{4}f_{2}(\widetilde{r})+\overline{x}_{5}f_{1}(\widetilde{r})=0,\\
\overline{y}_{1}f_{5}(\widetilde{r})-\overline{y}_{2}f_{4}(\widetilde{r})+\overline{y}_{3}f_{3}(\widetilde{r})-\overline{y}_{4}f_{2}(\widetilde{r})+\overline{y}_{5}f_{1}(\widetilde{r})=0
\end{cases}\label{eq:Roverline}
\end{equation}
in the linear subspace 
\[
\left\{ M_{q}\left(\begin{array}{c}
\widetilde{r}_{15}\\
\widetilde{r}_{24}\\
\widetilde{r}_{34}\\
\widetilde{r}_{35}
\end{array}\right)=\bm{o}\right\} \simeq\mP^{3-\rank M_{q}}\subset\mP(\widetilde{r}_{15},\widetilde{r}_{24},\widetilde{r}_{34},\widetilde{r}_{35}),
\]
where $\bm{o}$ is the zero vector. We can determine the fiber $\widetilde{g}_{\mathscr{R}}^{-1}(\mathsf{p})$
as follows according to $\rank M_{q}$:

If $\rank M_{q}=3$ at $\mathsf{p}$, then the $\widetilde{g}_{\mathscr{R}}$-fiber
over $\mathsf{p}$ consists of at most one point, hence $\widetilde{g}_{\mathscr{R}}$
must be an isomorphism over $\mathsf{p}$ by the Zariski main theorem
(note that $\overline{\mathscr{R}}$ is normal by Proposition \ref{prop:SingRpre}).

Assume that $\rank M_{q}=2$ at $\mathsf{p}$. We have shown that
$\overline{S}_{\mathscr{R}}$ coincides with ${\rm (Sing}\,\overline{\mathscr{R}})\cap\Pi$
in Proposition \ref{prop:SingRpre}. We can directly check that $(\cup_{i=1}^{5}\overline{\Gamma}_{i})\cap\Pi$
is contained in the locus of $\rank M_{q}=0$. Therefore the $\widetilde{g}_{\mathscr{R}}$-fiber
over $\mathsf{p}$ with $\rank M_{q}=2$ is contained in the regular
locus of $\widetilde{\mathscr{R}}$ by Proposition \ref{prop:QfacR}
and hence must be non-trivial. Since it is at most $\mP^{1}$, it
must coincides with $\mP^{1}$ (otherwise, it would consists of a
finite number of points, hence $\widetilde{g}_{\mathscr{R}}$ is an
isomorphism over $\mathsf{p}$ by the Zariski main theorem, a contradiction).

If $\rank M_{q}=1$ at $\mathsf{p}$, then we can directly check that
the fiber over $\mathsf{p}$ is isomorphic to $\mP^{2}$ using (\ref{eq:10loci}). 

Assume that $\rank M_{q}=0$ at $\mathsf{p}$. Then the 3rd--7th
equations of (\ref{eq:Roverline}) are identically zero and the first
two equations determine a quadric surface. Thus the fiber over $\mathsf{p}$
is a quadric surface.

Note that $\widetilde{\mathscr{R}}$ is terminal by Proposition \ref{prop:QfacR}.
Since the codimension of $\overline{S}_{\mathscr{R}}$ is $3$ in
$\overline{\mathscr{R}}$, we see that $\widetilde{g}_{\mathscr{R}}$
is a small contraction by the description of the $\widetilde{g}_{\mathscr{R}}$-fibers
as above. Therefore $\widetilde{g}_{\mathscr{R}}$ is a flopping contraction
since $\overline{\mathscr{R}}$ is $\mQ$-Gorenstein.
\end{proof}
By (\ref{eq:KG}) with $n=10$, we have $-K_{\overline{{\rm G}}_{\mathscr{R}}}=\sO_{\overline{{\rm G}}_{\mathscr{R}}}(20)$.
By (\ref{eq:linear overlineR}) and Proposition \ref{prop:RvarGor},
we obtain
\begin{equation}
-K_{\mathscr{\overline{R}}}=\sO_{\overline{\mathscr{R}}}(8).\label{eq:canonical-divisor-overlineR}
\end{equation}

We denote by $\sO_{\widetilde{\mathscr{R}}}(1)$ the pull-back of
the divisor $\sO_{\overline{\mathscr{R}}}(1)$, and we call it a\textit{
tautological divisor }of $\widetilde{\mathscr{R}}$. By (\ref{eq:canonical-divisor-overlineR})
and Proposition \ref{prop:hR}, we have 

\begin{equation}
-K_{\mathscr{\widetilde{R}}}=\sO_{\widetilde{\mathscr{R}}}(8).\label{eq:tildeR canonicaldiv}
\end{equation}
Hereafter let us study singularities of $\overline{\mathscr{R}}$
along the locus of points of $\Pi$ with $\rank M_{q}=2$. For this,
we need some preparations.

Let $\mathbb{E}_{\widetilde{\mathscr{R}}}\to\mP(\wedge^{2}V)$ be
the toric variety obtained as the quotient of $(\wedge^{2}V\setminus\{\bm{o}\})\times\left((V^{\oplus2}\oplus\mC)\setminus\{\bm{o}\}\right)$
by the $(\mC^{*})^{\oplus2}$-action with similar weights of coordinates
as in Table \ref{Table:FR}, where we take the same coordinates $\widetilde{x},\widetilde{y},w$
of $V^{\oplus2}\oplus\mC$ as in the table, and the coordinates $\widetilde{r}_{ij}$
of $\wedge^{2}V$ as usual. The toric variety $\mathbb{F}_{\widetilde{\mathscr{R}}}\to L_{\mathscr{R}}$ defined as in Subsection \ref{subsec:The-quasi--bundle R}
is nothing but the restriction of $\mathbb{E}_{\widetilde{\mathscr{R}}}\to\mP(\wedge^{2}V)$
over $L_{\mathscr{R}}$. We define the subscheme $\widetilde{{\rm G}}_{\mathscr{R}}$
of $\mathbb{E}_{\widetilde{\mathscr{R}}}$ by 
\begin{equation}
\begin{cases}
\widetilde{x}_{1}{\rm Pf}_{2345}(\widetilde{r})-\widetilde{x}_{2}{\rm Pf}_{1345}(\widetilde{r})+\widetilde{x}_{3}{\rm Pf}_{1245}(\widetilde{r})-\widetilde{x}_{4}{\rm Pf}_{1235}(\widetilde{r})+\widetilde{x}_{5}{\rm Pf}_{1234}(\widetilde{r})=0,\\
\widetilde{y}_{1}{\rm Pf}_{2345}(\widetilde{r})-\widetilde{y}_{2}{\rm Pf}_{1345}(\widetilde{r})+\widetilde{y}_{3}{\rm Pf}_{1245}(\widetilde{r})-\widetilde{y}_{4}{\rm Pf}_{1235}(\widetilde{r})+\widetilde{y}_{5}{\rm Pf}_{1234}(\widetilde{r})=0,
\end{cases}\label{eq:GRP7}
\end{equation}
and 
\begin{align}
 & (\widetilde{r}_{ij}\widetilde{q}_{kl}+\widetilde{r}_{kl}\widetilde{q}_{ij})-(\widetilde{r}_{ik}\widetilde{q}_{jl}+\widetilde{r}_{jl}\widetilde{q}_{ik})+(\widetilde{r}_{il}\widetilde{q}_{jk}+\widetilde{r}_{jk}\widetilde{q}_{il})+w{\rm Pf}_{ijkl}(\widetilde{r})=0\label{eq:GRtilde-rq}
\end{align}
for $(i,j,k,l)=(1,2,3,4),(1,2,3,5),(1,2,4,5),(1,3,4,5),(2,3,4,5)$.
The morphism $\widetilde{\mathscr{R}}\to L_{\mathscr{R}}$ is also
the restriction of $\widetilde{{\rm G}}_{\mathscr{R}}\to\mP(\wedge^{2}V)$
over $L_{\mathscr{R}}$. Similarly to Lemma \ref{lem:isomgtildeR},
we can show that the natural morphism $\psi\colon\widetilde{{\rm G}}_{\mathscr{R}}\to\overline{{\rm G}}_{\mathscr{R}}$
is surjective and is an isomorphism outside $\psi^{-1}(\Pi)$.
\begin{lem}
\label{lem:PE}It holds that $\psi^{-1}(\Pi\setminus\overline{\Delta}_{0})$
is isomorphic to $\mP(q^{*}\sE)$ over $\Pi\setminus\overline{\Delta}_{0}$,
where $\sE$ is defined as in Proposition \ref{prop:bundleE}, and
$q$ is the rational map $\Pi\simeq\mP(V^{\oplus2})\dashrightarrow{\rm G}(2,V)$
defined by $q_{ij}\,(1\leq i<j\leq5)$ (note that this is a morphism
on $\Pi\setminus\overline{\Delta}_{0}$).
\end{lem}

\begin{proof}
Let $(\widetilde{r};\widetilde{x},\widetilde{y},0)$ be a point of
$\psi^{-1}(\Pi)$. If $\widetilde{q}_{12}\not=0$ for example, then,
by (\ref{eq:GRtilde-rq}), we have 
\[
(\widetilde{r}_{12}\widetilde{q}_{kl}+\widetilde{r}_{kl}\widetilde{q}_{12})-(\widetilde{r}_{1k}\widetilde{q}_{2l}+\widetilde{r}_{2l}\widetilde{q}_{1k})+(\widetilde{r}_{1l}\widetilde{q}_{2k}+\widetilde{r}_{2k}\widetilde{q}_{1l})=0
\]
for $(k,l)=(3,4),$(3,5) or $(4,5)$, hence 
\[
\widetilde{r}_{kl}=-\widetilde{r}_{12}\widetilde{q}_{kl}/\widetilde{q}_{12}+(\widetilde{r}_{1k}\widetilde{q}_{2l}/\widetilde{q}_{12}+\widetilde{r}_{2l}\widetilde{q}_{1k}/\widetilde{q}_{12})-(\widetilde{r}_{1l}\widetilde{q}_{2k}/\widetilde{q}_{12}+\widetilde{r}_{2k}\widetilde{q}_{1l}/\widetilde{q}_{12}).
\]
 Substituting these into (\ref{eq:GRP7}), we see that (\ref{eq:GRtilde-rq})
implies (\ref{eq:GRP7}) if $\widetilde{q}_{12}\not=0$. Actually,
this is always the case if $\widetilde{q_{ab}}\not=0$ for some $a,b$
with $1\leq a<b\leq5$. Therefore, $\psi^{-1}(\Pi\setminus\overline{\Delta}_{0})$
is isomorphic to the subscheme in $\mP(\wedge^{2}V)\times(\Pi\setminus\overline{\Delta}_{0})$
defined by (\ref{eq:GRtilde-rq}) with $w=0$. Then, by Proposition
\ref{prop:bundleE}, this is isomorphic to $\mP(q^{*}\sE)$.
\end{proof}
\begin{prop}
\label{prop:hyp double pt} It holds that $\overline{\mathscr{R}}$
has hypersurface singularities with multiplicity $2$ along the locus
of points of $\Pi$ with $\rank M_{q}=2$.
\end{prop}

\begin{proof}
Let $\mathsf{p}$ be a point of $\Pi$ with $\rank M_{q}=2$. In the
proof of Proposition \ref{prop:hR}, we have shown that the $\widetilde{g}_{\mathscr{R}}^{-1}(\mathsf{p})$
is isomorphic to $\mP^{1}$ and it is contained in the smooth locus
of $\widetilde{\mathscr{R}}$. Then we may choose the linear forms
$L_{1}(r),\dots,L_{6}(r)$ such that $L_{\mathscr{R}}=\{L_{1}(\widetilde{r})=\cdots=L_{6}(\widetilde{r})=0\}\subset\mP(\wedge^{2}V)$,
and for 
\[
\widetilde{A}:=\widetilde{{\rm G}}_{\mathscr{R}}\cap\{L_{1}(\widetilde{r})=\cdots=L_{5}(\widetilde{r})=0\},\overline{A}:=\overline{{\rm G}}_{\mathscr{R}}\cap\{L_{1}(\overline{r})=\cdots=L_{5}(\overline{r})=0\},
\]
 and the restriction $\widetilde{g}_{A}\colon\widetilde{A}\to\overline{A}$
of $\psi\colon\widetilde{{\rm G}}_{\mathscr{R}}\to\overline{{\rm G}}_{\mathscr{R}}$,
it holds that $\widetilde{A}$ is smooth and the $\widetilde{g}_{A}$-exceptional
locus $E$ is a $\mP^{1}$-bundle over $\Pi$ near $\mathsf{p}$.
We study the morphism $\widetilde{g}_{A}$ locally near $\mathsf{p}$.
Unless otherwise stated, we will always work locally in this manner
in what follows. Since $\widetilde{g}_{A}$ is isomorphic over $\overline{A}\setminus\Pi$,
$\overline{A}$ is regular in codimension one. Since $\overline{A}$
is $\mQ$-Gorenstein by Proposition \ref{prop:RvarGor}, $\overline{A}$
is normal.

\vspace{5pt}

\noindent Claim. \textit{The variety $\overline{A}$ is smooth along
$\Pi,$ and $\widetilde{g}_{A}\colon\widetilde{A}\to\overline{A}$
is the blow-up along $\Pi$.}

\vspace{5pt}

\noindent \textit{Proof. }Let $\gamma$ be a $\widetilde{g}_{A}$-exceptional
curve. By the negativity lemma, $\gamma$ is numerically positive
for $-E$. We set $e:=-E\cdot\gamma>0$. Note that $\sN_{\gamma/E}\simeq\sO_{\mP^{1}}^{\oplus9}$
and $\sN_{E/\widetilde{A}}|_{\gamma}\simeq\sO_{\mP^{1}}(-e)$. Then,
by the normal bundle sequence
\[
0\to\sN_{\gamma/E}\to\sN_{\gamma/\widetilde{A}}\to\sN_{E/\widetilde{A}}|_{\gamma}\to0,
\]
it holds that $\sN_{\gamma/\widetilde{A}}\simeq\sO_{\mP^{1}}^{\oplus9}\oplus\sO_{\mP^{1}}(-e)$.
This implies that $-K_{\widetilde{A}}\cdot\gamma=2-e$. Since $\gamma$
is numerically trivial for $K_{\widetilde{\mathscr{IR}}}$ by Proposition
\ref{prop:hR}, and numerically positive for $\{L_{6}(\widetilde{r})=0\}$,
it holds that $\gamma$ is numerically positive for $-K_{\widetilde{A}}$.
Therefore we have $e=1$. Applying \cite[Prop.2.4]{AW2}, we see that
$\overline{A}$ is smooth along $\Pi$. Let $\mathcal{I}$ be the
ideal sheaf of $\Pi$ in $\overline{A}$. Note that $\sO_{\widetilde{A}}(-E)|_{\gamma}\simeq\sO_{\widetilde{A}}(-K_{\widetilde{A}})|_{\gamma}\simeq\sO_{\mP^{1}}(-1)$.
This implies that $\widetilde{g}_{A*}\sO_{\widetilde{A}}(-E)=\mathcal{I}$
and, by \cite[Thm.]{AW1}, $\widetilde{g}_{A}^{*}\widetilde{g}_{A*}\sO_{\widetilde{A}}(-E)\to\sO_{\widetilde{A}}(-E)$
is surjective. Therefore $\widetilde{g}_{A}^{-1}\mathcal{I}=\sO_{\widetilde{A}}(-E)$.
By the universal property of blow-up, we see that $\widetilde{g}_{A}$
factors through the blow-up of $\overline{A}$ along $\Pi$. By the
description of $\widetilde{g}_{A}$ and the Zariski main theorem,
$\widetilde{g}_{A}$ coincides with the blow-up of $\overline{A}$
along $\Pi$. $\hfill\square$

The claim implies that $\overline{\mathscr{R}}$ has a hypersurface
singularity at the image of $\gamma$. Now let us take local analytic
coordinates $x_{1},\dots,x_{11}$ of $\overline{A}$ at the image
of $\gamma$ such that $\Pi=\{x_{1}=x_{2}=0\}$ and $\overline{\mathscr{R}}=\{x_{1}^{a}p+x_{2}^{b}q=0\}$
with positive integers $a,b$ and polynomials $p,q$ which are not
divided by $x_{1}$ or $x_{2}$. Computing the blow-up of $\overline{A}$
along $\Pi$ explicitly, we see easily that $x_{1}^{a}p+x_{2}^{b}q$
has multiplicity two at the origin.
\end{proof}
We set $\widetilde{\Pi}:=\widetilde{\mathscr{R}}\cap\{w=0\}$. By
(\ref{eq:hR}), $\widetilde{\Pi}$ is nothing but the strict transform
of $\Pi$ since $\widetilde{g}_{\mathscr{R}}$ is a flopping contraction.
By the weights on the coordinates of $\widetilde{\mathscr{R}}$ as
in Table \ref{Table:FR}, we have the following:
\begin{lem}
\label{lem:2H-L R}The divisor $\widetilde{\Pi}$ on $\widetilde{\mathscr{R}}$
is linearly equivalent to $\sO_{\widetilde{\mathscr{R}}}(2)-\widetilde{f}_{\mathscr{R}}^{*}\sO_{L_{\mathscr{R}}}(1).$
\end{lem}

\subsubsection{\textbf{Glossary of notation\label{subsec:Glossary-of-notation}}}

In this subsection, we prepare notation we need to describe the key
variety $\mathscr{R}$ in Subsection \ref{subsec:Constructing--via unproj}. 

\begin{equation}
\begin{cases}
K_{1}(x,y,r):= & (3q_{12}+3q_{13}-3q_{14}-6q_{15}-2q_{24}-3q_{25}-q_{34})r_{15}\\
 & -2q_{15}r_{24}+2(q_{15}+q_{45})r_{34}+2(q_{14}+3q_{15}+q_{45})r_{35},\\
K_{2}(x,y,r):= & 2q_{24}r_{15}+(-3q_{12}-q_{13}+3q_{14}+2q_{15}+6q_{24}+3q_{25}+3q_{34})r_{24}\\
 & +2(-q_{12}+q_{14}+q_{24})r_{34}+2(-q_{14}-2q_{24})r_{35},\\
K_{3}(x,y,r):= & 2q_{34}r_{15}+(-2q_{13}+2q_{35})r_{24}\\
 & +(-q_{12}-3q_{13}+q_{14}+2q_{15}-2q_{24}+q_{25}-3q_{34}-2q_{45})r_{34}\\
 & +2(q_{14}-2q_{34}-q_{45})r_{35},\\
K_{4}(x,y,r):= & 2(-q_{23}-q_{34})r_{15}-2q_{35}r_{24}+2(-q_{25}+q_{45})r_{34}\\
 & +(q_{12}+q_{13}-q_{14}+2q_{15}-2q_{24}-q_{25}+q_{34}+2q_{45})r_{35}.
\end{cases}\label{eq:K}
\end{equation}

\begin{equation}
\begin{cases}
L_{1}(r):=-r_{15}(3r_{15}+2r_{24}-2r_{34}-6r_{35}),\\
L_{2}(r):=r_{24}(2r_{15}+3r_{24}+2r_{34}-4r_{35}),\\
L_{3}(r):=-r_{34}(-2r_{15}+2r_{24}+3r_{34}+4r_{35}),\\
L_{4}(r):=r_{35}(2r_{15}-2r_{24}+2r_{34}+r_{35}).
\end{cases}\label{eq:L}
\end{equation}
\begin{equation}
U_{i}(x,y,r):=1/30(30S_{i}(x,y)+5K_{i}(x,y,r)+2L_{i}(r))\,(1\leq i\leq4),\label{eq:U}
\end{equation}
where we recall $S_{i}(x,y)$ is defined in (\ref{eq:S}).

\begin{equation}
\begin{cases}
RF_{6}(x,y,r):=r_{0}r_{15}-30U_{1}(x,y,r),\\
RF_{7}(x,y,r):=r_{0}r_{24}-30U_{2}(x,y,r),\\
RF_{8}(x,y,r):=r_{0}r_{34}-30U_{3}(x,y,r),\\
RF_{9}(x,y,r):=r_{0}r_{35}-30U_{4}(x,y,r),
\end{cases}\label{eq:R2}
\end{equation}
where we set $r=(r_{0},r_{15},r_{24},r_{34},r_{35})$.

\subsubsection{\textbf{Constructing $\mathscr{R}$ from $\overline{\mathscr{R}}$
via unprojection\label{subsec:Constructing--via unproj}}}

Now we are in position to construct the key variety $\mathscr{R}$
via unprojection following \cite{PR} and \cite{P}. For this purpose,
we construct the following commutative diagram, which consists of
the resolutions of the ideals of $\Pi$ and $\overline{\mathscr{R}}$.
In the diagram, $P$ is the graded polynomial ring over $\mC$ with
variables $x,y,r_{15},r_{24},r_{34},r_{35}$, and $P^{n}\,(n\in\mN)$
is the direct product of $n$ copies of $P$ (though we consider the
resolutions of the ideal of $\Pi$ and $\overline{\mathscr{R}}$,
we eliminate $\overline{\phantom{a}}$ from the variables to produce
the variety $\mathscr{R}$ in a different ambient space):

\begin{equation}\label{eq:comp}
\xymatrix{0\ar[r]& P\ar[r]^{A_3}\ar[d]_{B_3}& P^5\ar[r]^{A_2}\ar[d]_{B_2}& P^5 \ar[r]^{A_1}\ar[d]_{B_1}& P\ar[d]_{\rm{id}}\\
P\ar[r]_{C_4}& P^4\ar[r]_{C_3} & P^6 \ar[r]_{C_2} & P^4\ar[r]_{C_1}& P.}
\end{equation}Constructions and explanations for this diagram follows:

\vspace{3pt}

\noindent $A_{i},B_{j},C_{k}\,(i,j,k\in\mN)$ are $P$-homomorphisms.
We identify these with the representation matrices with respect to
suitable bases of relevant free $P$-modules which we choose below.

\vspace{3pt}

\noindent The bottom line of (\ref{eq:comp}) is the graded Koszul
complex with respect to the ideal $(r_{15},r_{24},r_{34},r_{35})$
of $\Pi$. We denote the basis of $P^{4}$ in the second term from
the right by $\bm{e}_{1},\bm{e}_{2},\bm{e}_{3},\bm{e}_{4}$, and the
bases of $P$ in the most left term, $P^{4}$ in the second term from
the left, $P^{6}$ respectively as follows:
\begin{align*}
P\,\,: & \bm{e}_{1}\wedge\bm{e}_{2}\wedge\bm{e}_{3}\wedge\bm{e}_{4},\\
P^{4}: & \bm{e}_{2}\wedge\bm{e}_{3}\wedge\bm{e}_{4},\bm{e}_{1}\wedge\bm{e}_{3}\wedge\bm{e}_{4},\bm{e}_{1}\wedge\bm{e}_{2}\wedge\bm{e}_{4},\bm{e}_{1}\wedge\bm{e}_{2}\wedge\bm{e}_{3},\\
P^{6}: & \bm{e}_{1}\wedge\bm{e}_{2},\bm{e}_{1}\wedge\bm{e}_{3},\bm{e}_{1}\wedge\bm{e}_{4},\bm{e}_{2}\wedge\bm{e}_{3},\bm{e}_{2}\wedge\bm{e}_{4},\bm{e}_{3}\wedge\bm{e}_{4}.
\end{align*}
Then it holds that{\footnotesize
\begin{align*}
C_{1} & =\left(\begin{array}{cccc}
r_{15} & r_{24} & r_{34} & r_{35}\end{array}\right),\,C_{2}=\left(\begin{array}{cccccc}
-r_{24} & -r_{34} & -r_{35} & 0 & 0 & 0\\
r_{15} & 0 & 0 & -r_{34} & -r_{35} & 0\\
0 & r_{15} & 0 & r_{24} & 0 & -r_{35}\\
0 & 0 & r_{15} & 0 & r_{24} & r_{34}
\end{array}\right)\\
C_{3} & =\left(\begin{array}{cccc}
0 & 0 & r_{35} & r_{34}\\
0 & r_{35} & 0 & -r_{24}\\
0 & -r_{34} & -r_{24} & 0\\
r_{35} & 0 & 0 & r_{15}\\
-r_{34} & 0 & r_{15} & 0\\
r_{24} & r_{15} & 0 & 0
\end{array}\right),\,C_{4}=\left(\begin{array}{c}
r_{15}\\
-r_{24}\\
r_{34}\\
-r_{35}
\end{array}\right).
\end{align*}
}\vspace{3pt}

\noindent The top line of (\ref{eq:comp}) is the graded minimal
free resolution of the homogeneous ideal of $\mathscr{\overline{R}}$
due to \cite{BE}, where we may choose

{\footnotesize
\begin{align*}
A_{1} & =\left(\begin{array}{ccccc}
RF_{5} & -RF_{4} & RF_{3} & -RF_{2} & RF_{1}\end{array}\right),\\
A_{2} & =\left(\begin{array}{ccccc}
0 & q_{12} & q_{13}+r_{34}+r_{35} & q_{14} & q_{15}+r_{15}\\
 & 0 & q_{23}-r_{34}-r_{35} & q_{24}+r_{24} & q_{25}-r_{35}\\
 &  & 0 & q_{34}+r_{34} & q_{35}+r_{35}\\
 &  &  & 0 & q_{45}\\
 &  &  &  & 0
\end{array}\right),A_{3}=\left(\begin{array}{c}
RF_{5}\\
-RF_{4}\\
RF_{3}\\
-RF_{2}\\
RF_{1}
\end{array}\right),
\end{align*}
}where we set $q_{ij}:=q_{ij}(x,y)\,(1\leq i<j\leq5)$ and $RF_{k}:=RF_{k}(x,y,r)\,(1\leq k\leq5)$.

\vspace{3pt}

\noindent Now we construct the representation matrices $B_{1},B_{2},B_{3}$
such that the diagram (\ref{eq:comp}) is commutative. We set

\begin{align*}
& \qquad \qquad \quad B_3=\begin{pmatrix}{U}_{1}\\-{U}_{2}\\{U}_{3}\\-{U}_{4}\end{pmatrix},B_2=\\&{\tiny{\begin{pmatrix} 
1/3 {q}_{13}&\substack{ -2/3 {q}_{23}\\ + 1/3 ({r}_{34} + {r}_{35})}& 0& 2/3 {q}_{34} +1/3 {r}_{34}& 
-1/3{q}_{35}\\
\substack{1/3 {q}_{13} + {q}_{24} - 2/3 {q}_{15}\\ + 1/3 (-2 {r}_{15} + {r}_{24} + {r}_{34} +{r}_{35})}& \substack{-{q}_{24} + 1/3 {q}_{23} - 2/3 {q}_{25}\\
+ 1/3 (-2 {r}_{24} - {r}_{34} +{r}_{35})}& -{q}_{23} -{q}_{34} -   2/3 {q}_{35}& \substack{-{q}_{24} - 1/3 {q}_{34} - 2/3 {q}_{45}\\ - 1/3 (2 {r}_{24} + {r}_{34})}& -1/3 {q}_{35}\\
{q}_{24} + 2/3 {q}_{14} +1/3{r}_{24}& -1/3 {q}_{24}& -1/3 {q}_{34}& 0& \substack{-{q}_{24} - {q}_{34} -   2/3 {q}_{45}\\ - 1/3 ({r}_{24} + {r}_{34})}\\
1/3 {q}_{15}& \substack{-{q}_{15} - 2/3 {q}_{25}\\ - 1/3 ({r}_{15} - {r}_{35})}& 1/3 {q}_{35}& -{q}_{15} -   2/3 {q}_{45} -1/3 {r}_{15}& 0\\
{q}_{15} + 2/3{r}_{15}& \substack{{q}_{13} - {q}_{15}\\ - 1/3 ({r}_{15} - {r}_{34})}& -{q}_{13} + {q}_{35} - 2/3{r}_{34}& 0& -{q}_{15} - 2/3 {r}_{15}\\
0& \substack{-{q}_{12} + {q}_{14} + {q}_{24}\\ + (2 {r}_{24})/3}& \substack{{q}_{12} - {q}_{14} - {q}_{15} - {q}_{24} - {q}_{25} \\-   1/3 ({r}_{15} + {r}_{24})}& 0& -{q}_{15} - {q}_{45} - 2/3 {r}_{15}
\end{pmatrix}}}
\end{align*}

Then it holds that $C_{3}B_{3}=B_{2}A_{3}$ by direct calculations.

Finally setting

\begin{equation*}
B_1=\footnotesize{\begin{pmatrix}
0 & -{q}_{34}& \substack{{q}_{24}\\+1/3({r}_{24}+{r}_{34})}& -{q}_{23}+2/3{r}_{35}& -2/3{r}_{34}\\-{q}_{35} & 0 & {q}_{15}+2/3{r}_{15} & 0 & -{q}_{13}-2/3{r}_{34}\\
{q}_{25}-{q}_{45} &-{q}_{15}-{q}_{45}-{r}_{15} & -1/3{r}_{15} & {q}_{15}+{q}_{25}+{r}_{15} & \substack{{q}_{12}-{q}_{14}-{q}_{24}\\ +1/3(2{r}_{15}-{r}_{24})}\\
\substack{-{q}_{24}-{q}_{34}-{q}_{45}\\ -{r}_{24}-{r}_{34}}& {q}_{14}-{q}_{45}& {q}_{14} &\substack{-{q}_{12}-{q}_{13}+{q}_{15}+{q}_{25}\\+1/3({r}_{15}-3{r}_{34}-3{r}_{35})}&-{q}_{14}-{q}_{24}-{r}_{24}
\end{pmatrix}},
\end{equation*}we can directly check that $C_{2}B_{2}=B_{1}A_{2}$ and $C_{1}B_{1}=A_{1}$
hold. 

\vspace{3pt}

\noindent Thus we have the commutative diagram (\ref{eq:comp}).
Then we immediately obtain the following applying \cite[Thm.5.5 and 5.6]{P}
(cf.~\cite[Subsec.2.4]{PR}): 
\begin{prop}
\label{prop:Gorcodim4} We set $RF_{k}:=RF_{k}(x,y,r)\,(6\leq k\leq9)$,
which are derived from the map $B_{3}$ in the diagram (\ref{eq:comp}).
It holds that the ideal $(RF_{1},...,RF_{5},RF_{6},...,RF_{9})$ of
the polynomial ring $P[r_{0}]$ is Gorenstein of codimension $4$.
\end{prop}

\begin{defn}
We define 
\[
\mathscr{R}:=\{RF_{1}=\cdots=RF_{9}=0\}\subset\mP(1^{10},2^{4}),
\]
where the weights of entries of $x$, $y$ are $1$ and entries of
$r$ are $2$.
\end{defn}

\begin{cor}
\label{cor:RQGor}The subscheme $\mathscr{R}$ is $\mQ$-Gorenstein
and is of codimension $4$ in $\mP(1^{10},2^{4})$. Moreover, $\mathscr{R}$
has a $\nicefrac{1}{2}$-singularity at the $r_{0}$-point.
\end{cor}

\begin{proof}
The first statement follows by Proposition \ref{prop:Gorcodim4}.
We show the second statement. Putting $r_{0}=1$, we see that $\{RF_{6}=\cdots=RF_{9}=0\}$
is a quasi-smooth variety of codimension 4 in $\mP(1^{10},2^{4})$
near the $r_{0}$-point. Since $\mathscr{R}$ is of codimension 4
in $\mP(1^{10},2^{4})$, it holds that $\mathscr{R}=\{RF_{6}=\cdots=RF_{9}=0\}$
near the $r_{0}$-point. Then $\mathscr{R}$ has a $\nicefrac{1}{2}$-singularity
at the $r_{0}$-point since the weight of $r_{0}$ is $2$ and we
may choose orbifold coordinates $x,y$ whose entries have weight 1.
\end{proof}

\subsubsection{\textbf{$\mathfrak{S}_{6}$-action on $\mathscr{R}$\label{subsec:S6}}}

By the construction of $\mathscr{R}$, we can directly check that
the $\mathfrak{S}_{5}$-action extends to that of $\mathscr{R}$. 
\begin{prop}
The variety $\mathscr{R}$ has the $\mathfrak{S}_{6}$-action extending
the $\mathfrak{S}_{5}$-action.
\end{prop}

\begin{proof}
We set $\sigma_{i}:=(i\,i+1)\,(1\leq i\leq5)$. The actions of $\sigma_{1},\dots,\sigma_{5}$
on $V_{x}$ and $V_{y}$ are linear; the representation matrices for
the actions of $\sigma_{1},\dots,\sigma_{4}$ are the same as in (\ref{eq:standard S5}),
and that for the actions of $\sigma_{5}$ is
\begin{align}
{\footnotesize \left(\begin{array}{ccccc}
0 & 1 & 0 & -1 & 0\\
1 & 0 & 0 & 1 & 0\\
-1 & 0 & 0 & -1 & -1\\
0 & 0 & 0 & 1 & 0\\
0 & -1 & -1 & 0 & 0
\end{array}\right).}\label{eq:sigma5mat}
\end{align}
By the character table of $\mathfrak{S}_{6}$, we can check that $V_{x}$
and $V_{y}$ are (mutually isomorphic) $9$-dimensional irreducible
representations of $\mathfrak{S}_{6}$. 

Let $V_{r}$ be the $5$-dimensional vector space with the coordinates
$r=(r_{0},r_{15},r_{24},r_{34},r_{35})$. The actions of $\sigma_{1},\dots,\sigma_{4}$
on $V_{r}$ are linear; the representation matrices for the actions
of $\sigma_{1},\dots,\sigma_{4}$ are the following 
\begin{align*}
\sigma_{1}:{\footnotesize \left(\begin{array}{ccccc}
-1 & 0 & 0 & 0 & 0\\
0 & 1 & 0 & 0 & 0\\
0 & -1 & -1 & 0 & 0\\
0 & 1 & 0 & -1 & 0\\
0 & 0 & 0 & 0 & -1
\end{array}\right)},\, & \sigma_{2}:{\footnotesize \left(\begin{array}{ccccc}
0 & 0 & 0 & -1 & 0\\
-1 & -1 & 0 & 1 & 0\\
1 & 0 & -1 & -1 & 0\\
-1 & 0 & 0 & 0 & 0\\
0 & 0 & 0 & 0 & -1
\end{array}\right)},\\
\sigma_{3}:{\footnotesize \left(\begin{array}{ccccc}
0 & 0 & -1 & -1 & 0\\
0 & -1 & 0 & 0 & 0\\
-1 & 0 & 0 & 1 & 0\\
0 & 0 & 0 & -1 & 0\\
0 & 0 & 0 & 0 & -1
\end{array}\right)},\, & \sigma_{4}:{\footnotesize \left(\begin{array}{ccccc}
-1 & 0 & 1 & 1 & 0\\
0 & -1 & 0 & 0 & 0\\
0 & 0 & 0 & 1 & 0\\
0 & 0 & 1 & 0 & 0\\
0 & 0 & 0 & 0 & -1
\end{array}\right).}
\end{align*}
The action of $\sigma_{5}$ is determined by the matrix (\ref{eq:sigma5mat})
and the following nonlinear transformation:{\footnotesize
\begin{align*}
\left(\begin{array}{c}
r_{15}\\
r_{24}\\
r_{34}\\
r_{35}\\
r_{0}
\end{array}\right) & \mapsto\left(\begin{array}{ccccc}
-4/5 & -1/5 & 1/5 & 3/5 & 1/30\\
0 & -1 & 0 & 0 & 0\\
0 & 0 & -1 & 0 & 0\\
1/5 & -1/5 & 1/5 & -2/5 & 1/30\\
36/5 & -36/5 & 36/5 & 108/5 & 1/5
\end{array}\right)\left(\begin{array}{c}
r_{15}\\
r_{24}\\
r_{34}\\
r_{35}\\
r_{0}
\end{array}\right)\\
 & +\left(\begin{array}{c}
\nicefrac{1}{2}(q_{12}+q_{13}-q_{14}-2q_{15}+2q_{23}+q_{34}-q_{25})\\
-q_{14}-q_{24}\\
q_{14}-q_{34}-q_{45}\\
\nicefrac{1}{2}(-q_{12}-q_{13}-q_{14}+q_{25}+q_{34})\\
6(q_{12}+q_{13}+q_{15}-q_{23}-q_{24}-q_{25}-q_{34}+q_{45})
\end{array}\right).
\end{align*}
}We can directly check that the following hold:
\[
\sigma_{i}^{2}={\rm id}\,(1\leq i\leq5),\,\sigma_{i}\sigma_{i+1}\sigma_{i}=\sigma_{i+1}\sigma_{i}\sigma_{i+1}\,(1\leq i\leq4),\sigma_{i}\sigma_{j}=\sigma_{j}\sigma_{i}\,(|i-j|\geq2),
\]
and that $\sigma_{i}\,(1\leq i\leq5)$ stabilize $\mathscr{R}$. Thus
the proposition follows.
\end{proof}
\begin{rem}
By the character table of $\mathfrak{S}_{6}$, we can directly check
that the equations $RF_{1},\dots,RF_{9}$ generates the $9$-dimensional
irreducible representation of $\mathfrak{S}_{6}$ such that the trace
of the representation matrix of $\sigma_{1}=(12)$ is $3$. 
\end{rem}

\subsubsection{\textbf{The blow-up $\widehat{f}_{\mathscr{R}}\colon\widehat{\mathscr{R}}\to\mathscr{R}$
at the $r_{0}$-point\label{subsec:The-blow-up-} }}

The following construction of $\widehat{\mathscr{R}}$ is inspired
by \cite[Sect.4]{AZ} and \cite[Prop.3.7]{Ca}.

Let $V_{\mathscr{R}}\subset\wedge^{2}V$ be the $4$-dimensional subspace
defined as in Subsection \ref{subsec:The-quasi--bundle R}. We define
the toric variety $\mathbb{F}_{\mathscr{\widehat{\mathscr{R}}}}$
as the quotient of $(\mC^{2}\setminus\{\bm{o}\})\times\left((V_{\mathscr{R}}\oplus V^{\oplus2})\setminus\{\bm{o}\}\right)$
by the action of $(\mC^{*})^{2}$ with the weights given as in Table
\ref{eq:wtmatRhat}, where the entries of the second line are the
coordinates of $(\mC^{2}\setminus\{\bm{o}\})\times\left((V_{\mathscr{R}}\oplus V^{\oplus2})\setminus\{\bm{o}\}\right)$,
those of the third and the fourth lines are the weights of the coordinates
by the action of the first and the second factors of $(\mC^{*})^{2}$
respectively (we follow \cite[Sect.4]{AZ} as for the convention on
the weights). 
\begin{table}[h]
\begin{tabular}{|c|c|c|c|c|}
\hline 
space & \multicolumn{2}{c|}{$\mC^{2}\setminus\{\bm{o}\}$} & \multicolumn{2}{c|}{$(V_{\mathscr{R}}\oplus V^{\oplus2})\setminus\{\bm{o}\}$}\tabularnewline
\hline 
coordinates & $\widehat{w}$ & $\widehat{r}_{0}$ & $\widehat{x}$,$\widehat{y}$ & $\widehat{r}_{15},\widehat{r}_{24},\widehat{r}_{34},\widehat{r}_{35}$\tabularnewline
\hline 
1st weights & $0$ & $2$ & $1$ & $2$\tabularnewline
\hline 
2nd weights & $-1$ & $-1$ & $0$ & $1$\tabularnewline
\hline 
\end{tabular}

\caption{Toric variety $\mF_{\widehat{\mathscr{R}}}$}
\label{eq:wtmatRhat}
\end{table}

\begin{defn}
We set $\widehat{q}_{ij}:=q_{ij}(\widehat{x},\widehat{y})\,(1\leq i<j\leq5)$.
In the toric variety $\mathbb{F}_{\mathscr{\widehat{\mathscr{R}}}}$,
we define the subscheme $\widehat{\mathscr{R}}$ with the following
equations:
\begin{equation}
\begin{cases}
-\widehat{q}_{13}\widehat{r}_{24}+(\widehat{q}_{12}-\widehat{q}_{14}-\widehat{q}_{24})\widehat{r}_{34}-(\widehat{q}_{14}+\widehat{q}_{24})\widehat{r}_{35}+\widehat{w}f_{1}(\widehat{r})=0,\\
\widehat{q}_{23}\widehat{r}_{15}-(\widehat{q}_{15}+\widehat{q}_{25})\widehat{r}_{34}+(\widehat{q}_{12}+\widehat{q}_{13}-\widehat{q}_{15}-\widehat{q}_{25})\widehat{r}_{35}+\widehat{w}f_{2}(\widehat{r})=0,\\
\widehat{q}_{24}\widehat{r}_{15}+\widehat{q}_{15}\widehat{r}_{24}+\widehat{q}_{14}\widehat{r}_{35}+\widehat{w}f_{3}(\widehat{r})=0,\\
\widehat{q}_{34}\widehat{r}_{15}+(\widehat{q}_{15}+\widehat{q}_{45})\widehat{r}_{34}-(\widehat{q}_{14}-\widehat{q}_{45})\widehat{r}_{35}+\widehat{w}f_{4}(\widehat{r})=0,\\
-\widehat{q}_{35}\widehat{r}_{24}+(\widehat{q}_{25}-\widehat{q}_{45})\widehat{r}_{34}-(\widehat{q}_{24}+\widehat{q}_{34}+\widehat{q}_{45})\widehat{r}_{35}+\widehat{w}f_{5}(\widehat{r})=0,\\
\widehat{r}_{0}\widehat{r}_{15}-\left(30S_{1}(\widehat{x},\widehat{y})+5\widehat{w}K_{1}(\widehat{x},\widehat{y},\widehat{r})+2\widehat{w}^{2}L_{1}(\widehat{r})\right)=0,\\
\widehat{r}_{0}\widehat{r}_{24}-\left(30S_{2}(\widehat{x},\widehat{y})+5\widehat{w}K_{2}(\widehat{x},\widehat{y},\widehat{r})+2\widehat{w}^{2}L_{2}(\widehat{r})\right)=0,\\
\widehat{r}_{0}\widehat{r}_{34}-\left(30S_{3}(\widehat{x},\widehat{y})+5\widehat{w}K_{3}(\widehat{x},\widehat{y},\widehat{r})+2\widehat{w}^{2}L_{3}(\widehat{r})\right)=0,\\
\widehat{r}_{0}\widehat{r}_{35}-\left(30S_{4}(\widehat{x},\widehat{y})+5\widehat{w}K_{4}(\widehat{x},\widehat{y},\widehat{r})+2\widehat{w}^{2}L_{4}(\widehat{r})\right)=0,
\end{cases}\label{eq:RhatEq}
\end{equation}
which are obtained from the equations of $\mathscr{R}$ by modifying
with the variable $\widehat{w}$.
\end{defn}

Note that $\widehat{\mathscr{R}}$ contains the following copy of
$\mP^{9}$: 
\begin{equation}
\widehat{\Pi}:=\{\widehat{r}_{0}=1;\widehat{w}=\widehat{r}_{15}-30\widehat{S}_{1}=\widehat{r}_{24}-30\widehat{S}_{2}=\widehat{r}_{34}-30\widehat{S}_{3}=\widehat{r}_{35}-30\widehat{S}_{4}=0\}.\label{eq:Pihat}
\end{equation}
We have the morphism $\widehat{f}_{\mathscr{R}}\colon\widehat{\mathscr{R}}\to\mathscr{R}$
given by 
\begin{equation}
[\widehat{w},\widehat{r}_{0};\widehat{x},\widehat{y},\widehat{r}_{15},\widehat{r}_{24},\widehat{r}_{34},\widehat{r}_{35}]\mapsto[\widehat{r}_{0},\widehat{w}^{\frac{1}{2}}\widehat{x},\widehat{w}^{\frac{1}{2}}\widehat{y},\widehat{w}^{2}\widehat{r}_{15},\widehat{w}^{2}\widehat{r}_{24},\widehat{w}^{2}\widehat{r}_{34},\widehat{w}^{2}\widehat{r}_{35}].\label{eq:fR}
\end{equation}
Noting that the equations (\ref{eq:RhatEq}) become those of $\mathscr{R}$
by setting $\widehat{w}=1$, we can easily verify that $\widehat{f}_{\mathscr{R}}$
is an isomorphism outside $\widehat{\Pi}$ and contracts $\widehat{\Pi}$
to the $r_{0}$-point of $\mathscr{R}$. 
\begin{prop}
\label{prop:BlUp}The morphism $\widehat{f}_{\mathscr{R}}$ is the
blow-up at the $r_{0}$-point of $\mathscr{R}$ and the $\widehat{f}_{\mathscr{R}}$-exceptional
divisor is $\widehat{\Pi}$. In particular the normal bundle of $\widehat{\Pi}$
in $\widehat{\mathscr{R}}$ is $\sO_{\mP^{9}}(-2)$.
\end{prop}

\begin{proof}
Indeed, by the equation of $\widehat{\mathscr{R}}$, we see that,
$\widehat{\mathscr{R}}$ is isomorphic near $\widehat{\Pi}$ to an
open subset of $\mA^{1}\times\mP^{9}$ with coordinates $\widehat{w};\widehat{x},\widehat{y}$,
and similarly, $\mathscr{R}$ is isomorphic near the $r_{0}$-point
to the quotient of an open subset of the affine $10$-space with coordinates
$x,y$. Then, locally the map $\widehat{f}_{\mathscr{R}}$ is defined
by 
\[
(\widehat{w};\widehat{x},\widehat{y})\mapsto(\widehat{w}^{\frac{1}{2}}\widehat{x},\widehat{w}^{\frac{1}{2}}\widehat{y}).
\]

Moreover, on the $\widehat{x}_{k}$-chart ($1\leq k\leq5)$, this
coincides the map 
\[
(\widehat{w},\widehat{x}_{1},\dots,\check{\widehat{x}_{k}},\dots,\widehat{x}_{5},\widehat{y})\mapsto(\widehat{w}^{\frac{1}{2}}\widehat{x}_{1},\dots,\stackrel{k}{\check{\widehat{w}^{\frac{1}{2}}}},\dots,\widehat{w}^{\frac{1}{2}}\widehat{x}_{5},\widehat{w}^{\frac{1}{2}}\widehat{y}),
\]
where $\widehat{x}_{k}$ is deleted on the l.h.s., and $\widehat{w}^{\frac{1}{2}}$
is inserted in the $k$-th position on the r.h.s. We have also a similar
description of $\widehat{f}_{\mathscr{R}}$ on the $\widehat{y}_{l}$-chart
($1\leq k\leq5)$. Therefore this coincides the blow-up at the $r_{0}$-point. 
\end{proof}
\begin{cor}
\label{cor:The-subscheme-QGorwidehatR}The subscheme $\widehat{\mathscr{R}}$
is $\mQ$-Gorenstein.
\end{cor}

\begin{proof}
This follows from Corollary \ref{cor:RQGor} and Proposition \ref{prop:BlUp}.
\end{proof}

\subsubsection{\textbf{The flopped contraction $\widehat{g}_{\mathscr{R}}\colon\widehat{\mathscr{R}}\to\overline{\mathscr{R}}$
to $\widetilde{g}_{\mathscr{R}}$\label{subsec:The-flopped-contraction} }}

We have the morphism $\widehat{g}_{\mathscr{R}}\colon\widehat{\mathscr{R}}\to\overline{\mathscr{R}}$
given by 
\begin{equation}
[\widehat{w},\widehat{r}_{0};\widehat{x},\widehat{y},\widehat{r}_{15},\widehat{r}_{24},\widehat{r}_{34},\widehat{r}_{35}]\mapsto[\widehat{x},\widehat{y},\widehat{w}\widehat{r}_{15},\widehat{w}\widehat{r}_{24},\widehat{w}\widehat{r}_{34},\widehat{w}\widehat{r}_{35}].\label{eq:gR}
\end{equation}

\begin{rem}
As a morphism of VGIT of toric varieties, $\widehat{g}_{\mathscr{R}}$
should come also with the monomials $\widehat{r}_{0}\widehat{r}_{15},\widehat{r}_{0}\widehat{r}_{24},\widehat{r}_{0}\widehat{r}_{34},\widehat{r}_{0}\widehat{r}_{35}$
but we may erase them by the 6--9th equations in (\ref{eq:RhatEq}).
Accordingly, we obtain the weight matrix $\left(\begin{array}{ccc}
0 & 1 & 2\\
-1 & 0 & 1
\end{array}\right)$ from Table \ref{eq:wtmatRhat}, where the 1st, 2nd and 3rd column
of this matrix correspond to $\widehat{w}$, the collection of $\widehat{x}$
and $\widehat{y}$, and the collection of $\widehat{r}_{15},\widehat{r}_{24},\widehat{r}_{34},\widehat{r}_{35}$,
respectively. Subtracting twice the 2nd row from the 1st one and interchanging
the two rows, we obtain the matrix $\left(\begin{array}{ccc}
-1 & 0 & 1\\
2 & 1 & 0
\end{array}\right),$ which represents the weights of coordinates provided in Table \ref{Table:FR}.
This clarifies the relation between Tables \ref{Table:FR} and \ref{eq:wtmatRhat}.
\end{rem}

\begin{lem}
\label{lem:ghat Bir}The morphism $\widehat{g}_{\mathscr{R}}$ is
surjective and birational.
\end{lem}

\begin{proof}
To show that $\widehat{g}_{\mathscr{R}}$ is surjective, it suffices
to prove $\widehat{g}_{\mathscr{R}}$ is dominant since it is projective.
For this, we take a point $\mathsf{p}=(\overline{r}_{15},\overline{r}_{24},\overline{r}_{34},\overline{r}_{35},\overline{x},\overline{y})\in\mathscr{\overline{R}}$
with $\overline{r}_{15}\overline{r}_{24}\overline{r}_{34}\overline{r}_{35}\not=0$.
We substitute $\widehat{r}_{15},\widehat{r}_{24},\widehat{r}_{34},\widehat{r}_{35},\widehat{x},\widehat{y}$
with $\overline{r}_{15},\overline{r}_{24},\overline{r}_{34},\overline{r}_{35},\overline{x},\overline{y}$
and set $\widehat{w}=1$. Then the first 5 equations in (\ref{eq:RhatEq})
are satisfied since $\mathsf{p\in\overline{\mathscr{R}}}$. Note that
the last 4 equations are 
\begin{equation}
\widehat{r}_{0}\overline{r}_{15}-30\overline{U}_{1}=0,\,\widehat{r}_{0}\overline{r}_{24}-30\overline{U}_{2}=0,\,\widehat{r}_{0}\overline{r}_{34}-30\overline{U}_{3}=0,\,\widehat{r}_{0}\overline{r}_{35}-30\overline{U}_{4}=0,\label{eq:4eq}
\end{equation}
where we set $\overline{U}_{i}:=U_{i}(\overline{x},\overline{y},\overline{r})\,(1\leq i\leq4)$.
Now, we set only for the moment that $\overline{r}_{1}:=\overline{r}_{15}$,
$\overline{r}_{2}:=\overline{r}_{24}$, $\overline{r}_{3}:=\overline{r}_{34}$,
and $\overline{r}_{4}:=\overline{r}_{35}$. Since $\mathsf{p\in\overline{\mathscr{R}}}$,
we see that $\overline{r}_{i}\overline{U}_{j}-\overline{r}_{j}\overline{U}_{i}=0$
by the equality $C_{3}B_{3}=B_{2}A_{3}$ in (\ref{eq:comp}). Therefore
(\ref{eq:4eq}) uniquely determines the value of $\widehat{r}_{0}$,
and $[1,\widehat{r}_{0};\overline{r}_{15},\overline{r}_{24},\overline{r}_{34},\overline{r}_{35},\overline{x},\overline{y}]$
with this $\widehat{r}_{0}$ is mapped to $\mathsf{p}$. Thus $\widehat{g}_{\mathscr{R}}$
is dominant.

To show that $\widehat{g}_{\mathscr{R}}$ is birational, we again
take a point $\mathsf{p}=(\overline{r}_{15},\overline{r}_{24},\overline{r}_{34},\overline{r}_{35},\overline{x},\overline{y})\in\overline{\mathscr{R}}$
with $\overline{r}_{15}\overline{r}_{24}\overline{r}_{34}\overline{r}_{35}\not=0$.
It suffices to prove that $\widehat{g}_{\mathscr{R}}^{-1}(\mathsf{p})$
consists of only one point. Let $[\widehat{w},\widehat{r}_{0};\widehat{r}_{15},\widehat{r}_{24},\widehat{r}_{34},\widehat{r}_{35},\widehat{x},\widehat{y}]\in\widehat{g}_{\mathscr{R}}^{-1}(\mathsf{p})$.
Since $\widehat{w}\not=0$, we may assume that $\widehat{w}=1$ by
the $\mC^{*}$-action as in the 4th line of Table \ref{eq:wtmatRhat}.
Then we have $(\widehat{r}_{15},\widehat{r}_{24},\widehat{r}_{34},\widehat{r}_{35},\widehat{x},\widehat{y})=(\alpha^{2}\overline{r}_{15},\alpha^{2}\overline{r}_{24},\alpha^{2}\overline{r}_{34},\alpha^{2}\overline{r}_{35},\alpha\overline{x},\alpha\overline{y})$
for some $\alpha\not=0$. By the $\mC^{*}$-action as in the 3th line
of Table \ref{eq:wtmatRhat}, we may assume that $\alpha=1$. Now
the value of $\widehat{r}_{0}$ is uniquely determined as in the first
paragraph, we see that $\widehat{g}_{\mathscr{R}}^{-1}(\mathsf{p})$
consists of only one point.
\end{proof}
\begin{prop}
\label{prop:gR} The scheme $\widehat{\mathscr{R}}$ is a normal $\mQ$-factorial
variety with only terminal singularities and Picard number $2$. The
morphism $\widehat{g}_{\mathscr{R}}$ is the flopped contraction to
$\widetilde{g}_{\mathscr{R}}$. The $\widehat{g}_{\mathscr{R}}$-exceptional
locus coincides with 
\[
\widehat{\Lambda}=\{\widehat{r}_{15}=\widehat{r}_{24}=\widehat{r}_{34}=\widehat{r}_{35}=S_{1}(\widehat{x},\widehat{y})=S_{2}(\widehat{x},\widehat{y})=S_{3}(\widehat{x},\widehat{y})=S_{4}(\widehat{x},\widehat{y})=0\},
\]
which is the inverse image of $\overline{S}_{\mathscr{R}}$ defined
in Proposition \ref{prop:hR}. Any nontrivial $\widehat{g}_{\mathscr{R}}$-fiber
is isomorphic to $\mP^{1}$.
\end{prop}

\begin{proof}
Let $\mathsf{p}=[\overline{x},\overline{y},\overline{r}_{15},\overline{r}_{24},\overline{r}_{34},\overline{r}_{35}]$
be a point of $\overline{\mathscr{R}}$. By Lemma \ref{lem:ghat Bir}
and (\ref{eq:gR}), $\widehat{g}_{\mathscr{R}}$ is not isomorphic
over $\mathsf{p}$ if and only if $\mP^{1}(\widehat{w},\widehat{r}_{0})$
is contained in $\widehat{g}_{\mathscr{R}}^{-1}(\mathsf{p})$. By
the equations (\ref{eq:RhatEq}) of $\widehat{\mathscr{R}}$, we see
that this is equivalent to that $\overline{r}_{ij}=0,$ $\overline{S}_{1}=\overline{S}_{2}=\overline{S}_{3}=\overline{S}_{4}=0$,
namely, $\mathsf{p}\in\overline{S}_{\mathscr{R}}$. Thus $\widehat{g}_{\mathscr{R}}$
is isomorphic in codimension one since codimension of $\overline{S}_{\mathscr{R}}$
in $\overline{\mathscr{R}}$ is $3$, and any $\widehat{g}_{\mathscr{R}}$-fiber
over a point of $\overline{S}_{\mathscr{R}}$ is $\mP^{1}$. We also
see that the $\widehat{g}_{\mathscr{R}}$-exceptional locus coincides
with $\Lambda$.

This implies that $\widehat{\mathscr{R}}$ is regular in codimension
1 since so is $\mathscr{\overline{R}}$. Hence $\widehat{\mathscr{R}}$
is normal since $\widehat{\mathscr{R}}$ is $\mQ$-Gorenstein by Corollary
\ref{cor:The-subscheme-QGorwidehatR}. Moreover, since $\widehat{\mathscr{R}}\to\overline{\mathscr{R}}$
is isomorphic in codimension one and $\overline{\mathscr{R}}$ is
terminal by Proposition \ref{prop:hR} and $\widehat{\mathscr{R}}$
is normal, we conclude that $\widehat{\mathscr{R}}$ is terminal,
noting the discrepancies over $\widehat{\mathscr{R}}$ and $\overline{\mathscr{R}}$
of exceptional divisors are the same. 

Note that there are two morphisms $\widetilde{g}_{\mathscr{R}}$ and
$\widetilde{f}{}_{\mathscr{R}}$ from $\mathscr{\widetilde{R}}$ and
they are the only nontrivial morphisms from $\widetilde{\mathscr{R}}$
with connected fibers since the Picard number of $\widetilde{\mathscr{R}}$
is $2$ by Proposition \ref{prop:f'}. Therefore $\widehat{\mathscr{R}}$
cannot be isomorphic to $\widetilde{\mathscr{R}}$ by the existence
of the morphism $\widehat{f}_{\mathscr{R}}\colon\widehat{\mathscr{R}}\to\mathscr{R}$. 

Now we show that $\widehat{\mathscr{R}}$ is $\mQ$-factorial and
its Picard number is $2$. Let $\widehat{D}$ be a $\widehat{g}_{\mathscr{R}}$-ample
divisor on $\widehat{\mathscr{R}}$ and $\widetilde{D}$ its strict
transform on $\widetilde{\mathscr{R}}$. Since the relative Picard
number for $\widetilde{g}_{\mathscr{R}}$ is $1$ by Proposition \ref{prop:f'},
$\widetilde{\mathscr{R}}$ is $\mQ$-factorial by Proposition \ref{prop:QfacR},
and $\widehat{\mathscr{R}}\not\simeq\widetilde{\mathscr{R}}$ over
$\overline{\mathscr{R}}$, we see that $-\widetilde{D}$ is $\widetilde{g}_{\mathscr{R}}$-ample.
Running $\widetilde{D}$-MMP starting from $\widetilde{\mathscr{R}}$
over $\overline{\mathscr{R}}$ (it is possible by \cite{bchm}), we
get a relatively nef model $\widehat{g}'_{\mathscr{R}}\colon\widehat{\mathscr{R}}'\to\overline{\mathscr{R}}$
for $\widetilde{D}$. Since the relative Picard number for $\widetilde{g}_{\mathscr{R}}$
is $1$, so is for $\widehat{g}'_{\mathscr{R}}$, thus the strict
transform $\widehat{D}'$ on $\widehat{\mathscr{R}}'$ of $\widetilde{D}$
is actually $\widehat{g}'_{\mathscr{R}}$-ample. Thus it holds that
$\widehat{\mathscr{R}}\simeq\widehat{\mathscr{R}}'$ over $\overline{\mathscr{R}}$.
Therefore the assertion holds since $\widehat{\mathscr{R}}'$ is $\mQ$-factorial
and its Picard number is $2$.

Therefore $\widehat{g}_{\mathscr{R}}$ is the flopped contraction
to $\widetilde{g}_{\mathscr{R}}$.
\end{proof}
We denote by $\sO_{\widehat{\mathscr{R}}}(1)$ the pull-back of the
divisor $\sO_{\overline{\mathscr{R}}}(1)$. By (\ref{eq:tildeR canonicaldiv}),
we have 

\begin{equation}
-K_{\mathscr{\widehat{R}}}=\sO_{\widehat{\mathscr{R}}}(8).\label{eq:hatR canonicaldiv}
\end{equation}

\subsubsection{\textbf{The key variety $\mathscr{R}$ and its linear section\label{subsec:The-variety-}}}
\begin{prop}
\label{prop:CoseqR}The variety $\mathscr{R}$ is a terminal $\mQ$-factorial
$\mQ$-Fano variety with Picard number $1$. It holds that $-K_{\mathscr{R}}=\sO_{\mathscr{R}}(8)$,
$\deg\mathscr{R}=\sO_{\mathscr{R}}(1)^{10}=3$ and $\sO_{\mathscr{R}}(1)$
is numerically primitive.
\end{prop}

\begin{proof}
The first assertion holds by Propositions \ref{prop:BlUp} and \ref{prop:gR}.
The equality $-K_{\mathscr{R}}=\sO_{\mathscr{R}}(8)$ follows from
(\ref{eq:hatR canonicaldiv}).\textbf{ }Since $\overline{\mathscr{R}}=\overline{{\rm G}}\cap(2)^{6}$,
we have $\deg\overline{\mathscr{R}}=\sO_{\overline{\mathscr{R}}}(1)^{10}=5/2$
by (\ref{eq:degG}). Then, since $\widehat{g}_{\mathscr{R}}$ is a
flopping contraction and $\widehat{f}_{\mathscr{R}}$ is the blow-up
at a $\nicefrac{1}{2}$-singularity, we have $\deg\mathscr{R}=5/2+\nicefrac{1}{2}=3$.
Let $A$ be an ample Weil divisor on $\mathscr{R}$ such that its
numerical class generates $N^{1}(\mathscr{R}).$ Let $d$ be the integer
such that $\sO_{\mathscr{R}}(1)\equiv dA$. By $\sO_{\mathscr{R}}(1)^{10}=3$,
we have $A^{10}=3/d^{10}$. Since $2A^{10}\in\mN$, we have $d=1$. 
\end{proof}
\begin{cor}
\label{cor:numX}Let $X$ be a linear section of $\mathscr{R}$ with
only six $\nicefrac{1}{2}$-singularities. Then $X$ is a $\mQ$-Fano
$3$-fold with $(-K_{X})^{3}=3$ and $g(X)=1$ such that $-K_{X}$
is numerically primitive. 
\end{cor}

\begin{proof}
The assertions follows from Proposition \ref{prop:CoseqR} except
the assertion $g(X)=1$, which follows from the Riemann-Roch theorem
by using the fact that $(-K_{X})^{3}=3$ and $X$ has only six $\nicefrac{1}{2}$-singularities. 
\end{proof}
\begin{rem}
Actually $X$ as in Corollary \ref{cor:numX} has Picard number 1
and is of Type $R$. This is the content of Theorem \ref{thm:main}
(1-2). Unfortunately, we cannot directly show that the Picard number
of $X$ is $1$ from Proposition \ref{prop:CoseqR} since we do not
have an appropriate Lefschetz type theorem.
\end{rem}

\subsubsection{\textbf{Singularities of $\mathscr{R}$\label{subsec:Singularities-of}}}

For completeness, we determine the singularities of $\mathscr{R}$.
\begin{lem}
\label{lem:rank2 smooth} It holds that $\widehat{\mathscr{R}}$ is
smooth along the $\widehat{g}_{\mathscr{R}}$-inverse image of the
locus of points with $\rank M_{q}=2$, where the matrix $M_{q}$ is
defined in Lemma \ref{lem:Mq}.
\end{lem}

\begin{proof}
The assertion follows from Propositions \ref{prop:gR} and \ref{prop:hyp double pt},
and \cite[Prop.2.2 and Ex.2.3]{Ko}.
\end{proof}
\begin{prop}
\label{prop:SingR} The singular locus of $\mathscr{R}$ consists
of the following loci as in $(1)$--$(3):$

\begin{enumerate}[$(1)$]

\item The following $6$ points $\mathsf{\mathsf{p_{0}},p}_{1},\dots,\mathsf{p}_{5}\in\mP(r_{0},r_{15},r_{24},r_{34},r_{35})=\{x=y=\bm{o}\}:$

$\mathsf{p}_{1}:=\{r_{24}=r_{34}=r_{35}=0,r_{0}+6r_{15}=0\}$,

$\mathsf{p}_{2}:=\{r_{15}=r_{34}=r_{35}=r_{0}-6r_{24}=0\}$,

$\mathsf{p}_{3}:=\{r_{15}=r_{24}=r_{35}=r_{0}+6r_{34}=0\}$,

$\mathsf{p}_{4}:=\{r_{24}=r_{34}=r_{15}-r_{35}=r_{0}-6r_{15}=0\},$

$\mathsf{p}_{5}:=\{r_{15}=r_{34}+r_{35}=r_{24}+r_{34}=r_{0}+6r_{24}=0\}$,

$\mathsf{p}_{0}:=\{r_{15}=r_{24}=r_{34}=r_{35}=0\},$ 

where $\mathsf{p}_{i}\,(1\leq i\leq5)$ corresponds to $\widetilde{\mathsf{p}}_{i}\,(1\leq i\leq5)$
as in Proposition \ref{prop:QfacR}.

\item The following $15$ loci $\Gamma_{1},\dots,\Gamma_{5},\Delta_{1},\dots,\Delta_{10}$,
each of which is isomorphic to $\mP^{3}\!:$
\begin{align*}
\Gamma_{1}:=\{r_{24} & =r_{34}=r_{35}=r_{0}+6r_{15}=r_{15}+q_{15}\\
 & =x_{2}=x_{3}=x_{4}=y_{2}=y_{3}=y_{4}=0\}.\\
\Gamma_{2}:=\{r_{15} & =r_{34}=r_{35}=r_{0}-6r_{24}=r_{24}+q_{24}\\
 & =x_{1}=x_{3}=x_{5}=y_{1}=y_{3}=y_{5}=0\}.\\
\Gamma_{3}:=\{r_{15} & =r_{24}=r_{35}=r_{0}+6r_{34}=r_{34}+q_{34}\\
 & =x_{1}+x_{2}=x_{1}+x_{4}=x_{5}=y_{1}+y_{2}=y_{1}+y_{4}=y_{5}=0\}.\\
\Gamma_{4}:=\{r_{24} & =r_{34}=r_{15}-r_{35}=r_{0}-6r_{35}=r_{35}+q_{35}\\
 & =x_{4}=x_{1}+x_{2}=x_{1}-x_{3}+x_{5}=y_{4}=y_{1}+y_{2}=y_{1}-y_{3}+y_{5}=0\}.\\
\Gamma_{5}:=\{r_{15} & =r_{34}+r_{35}=r_{24}+r_{34}=r_{0}+6r_{24}=r_{24}+q_{24}\\
 & =x_{1}=x_{2}+x_{3}=x_{4}+x_{5}=y_{1}=y_{2}+y_{3}=y_{4}+y_{5}=0\}.
\end{align*}
\noindent Each of the $10$ loci below are contained in $\mP(x,y)=\{r_{0}=r_{15}=r_{24}=r_{34}=r_{35}=0\}\!:$
\begin{align*}
\Delta_{1}:= & \{x_{1}=x_{2}=x_{3}=y_{1}=y_{2}=y_{3}=0\}.\\
\Delta_{2}:= & \{x_{1}=x_{2}=x_{4}=y_{1}=y_{2}=y_{4}=0\}.\\
\Delta_{3}:= & \{x_{1}=x_{4}=x_{5}=y_{1}=y_{4}=y_{5}=0\}.\\
\Delta_{4}:= & \{x_{3}=x_{5}=x_{2}-x_{4}=y_{3}=y_{5}=y_{2}-y_{4}=0\}.\\
\Delta_{5}:= & \{x_{1}=x_{2}-x_{4}=x_{3}-x_{5}=y_{1}=y_{2}-y_{4}=y_{3}-y_{5}=0\}.\\
\Delta_{6}:= & \{x_{3}+x_{4}=x_{2}-x_{4}=x_{1}+x_{4}+x_{5}=y_{3}+y_{4}=y_{2}-y_{4}=y_{1}+y_{4}+y_{5}=0\}.\\
\Delta_{7}:= & \{x_{3}=x_{4}=x_{1}+x_{5}=y_{3}=y_{4}=y_{1}+y_{5}=0\}.\\
\Delta_{8}:= & \{x_{4}=x_{5}=x_{2}+x_{3}=y_{4}=y_{5}=y_{2}+y_{3}=0\}.\\
\Delta_{9}:= & \{x_{3}=x_{1}+x_{2}=x_{1}+x_{4}+x_{5}=y_{3}=y_{1}+y_{2}=y_{1}+y_{4}+y_{5}=0\}.\\
\Delta_{10}:= & \{x_{5}=x_{2}+x_{3}=x_{1}-x_{3}=y_{5}=y_{2}+y_{3}=y_{1}-y_{3}=0\}.
\end{align*}
(note that the loci $\Gamma_{i}$ correspond to $\widetilde{\Gamma}_{i}\,(1\leq i\leq5)$
as in the proof of Proposition \ref{prop:QfacR}, and the loci $\Delta_{j}$
correspond to the closure of $\overline{\Delta}_{j}^{o}\,(1\leq j\leq10)$
as in the proof of Lemma \ref{lem:Mq}).

These $15$ loci are transformed to each other by the $\mathfrak{S}_{6}$-action
on $\mathscr{R}$.

\item The locus $\Delta_{0}:=\{q_{ij}=0\,(1\leq i<j\leq5)\}\subset\mP(x,y)$.

\end{enumerate}

The variety $\mathscr{R}$ has $\nicefrac{1}{2}$-singularities at
the $6$ points as in $(1)$, and $c({\rm G}(2,5))$-singularities
along the $15$ loci as in $(2)$ outside the locus $\Delta_{0}$.
\end{prop}

\begin{proof}
We can directly check that the $15$ loci as in (2) are transformed
to each other by the $\mathfrak{S}_{6}$-action on $\mathscr{R}$. 

By Corollary \ref{cor:RQGor}, $\mathscr{R}$ has an $\nicefrac{1}{2}$-singularity
at $\mathsf{p}_{0}:=$the $r_{0}$-point. Thus, by Proposition \ref{prop:BlUp},
it suffices to show a similar statement for $\widehat{\mathscr{R}}$
along $\{\widehat{w}=1\}$. We denote by $\widehat{\Gamma}_{1},\dots,\widehat{\Gamma}_{5}$,
and $\widehat{\Delta}_{0},\widehat{\Delta}_{1},\dots,\widehat{\Delta}_{10}\subset\widehat{\mathscr{R}}$
the strict transforms of $\Gamma_{1},\dots,\Gamma_{5}$, and $\Delta_{0},\Delta_{1},\dots,\Delta_{10}$,
respectively. Let us determine the singularity of $\widehat{\mathscr{R}}$
separately on the $\widehat{g}_{\mathscr{R}}$-exceptional locus $\widehat{\Lambda}$
(cf.$\,$Proposition \ref{prop:gR}), and on $\widehat{\mathscr{R}}\setminus\widehat{\Lambda}$. 

First we study the singularities of $\widehat{\mathscr{R}}\setminus\widehat{\Lambda}$.
Since $\widehat{\mathscr{R}}\dashrightarrow\widetilde{\mathscr{R}}$
is isomorphic outside $\widehat{\Lambda}$, the description of the
singularities of $\widehat{\mathscr{R}}\setminus\widehat{\Lambda}$
is the same as that of the singularities of $\widetilde{\mathscr{R}}$
outside the inverse image of $\overline{S}_{\mathscr{R}}$. Since
$\overline{\Delta}_{0},\overline{\Delta}_{1},\dots,\overline{\Delta}_{10}\subset\Pi$
is contained in the locus of the points with $\rank M_{q}\leq1$ (Proposition
\ref{lem:Mq}), we see that the singular locus of $\widehat{\mathscr{R}}\setminus\widehat{\Lambda}$
consists of $(\cup_{i=1}^{5}\widehat{\Gamma}_{i})\setminus\widehat{\Lambda}$
and the $5$ points corresponding to $\mathsf{p}_{1},\dots,\mathsf{p_{5}}$,
and the description of the singularities along them is the desired
one by Proposition \ref{prop:QfacR} and its proof.

Now we show that the variety $\mathscr{R}$ has $c({\rm G}(2,5))$-singularities
along the $15$ loci as in $(2)$ outside the locus $\Delta_{0}$.
In fact, we have only to show this along $\Gamma_{1}$ by the result
of the first paragraph. By the third paragraph, we see that $\widehat{\mathscr{R}}$
has $c({\rm G}(2,5))$-singularities along $\widehat{\Gamma}_{1}\setminus\widehat{\Lambda}$.
Note that $\widehat{\Gamma}_{1}\cap\widehat{\Lambda}$ is contained
in $\widehat{\Delta}_{0}$. Therefore $\mathscr{R}$ has $c({\rm G}(2,5))$-singularities
along $\Gamma_{1}\setminus\Delta_{0}$ as desired. 

Finally we check the singularities of $\widehat{\mathscr{R}}$ in
$\widehat{\Lambda}$. Since we can check that $\widehat{\Lambda}\cap(\cup_{i=1}^{5}\widehat{\Gamma}_{i})$
is contained in $\widehat{\Delta}_{0}$, we have only to show that
the singular locus of $\widehat{\mathscr{R}}$ contained in $\widehat{\Lambda}$
coincides with $\widehat{\Delta}_{0}\cup(\cup_{i=1}^{10}\widehat{\Delta}_{i})$.
By Lemma \ref{lem:rank2 smooth}, the singular locus of $\widehat{\mathscr{R}}$
contained in $\widehat{\Lambda}$ is actually contained in the $\widehat{g}_{\mathscr{R}}$-inverse
image of the locus of the points with $\rank M_{q}\leq1$. Therefore,
by the $\mathfrak{S}_{5}$-action on $\widehat{\mathscr{R}}$ and
$\overline{\mathscr{R}}$, we have only to consider the problem over
$\overline{\Delta}_{0}$ and one of the $10$ loci $\overline{\Delta}_{1}^{o},\dots,\overline{\Delta}_{10}^{o}$
as in (\ref{eq:10loci}), say, over $\overline{\Delta}_{1}^{o}$.
The $\widehat{g}_{\mathscr{R}}$-inverse image of the closure of $\overline{\Delta}_{1}^{o}$
is 
\begin{equation}
\{\widehat{r}_{15}=\widehat{r}_{24}=\widehat{r}_{34}=\widehat{r}_{35}=\widehat{x}_{1}=\widehat{x}_{2}=\widehat{x}_{3}=\widehat{y}_{1}=\widehat{y}_{2}=\widehat{y}_{3}=0\}.\label{eq:Delta1}
\end{equation}
To determine the singular locus of $\widehat{\mathscr{R}}$ along
this locus, we return to the consideration on $\mathscr{R}$; it is
easier to determine the singular locus contained in 
\begin{equation}
\{r_{15}=r_{24}=r_{34}=r_{35}=x_{1}=x_{2}=x_{3}=y_{1}=y_{2}=y_{3}=0\},\label{eq:Delta1V2}
\end{equation}
which is the image of the locus (\ref{eq:Delta1}) on $\mathscr{R}$
and contains $\Delta_{1}$. Computing the linear part of the equations
of $\mathscr{R}$ at points of the locus (\ref{eq:Delta1V2}), we
see that $\mathscr{R}$ is nonsingular along $r_{0}\not=0$ near the
locus (\ref{eq:Delta1V2}). Therefore, near the locus (\ref{eq:Delta1V2}),
$\mathscr{R}$ is singular only along $\Delta_{1}$. Similarly, we
can conclude that $\mathscr{R}$ is singular in the locus $\{r_{15}=r_{24}=r_{34}=r_{35}=0,q_{ij}=0\,(1\leq i<j\leq5)\}$
only along $\Delta_{0}$ and $\mathsf{p}_{0}$. 
\end{proof}
\begin{rem}
\label{rem:P2P2Fib}We can also describe the singularities along the
locus as in (3). Here we only give a sketch of the description and
publish details elsewhere. For this, let $\mP^{9}$ be the projective
space with coordinates $y,r_{0},r_{15},r_{24},r_{34},r_{35}$ and
$\mA^{5}$ the affine space with coordinates $x$ (we can do the same
investigation by exchanging $x$ and $y$). We consider the quasi-projective
variety $\mathfrak{R}$ in $\mP^{9}\times\mA^{5}$ with the same equations
as $\mathscr{R}$. Let 

\begin{align*}
D_{R}(x):= & x_{1}^{2}x{}_{2}^{2}+2x{}_{1}^{2}x_{2}x_{3}+x{}_{1}^{2}x{}_{3}^{2}-2x{}_{1}^{2}x_{2}x_{4}-2x{}_{1}^{2}x_{3}x_{4}\\
 & +2x_{1}x_{2}x_{3}x_{4}+2x_{1}x{}_{3}^{2}x_{4}+x{}_{1}^{2}x{}_{4}^{2}-2x_{1}x_{3}x{}_{4}^{2}\\
 & +x_{3}^{2}x{}_{4}^{2}+2x_{1}x{}_{2}^{2}x_{5}+2x_{1}x_{2}x_{3}x_{5}-2x_{1}x_{2}x_{4}x_{5}\\
 & -4x_{1}x_{3}x_{4}x_{5}-2x_{2}x_{3}x_{4}x_{5}+x{}_{2}^{2}x{}_{5}^{2},
\end{align*}
which is an $\mathfrak{S}_{6}$-invariant quartic polynomial. Then
we can show that the natural projection $\mathfrak{R}\to\mA^{5}$
is $\mathfrak{S}_{6}$-equivariant, and the fiber of this over any
point outside $\{D_{R}=0\}$ is isomorphic to the 5-dimensional cone
over $\mP^{2}\times\mP^{2}.$ The proof of this fact is similar to
the one of \cite[Prop.4.5]{Taka5}. The key idea is to define the
quadratic Jordan algebra of a cubics form from the equations of $\mathfrak{R}$
(cf.{[}ibid., $\S$3{]}). Using this fact, we can describe the singularities
along the locus as in (3) in a similar way to {[}ibid., Prop.4.11{]}.
\end{rem}

\section{\textbf{Proof of Theorem \ref{thm:main} (1)\label{subsec:Proof-of-Theorem R}}}

\subsection{Preliminary results}

\begin{con}\label{con:Hi}

Let
\begin{align*}
\mathsf{H}_{i} & :=\{\sum_{j=1}^{5}a_{ij}\mathsf{x}_{j}+\sum_{j=1}^{5}b_{ij}\mathsf{y}_{j}=0\}\subset\mP(V^{*})\times\mP(V^{\oplus2}),\\
\widetilde{H}_{i} & :=\{\sum_{j=1}^{5}a_{ij}\widetilde{x}_{j}+\sum_{j=1}^{5}b_{ij}\widetilde{y}_{j}=0\}\subset\mF_{\widetilde{\mathscr{R}}}\quad(1\leq i\leq7),
\end{align*}
where $a_{ij},b_{ij}\in\mC$. By the same equations with changing
the coordinates suitably, we can also define the divisors corresponding
to $\mathsf{H}_{i}$ in the ambient space $\mP(1^{10},2^{4})$ of
$\overline{\mathscr{R}}$, the toric variety $\mF_{\widehat{\mathscr{R}}}$
and the ambient space $\mP(1^{10},2^{5})$ of $\mathscr{R}$, which
we denote by $\overline{H}_{i}$, $\widehat{H}_{i}$ and $H_{i}$,
respectively.

Associated with these, we give names to objects in the restriction
of the diagram (\ref{eq:Sarkisov}): we define
\begin{align*}
\mathsf{X}_{H} & :=\mathsf{R}\cap\bigcap_{i=1}^{7}\mathsf{H}_{i},\,\widetilde{X}_{H}:=\widetilde{\mathscr{R}}\cap\bigcap_{i=1}^{7}\widetilde{H}_{i},\,\overline{X}_{H}:=\overline{\mathscr{R}}\cap\bigcap_{i=1}^{7}\overline{H}_{i},\\
\widehat{X}_{H} & :=\widehat{\mathscr{R}}\cap\bigcap_{i=1}^{7}\widehat{H}_{i},\,and\,\,X_{H}:=\mathscr{R}\cap\bigcap_{i=1}^{7}H_{i},
\end{align*}
and let 
\[
\widetilde{f}_{H}\colon\widetilde{X}_{H}\to L_{\mathscr{R}},\,\widetilde{g}_{H}\colon\widetilde{X}_{H}\to\overline{X}_{H},\,\widehat{g}_{H}\colon\widehat{X}_{H}\to\overline{X}_{H},\,and\,\,\widehat{f}_{H}\colon\widehat{X}_{H}\to X_{H}
\]
 be the restrictions of $\widetilde{f}_{\mathscr{R}}$, $\widetilde{g}_{\mathscr{R}}$,
$\widehat{g}_{\mathscr{R}}$, and $\widehat{f}_{\mathscr{R}}$, respectively.
\end{con}

In the following two lemmas, we compare the two schemes $\widetilde{X}_{H}$
and $\mathsf{X}_{H}$. By Proposition \ref{prop:Rtilde and P(Omega2)},
these are isomorphic outside the loci $(\widetilde{f}{}_{\mathscr{R}})^{-1}(\cup_{1\leq i<j\leq5}\,\ell_{ij})$
and $\varphi^{-1}(\cup_{1\leq i\leq5}\,\mathsf{P}_{i})$, respectively.
Moreover, by the same proposition, we see that $\widetilde{X}_{H}\to L_{\mathscr{R}}$
is isomorphic near the fiber over a point of $\ell_{ij}\setminus\{\mathsf{q}_{i},\mathsf{q_{j}}\}$
if and only of $\mathsf{X}_{H}\to{\rm S}_{3}$ is isomorphic near
the fiber over $\mathsf{q}_{ij}$. 

The main point is to compare the schemes $\widetilde{X}_{H}$ and
$\mathsf{X}_{H}$ on the $\widetilde{f}_{H}$-fiber over the point
$\mathsf{q}_{i}$ and the inverse image of the plane $\mathsf{P}_{i}$
by $\mathsf{X}_{H}\to{\rm S}_{3}$ respectively (cf.~Proposition
\ref{prop:Segre}). We investigate the relation between $\bigcap_{i=1}^{7}\widetilde{H}_{i}$
and $\bigcap_{i=1}^{7}\mathsf{H}_{i}$ over $\mathsf{q}_{i}$ and
$\mathsf{P}_{i}$ respectively. For this, using the $\mathfrak{S}_{5}$-action
on $L_{\mathscr{R}},$ we may assume that $\mathsf{q}_{i}$ is $\mathsf{q}_{1}=$the
$\widetilde{r}_{15}$-point. Then, by an explicit calculation, we
see that $\mathsf{P}_{i}=\{\mathsf{z}_{1}=\mathsf{z}_{5}=0\}$. Then
the fiber of $\widetilde{{\rm \mathscr{R}}}\to L_{\mathscr{R}}$ over
$\mathsf{q}_{1}$ is 
\[
\mathsf{C}_{\mathsf{q}_{1}}:=\left\{ \rank\left(\begin{array}{ccc}
\widetilde{x}_{2} & \widetilde{x}_{3} & \widetilde{x}_{4}\\
\widetilde{y}_{2} & \widetilde{y}_{3} & \widetilde{y}_{4}
\end{array}\right)\leq1\right\} \subset\mP(\widetilde{x},\widetilde{y},w),
\]
and $S_{\mathsf{P_{1}}}:=\mP(\Omega_{\mP(V^{*})}(1)^{\oplus2}|_{\mathsf{P}_{1}})$
coincides with 
\[
\left\{ \mathsf{z}_{2}\mathsf{x}_{2}+\mathsf{z}_{3}\mathsf{x}_{3}+\mathsf{z}_{4}\mathsf{x}_{4}=0,\mathsf{z}_{2}\mathsf{y}_{2}+\mathsf{z}_{3}\mathsf{y}_{3}+\mathsf{z}_{4}\mathsf{y}_{4}=0\right\} \subset\mathsf{P}_{1}\times\mP(V^{\oplus2}).
\]

\begin{lem}
\label{lem:FibCorresp}The following are equivalent:

\begin{enumerate}[$(1)$]

\item $\mathsf{C}_{\mathsf{q}_{1}}\cap\bigcap_{i=1}^{7}\widetilde{H}_{i}$
is $1$-dimensional and $\rank M_{15}=4$, where we set 
\[
M_{15}:=\left(\begin{array}{cccc}
a_{i1} & a_{i5} & b_{i1} & b_{i5}\end{array}\,(1\leq i\leq7)\right).
\]

\item $S_{\mathsf{P}_{1}}\cap\bigcap_{i=1}^{7}\mathsf{H}_{i}$ is
irreducible and 2-dimensional.

\end{enumerate}
\end{lem}

\begin{proof}
Let $\mathsf{C}_{\mathsf{q}_{1}}'$ be the variety in $\mP(\widetilde{x},\widetilde{y})$
obtained from $\mathsf{C}_{\mathsf{q}_{1}}$ by the projection from
the $w$-point, and $\widetilde{H}_{i}'$ the hypersurface in $\mP(\widetilde{x},\widetilde{y})$
with the same equation as $\widetilde{H}_{i}$. It is easy to verify
that (1) is equivalent to the following:

\[
(1)'\ \mathsf{C}_{\mathsf{q}_{1}}'\cap\bigcap_{i=1}^{7}\widetilde{H}_{i}'\ \text{is}\ 0\text{-dimensional and}\ \rank M_{15}=4.
\]
We prove the equivalence between (1)' and (2). 

Supposing (2) holds, we show that $\rank M_{15}=4$. Indeed, if $\rank M_{15}<4$,
then changing the coordinates of $x_{1},x_{5},y_{1},y_{5}$ if necessary,
we may assume that $b_{i5}=0\,(1\leq i\leq7)$. Then the condition
{\footnotesize
\begin{equation}
\rank\left(\begin{array}{cccccccccc}
z_{2} & z_{3} & z_{4} & 0 & 0 & 0 & 0 & 0 & 0 & 0\\
0 & 0 & 0 & z_{2} & z_{3} & z_{4} & 0 & 0 & 0 & 0\\
a_{12} & a_{13} & a_{14} & b_{12} & b_{13} & b_{14} & a_{11} & a_{15} & b_{11} & 0\\
\vdots & \vdots & \vdots & \vdots & \vdots & \vdots & \vdots & \vdots & \vdots & \vdots\\
a_{72} & a_{73} & a_{74} & b_{72} & b_{73} & b_{74} & a_{71} & a_{75} & b_{71} & 0
\end{array}\right)\leq8\label{eq:9times10 mat}
\end{equation}
}is equivalent to the vanishing of the determinant $D$ of $9\times9$
matrix obtained by deleting the $10$ th column from the matrix as
in (\ref{eq:9times10 mat}). Thus the fiber of $S_{\mathsf{P}_{1}}\cap\bigcap_{i=1}^{7}\mathsf{H}_{i}\to\mathsf{P}_{1}$
is positive dimensional over the curve $\{D=0\}$ in $\mathbf{\mathsf{P}}_{1}$.
This implies that $\dim(S_{\mathsf{P_{1}}}\cap\bigcap_{i=1}^{7}\mathsf{H}_{i})\geq3$
or $S_{\mathsf{P_{1}}}\cap\bigcap_{i=1}^{7}\mathsf{H}_{i}$ is reducible,
a contradiction to the assumption of (2). Thus $\rank M_{15}=4$. 

Hence, to prove the equivalence of (1)' and (2), we may assume that
$\rank M_{15}=4$. By row operations, we may further assume that $M_{15}=\left(\begin{array}{c}
O\\
I
\end{array}\right),$ where $I$ is the $4\times4$ identity matrix and $O$ is the zero
matrix (note that we can change each of $\mathsf{H}_{i}$ and $\widetilde{H}_{i}'$
without changing $\bigcap_{i=1}^{7}\mathsf{H}_{i}$ and $\bigcap_{i=1}^{7}\widetilde{H}_{i}'$).
Then the conditions (1)' and (2) are equivalent to the following conditions
(i) and (ii) respectively:

\begin{enumerate}[(i)]

\item
\begin{align*}
\left\{ \rank\left(\begin{array}{ccc}
\widetilde{x}_{2} & \widetilde{x}_{3} & \widetilde{x}_{4}\\
\widetilde{y}_{2} & \widetilde{y}_{3} & \widetilde{y}_{4}
\end{array}\right)\leq1\right\} \cap\bigcap_{i=1}^{3}\left\{ \sum_{j=2}^{4}a_{ij}\widetilde{x}_{j}+\sum_{j=2}^{4}\widetilde{b}_{ij}y_{j}=0\,(1\leq i\leq3)\right\} \\
\subset\mP(\widetilde{x}_{2},\widetilde{x}_{3},\widetilde{x}_{4},\widetilde{y}_{2},\widetilde{y}_{3},\widetilde{y}_{4})\simeq\mP^{5}
\end{align*}
is $0$-dimensional. 

\item 
\begin{equation}
\mP(\Omega_{\mathsf{P_{1}}}(1)^{\oplus2})\cap\bigcap_{i=1}^{3}\left\{ \sum_{j=2}^{4}a_{ij}\mathsf{x}_{j}+\sum_{j=2}^{4}b_{ij}\mathsf{y}_{j}=0\,(1\leq i\leq3)\right\} \label{eq:Omegasection}
\end{equation}
 is irreducible and 2-dimensional (note that (\ref{eq:Omegasection})
coincides with $S_{\mathsf{P}_{1}}\cap\bigcap_{i=1}^{7}\mathsf{H}_{i}$).

\end{enumerate}The condition (i) (resp.$\,$(ii)) is equivalent to
Lemma \ref{lem:AL} (c) (resp.$\,$(a)) below. Therefore this lemma
follows from Lemma \ref{lem:AL}.
\end{proof}
\begin{lem}
\label{lem:AL}Let $L$ be a $3$-dimensional linear space and $A$
a $2$-dimensional linear space. Let $p_{1}\colon\mP(A\otimes\Omega_{\mP(L)}(1))\to\mP(L)$
and $p_{2}\colon\mP(A\otimes\Omega_{\mP(L)}(1))\to\mP(A\otimes L^{*})$
be the natural morphisms. It holds that $p_{2}$ is the blow-up of
$\mP(A\otimes L^{*})$ along the Segre variety $\mP(A)\times\mP(L^{*})$
(note that $\mP(A)\times\mP(L^{*})$ is defined by the $2\times2$
minors of a $2\times3$ matrix by choosing coordinates of $\mP(A\otimes L^{*})$
as in the proof). Let $M$ be a $3$-dimensional subspace of $A\otimes L^{*}$.
The following are equivalent:

\begin{enumerate}[(a)]

\item $p_{2}^{-1}(\mP(M))$ is irreducible and $2$-dimensional.

\item $p_{2}^{-1}(\mP(M))\to\mP(M)$ has a finite number of positive
dimensional fibers.

\item $(\mP(A)\times\mP(L^{*}))\cap\mP(M)$ is $0$-dimensional.

\item $p_{2}^{-1}(\mP(M))\to\mP(L)$ is birational and has a finite
number of positive dimensional fibers and they are $1$-dimensional.

\end{enumerate}
\end{lem}

\begin{proof}
Let $z_{1},z_{2},z_{3}$ be coordinates of $L$ and $x_{1},x_{2},x_{3},y_{1},y_{2},y_{3}$
are coordinates of $A\otimes L^{*}$, where $\text{\ensuremath{x_{1}},\ensuremath{x_{2}},\ensuremath{x_{3}}}$
and $y_{1},y_{2},y_{3}$ are coordinates of copies of $L^{*}$. Then
\[
\mP(A\otimes\Omega_{\mP(L)}(1))=\left\{ \sum_{i=1}^{3}x_{i}z_{i}=0,\sum_{i=1}^{3}y_{i}z_{i}=0\right\} \subset\mP(L)\times\mP(A\otimes L^{*}).
\]
A local computation shows that $p_{2}$ is actually the blow-up of
$\mP(A\otimes L^{*})$ along $\mP(A)\times\mP(L^{*})$. 

$(a)\Rightarrow(b)$ This follows since $p_{2}$ is the blow-up of
$\mP(A\otimes L^{*})$ along $\mP(A)\times\mP(L^{*})$.

$(b)\Rightarrow(a)$ Note that a positive dimensional fiber of $p_{2}^{-1}(\mP(M))\to\mP(M)$
is $1$-dimensional since $p_{2}$ is the blow-up of $\mP(A\otimes L^{*})$
along $\mP(A)\times\mP(L^{*})$. Therefore, by $(b)$, $p_{2}^{-1}(\mP(M))$
is $2$-dimensional and then is a complete intersection in $\mP(A\otimes\Omega_{\mP(L)}(1))$.
Hence $p_{2}^{-1}(\mP(M))$ is purely 2-dimensinal. By $(b)$ again,
$p_{2}^{-1}(M)$ has only one 2-dimensional component. Thus $p_{2}^{-1}(M)$
is irreducible.

The equivalence $(b)\Leftrightarrow(c)$ follows since $p_{2}$ is
the blow-up of $\mP(A\otimes L^{*})$ along the Segre variety $\mP(A)\times\mP(L^{*})$.

$(a)\Rightarrow(d)$ This follows since any fiber of $p_{2}^{-1}(\mP(M))\to\mP(L)$
is linear.

$(d)\Rightarrow(a)$ The proof is identical to that of $(b)\Rightarrow(a)$.
\end{proof}
\begin{lem}
\label{lem:unique}Let $V$ be a $5$-dimensional vector space, and
${\rm S}_{3}$ and $T$ a Segre cubic and a cubic scroll in $\mP(V^{*})$
respectively. Suppose that $C_{{\rm S}}:={\rm S}_{3}\cap T$ is a
smooth curve of genus $7$ and degree $9$. Let $C$ be the birational
image of $C_{{\rm S}}$ on $L_{\mathscr{R}}$ by the birational map
$\mu_{{\rm S}}\colon L_{\mathscr{R}}\dashrightarrow{\rm S}_{3}$ as
in Subsection \ref{subsec:LR and S3}. We take \textbf{$\bigcap_{i=1}^{7}\mathsf{H}_{i}$}
for $T$ as in Proposition \ref{prop:The-cubic-scroll zero} and follow
Convention \ref{con:Hi}. 

\begin{enumerate}[$(A)$]

\item The following hold:

\begin{enumerate}[$(1)$]

\item $\widetilde{X}_{H}$ is irreducible and normal,

\item $-K_{\widetilde{X}_{H}}$ is $\mQ$-Cartier and relatively
ample over $L_{\mathscr{R}}$, and 

\item the morphism $\widetilde{f}_{H}\colon\widetilde{X}_{H}\to L_{\mathscr{R}}$
has only $1$-dimensional non-trivial fibers and coincides with the
blow-up along $C$ outside the fibers over $\mathsf{q}_{1},\dots,\mathsf{q}_{5}$. 

\end{enumerate}

\item If $\widetilde{X}$ and $\widetilde{X}\to L_{\mathscr{R}}$
satisfies the same properties as $\widetilde{X}_{H}$ and $\widetilde{X}_{H}\to L_{\mathscr{R}}$
stated in $(A)$ $(1)$--$(3),$ then $\widetilde{X}$ and $\widetilde{X}_{H}$
are isomorphic over $L_{\mathscr{R}}$.

\end{enumerate}
\end{lem}

\begin{proof}
(A) (3). By Proposition \ref{prop:The-cubic-scroll zero}, the naturally
induced morphism $\mathsf{X}_{H}\to{\rm S}_{3}$ is the blow-up along
$C_{{\rm S}}$. Therefore $\widetilde{f}_{H}\colon\widetilde{X}_{H}\to L_{\mathscr{R}}$
is the blow-up along $C$ possibly outside $\widetilde{f}_{H}^{-1}(\cup_{1\leq i<j\leq5}\ell_{ij})$
by Proposition \ref{prop:Rtilde and P(Omega2)}. Let us check that
$\widetilde{f}_{H}$ is an isomorphism near the fiber over any point
of $\ell_{ij}\setminus \{\mathsf{q}_{i},\mathsf{q_{j}}\}$. If $\widetilde{f}_{H}$
is not an isomorphism near the fiber over a point of $\ell_{ij}\setminus\{\mathsf{q}_{i},\mathsf{q_{j}}\}$,
then $\mathsf{X}_{H}\to{\rm S}_{3}$ is not isomorphic near the fiber
over $\mathsf{q}_{ij}$ by Proposition \ref{prop:Rtilde and P(Omega2)}.
This implies that $\mathsf{q}_{ij}\in C_{S}$, which is a contradiction
since we can show that ${\rm S}_{3}$ is smooth near $C_{{\rm S}}$
in the same way as the proof of Proposition \ref{prop:KeyLemma} (A)
(3). We have shown that $\widetilde{f}_{H}$ is the blow-up of $L_{\mathscr{R}}$
along $C$ outside the fibers over $\mathsf{q}_{1},\dots\mathsf{q}_{5}$.
It remains to show that the fibers of $\widetilde{f}_{H}$ over $\mathsf{q}_{1},\dots\mathsf{q}_{5}$
are $1$-dimensional. We show this applying Lemma \ref{lem:FibCorresp}
with interpreting $\mathsf{q}_{1}$ in there as $\mathsf{q}_{i}\,(1\leq i\leq5)$.
Note that the fiber of $\widetilde{f}_{H}$ over $\mathsf{q}_{i}$
is nothing but $\mathsf{C}_{\mathsf{q}_{1}}\cap\bigcap_{i=1}^{7}\widetilde{H}_{i}$
as in the lemma. Since $\mathsf{X}_{H}\to{\rm S}_{3}$ is the blow-up
along $C_{{\rm S}}$, the condition (2) of the lemma holds since $\mathsf{P}_{i}\cap C_{{\rm S}}$
consists of a finite number of points (cf.~Lemma \ref{lem:AL} (d)).
Therefore, the fiber of $\widetilde{f}_{H}$ is 1-dimensional over
$\mathsf{q}_{i}$ by the lemma. 

(1). Now we show that $\widetilde{X}_{H}$ is irreducible and normal.
First we note that the codimension of $\widetilde{X}_{H}$ in $\widetilde{\mathscr{R}}$
is $7$, which is equal to the number of divisors $\widetilde{H}_{i}$.
Since $\widetilde{\mathscr{R}}$ is quasi-smooth at non-Gorenstein
points by Proposition \ref{prop:QfacR}, it holds that $\widetilde{X}_{H}$
is Cohen-Macaulay at non-Gorenstein points. Since $\widetilde{\mathscr{R}}$
is Cohen-Macaulay by the same proposition, and $\widetilde{H}_{i}$
define Cartier divisors at Gorenstein points of $\widetilde{\mathscr{R}},$
$\widetilde{X}_{H}$ is Gorenstein at Gorenstein points of $\widetilde{\mathscr{R}}$.
Therefore $\widetilde{X}_{H}$ is Cohen-Macaulay. Thus, by unmixedness
and (3), $\widetilde{X}_{H}$ is purely $3$-dimensional, hence is
irreducible. Note that $\widetilde{X}_{H}$ is regular in codimension
one by (3). Therefore $\widetilde{X}_{H}$ is normal. 

(2). Since $\widetilde{H}_{i}$ are tautological divisors and $-K_{\widetilde{\mathscr{R}}}=\sO_{\mathscr{\widetilde{R}}}(8)$,
we have $-K_{\widetilde{X}_{H}}=\sO_{\mathscr{\widetilde{R}}}(1)|_{\widetilde{X}_{H}}$,
which is $\mQ$-Cartier and is relatively ample over $L_{\mathscr{R}}$. 

\vspace{3pt}

\noindent (B) The assertion follows from \cite[Lem.5.5]{Taka3} since
$\widetilde{X}_{H}$ and $\widetilde{X}$ are isomorphic in codimension
$1$.
\end{proof}
\begin{lem}
\textbf{\label{lem:RestRToX}} We take divisors $\widetilde{H}_{i}\subset\widetilde{\mathscr{R}}\,(1\leq i\leq7)$
as in Convention \ref{con:Hi}, and follows the convention. Assume
that \begin{enumerate}[$(i)$]

\item $\widetilde{X}{}_{H}$ is an irreducible normal $3$-fold with
only five $\nicefrac{1}{2}$-singularities, 

\item $\widetilde{X}_{H}$ has Picard number $2$, and 

\item the birational morphism $\widetilde{g}_{H}\colon\widetilde{X}{}_{H}\to\overline{X}_{H}$
is a flopping contraction.

\end{enumerate}

Then $X_{H}$ is a Type $R$ $\mQ$-Fano $3$-fold for the point $\widehat{f}_{\mathscr{R}}(\widehat{\Pi})$,
and the Sarkisov diagram (\ref{eq:Sarkisov}) restricts to (\ref{eq:3-foldSarkisov}).
\end{lem}

\begin{proof}
We recall the situation of Proposition \ref{prop:hR} and its proof,
and use the notation there. Since the $\widetilde{g}_{H}$-exceptional
locus consists of a finite number of $\mP^{1}$'s by (iii), it holds
that $\overline{S}_{\mathscr{R}}\cap\bigcap_{i=1}^{7}\overline{H}_{i}$
consist of a finite number of points in the rank 2 locus of the matrix
$M_{q}$ by the proposition. Therefore the $\widehat{g}_{H}$-exceptional
locus also consists of a finite number of $\mP^{1}$'s by Proposition
\ref{prop:gR}. 

We can show that $\widehat{X}_{H}$ is irreducible and normal in the
same way as proving that $\widetilde{X}{}_{H}$ is irreducible and
normal in the proof of Lemma \ref{lem:unique}. In the same way to
show Proposition \ref{prop:gR}, we see that $\widehat{X}_{H}\dashrightarrow\widetilde{X}{}_{H}$
is the flop for the flopping contraction $\widetilde{g}_{H}$. In
particular, $\widehat{X}_{H}$ has Picard number $2$ and only five
$\nicefrac{1}{2}$-singularities by (i), (ii) and \cite{Ko}. Since
$\widehat{E}:=\widehat{\Pi}\cap\widehat{X}_{H}$ is $\mP^{2}$ with
normal bundle $\sO_{\mP^{2}}(-2)$ by Proposition \ref{prop:BlUp},
it holds that the restriction $\widehat{f}_{H}\colon\widehat{X}_{H}\to X_{H}$
of $\widehat{f}_{\mathscr{R}}$ to $\widehat{X}_{H}$ is the blow-up
at a $\nicefrac{1}{2}$-singularity. Therefore $X_{H}$ has only six
$\nicefrac{1}{2}$-singularities and Picard number 1. Taking Corollary
\ref{cor:numX} into account, it remains to show that $\widetilde{f}_{H}$
has the desired properties. For simplicity of notation, we denote
$\widetilde{X}_{H}$ by $\widetilde{X}$ and $L_{\mathscr{R}}\simeq\mP^{3}$
by $Y$. Since $\dim\widetilde{X}=\dim Y$, $\widetilde{f}_{H}$ is
birational. Let $\widetilde{E}$ be the $\widetilde{f}_{H}$-exceptional
divisor. If $\widetilde{f}_{H}$ contracts $\widetilde{E}$ to a point,
then $\widetilde{f}_{H}$ is the weighted blow-up at a smooth point
with some weight by \cite{Ka1}. Then, however, $\widetilde{X}$ cannot
have five $\nicefrac{1}{2}$-singularities, a contradiction. Thus
$\widetilde{f}_{H}$ contracts $\widetilde{E}$ to a curve, which
we denote by $C$. Let $\widetilde{L}$ be the $\widetilde{f}_{H}$-pull-back
of $\sO_{\mP^{3}}(1)$ and $\widehat{E}_{\widetilde{X}}:=\widetilde{\Pi}\cap\widetilde{X}$.
Note that\textbf{ }$\widehat{E}_{\widetilde{X}}$ is the strict transform
of $\widehat{E}$. Since $\widehat{f}_{H}$ is the blow-up at a $\nicefrac{1}{2}$-singularity,
we see that $(-K_{\widehat{X}_{H}})^{3}=\nicefrac{5}{2}$, $(-K_{\widehat{X}_{H}})^{2}\widehat{E}=1$
and $(-K_{\widehat{X}_{H}})\widehat{E}^{2}=-2$ by a standard computation.
Using a property of the flop $\widehat{X}_{H}\dashrightarrow\widetilde{X}$,
we obtain that 
\begin{equation}
(-K_{\widetilde{X}})^{3}=\nicefrac{5}{2},\,(-K_{\widetilde{X}})^{2}\widehat{E}_{\widetilde{X}}=1\ \text{and}\ (-K_{\widetilde{X}})\widehat{E}_{\widetilde{X}}^{2}=-2\label{eq:KE}
\end{equation}
 (cf.~\cite[Part I, Lem.3.1]{Taka2}). By Lemma \ref{lem:2H-L R},
it holds that\textbf{ 
\begin{equation}
\widetilde{L}=2(-K_{\widetilde{X}})-\widehat{E}_{\widetilde{X}}.\label{L=00003D2-K-E}
\end{equation}
}Since $\widetilde{L}^{3}=1$, we have
\[
1=\left(2(-K_{\widetilde{X}})-\widehat{E}_{\widetilde{X}}\right)^{3}=-4-\widehat{E}_{\widetilde{X}}^{3},\text{i.e.},\widehat{E}_{\widetilde{X}}^{3}=-5
\]
by (\ref{eq:KE}). By (\ref{L=00003D2-K-E}) and $-K_{\widetilde{X}}=4\widetilde{L}-\widetilde{E}$,
we have $\widetilde{E}=7(-K_{\widetilde{X}})-4\widehat{E}_{\widetilde{X}}$,
and then, by (\ref{eq:KE}), we obtain
\begin{align*}
 & -\widetilde{E}^{2}\widetilde{L}=-\left(7(-K_{\widetilde{X}})-4\widehat{E}_{\widetilde{X}}\right){}^{2}\left(2(-K_{\widetilde{X}})-\widehat{E}_{\widetilde{X}}\right)=12,
\end{align*}
which implies that $\deg C=12$. We also have
\[
(-K_{\widetilde{E}})^{2}=(-K_{\widetilde{X}}-\widetilde{E})^{2}\widetilde{E}=\left(6(-K_{\widetilde{X}})-4\widehat{E}_{\widetilde{X}}\right){}^{2}\left(7(-K_{\widetilde{X}})-4\widehat{E}_{\widetilde{X}}\right)=-138,
\]
which is equal to $8\left(1-p_{g}(C)\right)-2m-18m'$ as in \cite[Prop.7.1]{Taka3},
where $m$ and $m'$ are the numbers of fibers of types (A) and (B)
in there, respectively. Since $Y=\mP^{3}$ is smooth, all the $\nicefrac{1}{2}$-singularities
are contained in $\widetilde{E}$, and it holds that $m=0$ and $m'=5$,
and $p_{g}(C)=7$. Since $m'=5$, $C$ has $5$ triple points at $\mathsf{q}_{1},\dots,\mathsf{q}_{5}$
which can be resolved by the blow-up at them. Therefore, $X$ is of
Type $R$ for for the point $\widehat{f}_{\mathscr{R}}(\widehat{\Pi})$.
\end{proof}

\subsection{Proof of Theorem \ref{thm:main} (1-1)}

Let $X$ be a Type $R$ $\mQ$-Fano $3$-fold. We show that $X$ is
a linear section of $\mathscr{R}$ of codimension $7$. 

We consider the diagram (\ref{eq:3-foldSarkisov}). We identify $Y\simeq\mP^{3}$
with $L_{\mathscr{R}}$. We define $C_{{\rm S}}$ as the birational
model of $C$ on ${\rm S}_{3}$ by the birational map $\mu_{{\rm S}}\colon L_{\mathscr{R}}\dashrightarrow{\rm S}_{3}$
as in Proposition \ref{prop:KeyLemma}. Then, by the proposition,
we may write $C_{{\rm S}}={\rm S}_{3}\cap T$, where $T$ is a cubic
scroll in $\mP(V^{*})$. We take \textbf{$\bigcap_{i=1}^{7}\mathsf{H}_{i}$}
for $T$ as in Proposition \ref{prop:The-cubic-scroll zero}, and
follow Convention \ref{con:Hi}. By Lemma \ref{lem:unique} (B), $\widetilde{X}$
as in (\ref{eq:3-foldSarkisov}) is isomorphic to $\widetilde{X}_{H}$
over $L_{\mathscr{R}}$. Hence we have $\widetilde{X}=\widetilde{\mathscr{R}}\cap\bigcap_{i=1}^{7}\widetilde{H}_{i}$.
By the assumptions of Theorem \ref{thm:main} (1-1), $\widetilde{X}$
satisfies the assumptions of Lemma \ref{lem:RestRToX}. Since we can
go backwards along the Sarkisov diagram uniquely, $X$ coincides with
$X_{H}$ constructed in Lemma \ref{lem:RestRToX}. Thus we have shown
Theorem \ref{thm:main} (1-1). $\hfill\square$

\subsection{Proof of Theorem \ref{thm:main} (1-2)}
\begin{lem}
\label{lem:Construct-example}There exists a linear section of $\mathscr{R}$
of codimension $7$ that is a Type $R$ $\mQ$-Fano $3$-fold.
\end{lem}

\begin{proof}
We construct such an example starting from the cubic scroll $T\subset\mP^{4}$
as in Proposition \ref{prop:GenExOf f'} (A). We use the notation
as in there. Since $T$ as in the proposition satisfies the assumptions
as in Lemma \ref{lem:unique}, the assertion (A) of the lemma holds.
Since $\widetilde{X}$ and the morphism $\widetilde{f}\colon\widetilde{X}\to\mP^{3}$
constructed in the proposition satisfies the assumption of the assertion
(B) of the lemma, we may write $\widetilde{X}=\widetilde{\mathscr{R}}\cap\bigcap_{i=1}^{7}\widetilde{H}_{i}$
as in the lemma. We follow Convention \ref{con:Hi} but we omit the
subscript $H$ from the symbols for simplicity. An important point
to note is that the construction in the proof of the proposition shows
that the Picard number of $\widetilde{X}$ is $2.$ By the explicit
construction as in the proposition, we are ready to apply Lemma \ref{lem:RestRToX}
once we check the following:

\vspace{5pt}

\noindent Claim. \textit{It holds that $\widetilde{g}\colon\widetilde{X}\to\overline{X}$
is a flopping contraction.}

\vspace{5pt}

\noindent \textit{Proof.} Let $\widetilde{L}$ be the $\widetilde{f}_{H}$-pull-back
of $\sO_{\mP^{3}}(1)$. Suppose for contradiction that $\widetilde{g}$
is not a flopping contraction. Then it is a crepant divisorial contraction.
Let $\widetilde{G}$ be the $\widetilde{g}$-exceptional divisor.
Since the Picard number of $\widetilde{X}$ is $2$, and $-K_{\widetilde{X}}$
and $\widetilde{L}$ are numerically independent, we may write $\widetilde{G}\equiv p(-K_{\widetilde{X}})-q\widetilde{L}$
with some $p,q\in\mQ$. Since $(-K_{\widetilde{X}})^{3}=5/2$ by Proposition
\ref{prop:GenExOf f'} (B) (4), and $(-K_{\widetilde{X}})^{2}\widetilde{G}=0$,
we have 
\begin{equation}
0=5/2\,p-q(-K_{\widetilde{X}})^{2}\widetilde{L}.\label{eq:pq}
\end{equation}
We note that $\widetilde{E}\widetilde{L}^{2}=0$ and $\widetilde{E}^{2}\widetilde{L}=-\deg C=-12$.
Then, since $-K_{\widetilde{X}}=4\widetilde{L}-\widetilde{E}$, we
obtain

\begin{align}
 & (-K_{\widetilde{X}})^{2}\widetilde{E}+2(-K_{\widetilde{X}})\widetilde{E}^{2}+\widetilde{E}^{3}=0,\label{eq:K2E+2KE2+E3}\\
 & (-K_{\widetilde{X}})\widetilde{E}^{2}+\widetilde{E}^{3}=-48.\label{KE2+E3}
\end{align}
Since $(-K_{\widetilde{E}})^{2}=(-K_{\widetilde{X}}-\widetilde{E})^{2}\widetilde{E}=-138$
by Proposition \ref{prop:GenExOf f'} (B) (4), we have 
\begin{equation}
(-K_{\widetilde{X}})^{2}\widetilde{E}-2(-K_{\widetilde{E}})\widetilde{E}^{2}+\widetilde{E}^{3}=-138.\label{eq:KwE-2KE^2+E3}
\end{equation}
Therefore, by (\ref{eq:K2E+2KE2+E3}), (\ref{KE2+E3}) and (\ref{eq:KwE-2KE^2+E3}),
we obtain $(-K_{\widetilde{X}})^{2}\widetilde{E}=\nicefrac{27}{2}$.
Since $-K_{\widetilde{X}}=4\widetilde{L}-\widetilde{E}$, we have
$(-K_{\widetilde{X}})^{2}\widetilde{L}=\nicefrac{1}{4}\left((-K_{\widetilde{X}})^{3}+(-K_{\widetilde{X}})^{2}\widetilde{E}\right)=4$.
From this and (\ref{eq:pq}), we have 
\begin{equation}
\widetilde{G}\equiv p\left((-K_{\widetilde{X}})-\nicefrac{5}{8}\widetilde{L}\right).\label{eq:Gpq}
\end{equation}
 Pushing this down to $\mP^{3}$, we have $\widetilde{f}{}_{*}\widetilde{G}\equiv\nicefrac{27p}{8}\widetilde{L}.$
Thus $m:=\nicefrac{27p}{8}$ is a natural number. If $\widetilde{g}(\widetilde{G})$
is a curve, then we have $\widetilde{G}\cdot\gamma=-2$ for a general
$\widetilde{g}|_{\widetilde{G}}$-fiber $\gamma$. Thus, by (\ref{eq:Gpq}),
it holds that 
\[
-2=-\nicefrac{5p}{8}(\widetilde{L}\cdot\gamma)=-\nicefrac{5m}{27}(\widetilde{L}\cdot\gamma),
\]
which has no integral solutions for $m$ and $\widetilde{L}\cdot\gamma$,
a contradiction. If $\widetilde{g}(\widetilde{G})$ is a point, then,
by the proof of Proposition \ref{prop:hR}, $\widetilde{G}$ is $\mP^{2}$
or a quadric surface. Let $\gamma$ be a line on $\widetilde{G}$.
Then $\widetilde{G}\cdot\gamma=K_{\widetilde{G}}\cdot\gamma=-3$ or
$-2$. Then we derive a contradiction as in the case that $\widetilde{g}(\widetilde{G})$
is a curve.$\hfill\square$

Therefore, by Lemma \ref{lem:RestRToX}, we obtain a Type $R$ $\mQ$-Fano
$3$-fold $X_{H}$ as desired.
\end{proof}
\noindent \textit{Proof of Theorem} \ref{thm:main} \textit{(1-2).}
Let $X$ be a linear section of $\mathscr{R}$ of codimension $7$
with only six $\nicefrac{1}{2}$-singularities. By \cite[Thm.4.1]{To}
and Lemma \ref{lem:Construct-example}, $X$ has Picard number 1.
Then, by Corollary \ref{cor:numX}, $X$ is a prime $\mQ$-Fano $3$-fold
of $g(X)=1$. 

It remains to show that $X$ is of Type $R$ for any $\nicefrac{1}{2}$-singularity
$p$. By the $\mathfrak{S}_{6}$-action on $\mathscr{R}$, we can
construct the diagram (\ref{eq:Sarkisov}) for $p\in\mathscr{R}$.
We write $X=\mathscr{R}\cap\bigcap_{i=1}^{7}H_{i}$ and as in Convention
\ref{con:Hi}, we define the divisors $\overline{H}_{i}\subset\mP(1^{10},2^{4})$,
etc.~corresponding to $H_{i}$. We also define $\overline{X}_{H}$,
$\widehat{f}_{H}$, etc., as in Convention \ref{con:Hi} but we omit
the subscript $H$ from the symbols for simplicity. By the local computation
for $\widehat{f}_{\mathscr{R}}$ using (\ref{eq:fR}), we have $\widehat{H}_{i}=\widehat{f}_{\mathscr{R}}^{*}H_{i}-\nicefrac{1}{2}\widehat{\Sigma}$.
On the other hand, since $X$ has a $\nicefrac{1}{2}$-singularity
at $p$, we see that $\widehat{H}_{i}$ are the strict transforms
of $H_{i}$. Therefore $\widehat{f}\colon\widehat{X}\to X$ is the
blow-up at $p$. In particular, the Picard number of $\widehat{X}$
is $2$. By the local computation for $\widehat{g}_{\mathscr{R}}$
using (\ref{eq:gR}), we have $\widehat{H}_{i}=\widehat{g}_{\mathscr{R}}^{*}H_{i}$.
Since $\widehat{g}_{\mathscr{R}}$ is a flopping contraction by Proposition
\ref{prop:gR}, we have $-K_{\widehat{\mathscr{R}}}=\widehat{g}_{\mathscr{R}}^{*}(-K_{\overline{\mathscr{R}}})$.
This implies that $-K_{\widehat{X}}=\widehat{g}^{*}(-K_{\overline{X}})$.

\vspace{5pt}

\noindent Claim. \textit{It holds that $\widehat{g}\colon\widehat{X}\to\overline{X}$
is a flopping contraction.}

\vspace{5pt}

\noindent \textit{Proof.} Suppose for contradiction that $\widehat{g}$
is not a flopping contraction. Then it is a crepant divisorial contraction.
Let $\widehat{G}$ be the $\widehat{g}$-exceptional divisor. Since
the Picard number of $\widehat{X}$ is $2$ and $-K_{\widehat{X}}$
and $\widehat{E}$ are numerically independent, we may write $\widehat{G}\equiv p(-K_{\widehat{X}})-q\widehat{E}$
with some $p,q\in\mQ$. Since $(-K_{\widehat{X}})^{3}=5/2$, $(-K_{\widehat{X}})^{2}\widehat{E}=1$,
and $(-K_{\widehat{X}})^{2}\widehat{G}=0$, we have 
\begin{equation}
\widehat{G}\equiv p\left((-K_{\widehat{X}})-5/2\,\widehat{E}\right).\label{eq:Gpq-1}
\end{equation}
 Pushing this down to $X$, we have $\widehat{f}_{*}\widehat{G}\equiv p(-K_{\widehat{X}}).$
Thus it holds that $p\in\mN$ since $-K_{\widehat{X}}$ is numerically
primitive. By Proposition \ref{prop:gR}, $\widehat{g}(\widehat{G})$
is a curve. Hence we have $\widehat{G}\cdot\gamma=-2$ for a general
non-trivial $\widehat{g}$-fiber $\gamma$. Thus, by (\ref{eq:Gpq-1}),
it holds that 
\[
-2=-5p/2\,(\widehat{E}\cdot\gamma).
\]
This has no integral solutions for $p$ and $(\widehat{E}\cdot\gamma)$,
a contradiction.$\hfill\square$

\vspace{5pt}

By the claim, $\widehat{g}$-exceptional locus consists of a finite
number of $\mP^{1}$'s. Note that $\mathscr{R}$ has Gorenstein singularities
along $\Delta_{0}\cup\bigcup_{i=1}^{10}\Delta_{i}$ by Proposition
\ref{prop:SingR} and it dominates the rank $\leq1$ locus of the
matrix $M_{q}$ in $\overline{S}_{\mathscr{R}}$ as in Lemma \ref{lem:Mq}.
Therefore, since $\widehat{X}$ has only $\nicefrac{1}{2}$-singularities,
it holds that $\overline{S}_{\mathscr{R}}\cap\bigcap_{i=1}^{7}\overline{H}_{i}$
consist of a finite number of points in the rank 2 locus of the matrix
$M_{q}$. Therefore the $\widetilde{g}$-exceptional locus also consists
of a finite number of $\mP^{1}$'s by the proof of Proposition \ref{prop:hR}.
By the proof of Proposition \ref{prop:gR}, we see that $\widehat{X}\dashrightarrow\widetilde{X}$
is the flop for $\widehat{g}$. Therefore the conditions (i)--(iii)
of Lemma \ref{lem:RestRToX} are satisfied, and then it follows that
$X$ is of Type $R$ for $p$. 

$\hfill\square$

\section{\textbf{The key variety $\mathscr{IR}$ for Type $I\!R$ $\mQ$-Fano
$3$-folds\label{sec:TypeIRGen}}}

In this section, we construct the key variety $\mathscr{IR}$ for
a Type $I\!R$ $\mQ$-Fano $3$-fold $X$. We construct it via the
Sarkisov link as in the case of Type $R$ (Subsection \ref{subsec:IR}).
We often divide the discussions into General and Special Cases, which
correspond to the case of\textbf{ $\rank\mathcal{Q}_{B}=6$} and the
case of\textbf{ $\rank\mathcal{Q}_{B}=5$} as in Subsection \ref{subsec:TypeIR},
respectively. 

\subsection{Birational geometry of the cone $\overline{{\rm G}}_{\mathscr{I}\mathscr{R}}$
over ${\rm G}(2,5)$\label{subsec:Birational-geometry-of IR Goverline}}

Before the construction of $\mathscr{IR}$, we study in detail birational
geometry of the suitable cone $\overline{{\rm G}}_{\mathscr{IR}}$
over ${\rm G}(2,5)$ in Subsection \ref{subsec:Birational-geometry-of IR Goverline}.
Among other things, we define a birational morphism ${\rm \widetilde{G}}_{\mathscr{IR}}\to\overline{{\rm G}}_{\mathscr{IR}}$
and investigate its properties. 

\subsubsection{\textbf{The setting\label{subsec:The-settingIR Gen}}}

Let $V$ be a $5$-dimensional vector space with coordinates $\overline{p}:=(\overline{p}_{1},\dots,\overline{p}_{5})$,
and $W$ a $4$-dimensional vector space with coordinates $\overline{x}:=(\overline{x}_{1},\dots,\overline{x}_{4})$.
We consider the cone $\overline{{\rm G}}$ over ${\rm G}(2,5)$ as
in Section \ref{sec:ConeG(2,5)} under this situation. We denote the
cone by $\overline{{\rm G}}_{\mathscr{IR}}$. In both General and
Special Cases, we need a smooth quadric $3$-fold $Q$ in $\mP(V)$.
Explicitly, we may assume that 
\[
Q=\{\overline{p}_{1}\overline{p}_{2}+\overline{p}_{3}\overline{p}_{4}+\overline{p}_{5}^{2}=0\}\subset\mP(V).
\]

\vspace{3pt}

\noindent\textbf{General Case.} It is well-known that lines in $Q$
are parameterized by the doubly Veronese embedded $\mP^{3}=\mP(W)$
in ${\rm G}(2,V),$ where $W$ is a $4$-dimensional vector space.
By an explicit calculation, we see that the double Veronese embedding
is defined by 

\begin{align}
\mP(W) & \ni[\overline{x}]\mapsto\nonumber \\
 & \left[\begin{array}{ccccc}
0 & -\overline{x}_{1}\overline{x}_{4}-\overline{x}_{2}\overline{x}_{3} & -\overline{x}_{3}^{2} & \overline{x}_{4}^{2} & \overline{x}_{3}\overline{x}_{4}\\
 & 0 & \overline{x}_{1}^{2} & -\overline{x}_{2}^{2} & \overline{x}_{1}\overline{x}_{2}\\
 &  & 0 & \overline{x}_{1}\overline{x}_{4}-\overline{x}_{2}\overline{x}_{3} & \overline{x}_{1}\overline{x}_{3}\\
 &  &  & 0 & -\overline{x}_{2}\overline{x}_{4}\\
 &  &  &  & 0
\end{array}\right]\in{\rm G}(2,V).\label{eq:twomat}
\end{align}
We denote this embedding by $q\colon\mP(W)\to{\rm G}(2,V)$. We denote
by $q_{ij}(\overline{x})\,(1\leq i<j\leq5)$ the $(i,j)$-entries
of the skew-symmetric matrix in (\ref{eq:twomat}). For $W$ and $q_{ij}(\overline{x})$
defined here, we consider the constructions as in Section \ref{sec:ConeG(2,5)}. 

\vspace{10pt}

\noindent\textbf{Special Case.} As in General Case, let $\mP(W)$
be the projective $3$-space parameterizing lines contained in $Q$. 

We consider the morphism $q\colon\mP(W)\to{\rm G}(2,V)$ defined by 

\begin{align}
\mP(W) & \ni[\overline{x}]\mapsto\nonumber \\
 & \left[\begin{array}{ccccc}
0 & -\overline{x}_{1}\overline{x}_{4}-\overline{x}_{2}\overline{x}_{3} & -\overline{x}_{3}^{2} & \overline{x}_{4}^{2} & 0\\
 & 0 & \overline{x}_{1}^{2} & -\overline{x}_{2}^{2} & 0\\
 &  & 0 & \overline{x}_{1}\overline{x}_{4}-\overline{x}_{2}\overline{x}_{3} & 0\\
 &  &  & 0 & 0\\
 &  &  &  & 0
\end{array}\right]\in{\rm G}(2,V).\label{eq:twomat-1}
\end{align}
We denote by $q_{ij}(\overline{x})\,(1\leq i<j\leq5)$ the $(i,j)$-entries
of the skew-symmetric matrix in (\ref{eq:twomat-1}). For $W$ and
$q_{ij}(\overline{x})$ defined here, we consider the constructions
in Section \ref{sec:ConeG(2,5)}.

Let us clarify the meaning of the map $q$ as in (\ref{eq:twomat-1}).
We set $V':=\{\overline{p}_{5}=0\}\subset V$. Note that the $q$-image
of $\mP(W)$ in ${\rm G}(2,V)$ is contained in ${\rm G}(2,V').$
Let 
\[
\delta\colon Q\to\mP(V')
\]
be the projection from the point $(0:0:0:0:1)$ ($\not\in Q$). We
see that $\delta$ is the double cover of $\mP(V')$ branched along
the quadric surface 
\[
Q'=\{\overline{p}_{1}\overline{p}_{2}+\overline{p}_{3}\overline{p}_{4}=0\}\subset\mP(V').
\]

\begin{prop}
\label{prop:q(P(W))} The following assertions hold:

\begin{enumerate}[$(1)$]

\item $q(\mP(W))$ is the parameter space of tangent lines of $Q'$
and lines contained in $Q'$.

\item $q(\mP(W))$ is the intersection between ${\rm G}(2,V')\simeq{\rm G}(2,4)$
and the quadric $\{(p_{12}+p_{34})^{2}=4p_{13}p_{24}\}$, and has
ordinary double points along the two disjoint conics $q_{1}$ and
$q_{2}$ in $q(\mP(W))$ parameterizing lines on $Q'$. 

\item The map $q\colon\mP(W)\to q(\mP(W))$ is a double cover branched
along $q_{1}\cup q_{2}$.

\end{enumerate}
\end{prop}

\begin{proof}
(1). Comparing (\ref{eq:twomat-1}) with (\ref{eq:twomat}), the $q$-image
of $\mP(W)$ in ${\rm G}(2,V)$ parameterizes the lines in $\mP(V')$
which are the $\delta$-images of lines in the quadric $3$-fold $Q$.
Therefore a point of $q(\mP(W))$ corresponds to a tangent line of
$Q'$ or a line contained in $Q'$. Conversely, we see that a tangent
line of $Q'$ or a line contained in $Q'$ corresponds to a point
of $q(\mP(W))$ by a direct computation.

The proof of (2) and (3) also follows from direct computations. 
\end{proof}
\begin{rem}
\label{rem:overlineGto PU}We present here some remarks that may serve
as a guide for what follows. 

Let $V_{4}\subset V$ be a $4$-dimensional subspace. 

\noindent\textbf{General Case.} If the quadric surface $Q\cap\mP(V_{4})$
is smooth (resp.~singular), then $q^{-1}({\rm G}(2,V_{4}))$ consists
of the two lines in $\mP(W)$ (resp.~the line in $\mP(W)$) parameterizing
lines in $Q\cap\mP(V_{4})$. This observation and the consideration
as in Subsection \ref{subsec:GtoPwedgeV} indicate that the double
cover of $\mP(V^{*})$ branched along $Q^{*}$ can be identified with
${\rm G}(2,W)$, and a rational map $\overline{{\rm G}}_{\mathscr{I}\mathscr{R}}\dashrightarrow\mP(\sU)$
is induced, where we recall that $\sU$ is the universal subbundle
of rank $2$ on ${\rm G(2,W)}$.

\vspace{3pt}

\noindent\textbf{Special Case. }Assume that $V_{4}\not=V'$. Using
Proposition \ref{prop:q(P(W))}, we see that, if $Q'\cap\mP(V_{4})$
is smooth (resp.~singular), $q^{-1}({\rm G}(2,V_{4}))$ consists
of the two lines in $\mP(W)$ (resp.~the one line in $\mP(W)$) parameterizing
lines in the hyperplane section $\delta^{-1}(\mP(V_{4})\cap\mP(V'))$
of $Q$. This observation and the consideration as in Subsection \ref{subsec:GtoPwedgeV}
indicate that a map $\widetilde{(Q')^{*}}\setminus\{\widetilde{v}\}\to{\rm G}(2,W)$
is induced and is constant on rulings of $\widetilde{(Q')^{*}}\setminus\{\widetilde{v}\}$,
and a map $\overline{{\rm G}}_{\mathscr{I}\mathscr{R}}\dashrightarrow\mP(\sU')$
is also induced, where $\sU'$ is the pull-back by $\widetilde{(Q')^{*}}\setminus\{\widetilde{v}\}\to{\rm G}(2,W)$
of $\sU$ (note that $\sU'$ defined here can be identified with $\sU'$
defined in (\ref{eq:U'}) by regarding $\widetilde{(Q')^{*}}$ with
$\mathcal{Q}_{B}$). 
\end{rem}

\subsubsection{\textbf{The double cover of $\mP(V^{*})$\label{subsec:The-double-cover G(2,W)} }}

We define the double cover of $\mP(V^{*}),$ which is a quadric $4$-fold.
We relate it with ${\rm G}(2,W)\subset\mP(\wedge^{2}W)$ and clarify
the meaning of the relation. Let $\mathsf{m}_{ij}\,(1\leq i<j\leq4)$
be the Pl\"ucker coordinates of $\mP(\wedge^{2}W)$. 

\vspace{3pt}

\noindent\textbf{General Case.} Let $V_{4}\subset V$ be a $4$-dimensional
subspace of $V$. We consider $\mP(V_{4})$ as a point of $\mP(V^{*})$.
Let $Q^{*}\subset\mP(V^{*})$ be the quadric $3$-fold which is dual
to $Q$. By an explicit calculation, we see that 
\[
Q^{*}=\{4\mathsf{z}_{1}\mathsf{z}_{2}+4\mathsf{z}_{3}\mathsf{z}_{4}+\mathsf{z}_{5}^{2}=0\},
\]
where $\mathsf{z}_{1},\dots,\mathsf{z}_{5}$ are the dual coordinates
of $\mP(V^{*})$. We denote the double cover of $\mP(V^{*})$ branched
along $Q^{*}$ by 
\[
\{\mathsf{s}^{2}=4(\mathsf{z}_{1}\mathsf{z}_{2}+\mathsf{z}_{3}\mathsf{z}_{4}+\mathsf{z}_{5}^{2})\}.
\]
 Then we identify the double cover of $\mP(V^{*})$ with 
\[
{\rm G}(2,W)=\{\mathsf{m}_{12}\mathsf{m}_{34}-\mathsf{m}_{13}\mathsf{m}_{24}+\mathsf{m}_{14}\mathsf{m}_{23}=0\}
\]
by setting
\begin{equation}
\begin{cases}
\mathsf{m}_{12}=\mathsf{z}_{1},\mathsf{m}_{13}=\mathsf{z}_{4},\mathsf{m}_{14}=-\nicefrac{1}{2}(\mathsf{z}_{5}-\mathsf{s}),\\
\mathsf{m}_{23}=\nicefrac{1}{2}(\mathsf{z}_{5}+\mathsf{s}),\mathsf{m}_{24}=\mathsf{z}_{3},\mathsf{m}_{34}=-\mathsf{z}_{2}.
\end{cases}\label{eq:doublecoverG2W}
\end{equation}

\vspace{10pt}

\noindent\textbf{Special Case. }Let $(Q')^{*}\subset\mP(V^{*})$
be the subvariety parameterizing hyperplanes $\mP(V_{4})$ in $\mP(V)$
such that $Q'\cap\mP(V_{4})$ are singular conics or $V_{4}=V'$.
By an explicit calculation, we see that 
\[
(Q')^{*}=\{\mathsf{z}_{1}\mathsf{z}_{2}+\mathsf{z}_{3}\mathsf{z}_{4}=0\}\subset\mP(V^{*}),
\]
where $\mathsf{z}_{1},\dots,\mathsf{z}_{5}$ are the dual coordinates.
This is the cone over $\mP^{1}\times\mP^{1}$ and $V'$ corresponds
to the $\mathsf{z}_{5}$-point, which is the vertex $v$ of $(Q')^{*}$.
Let $\widetilde{(Q')^{*}}$ be the double cover of $\mP(V^{*})$ branched
along $(Q')^{*}$ and $\widetilde{v}$ the inverse image of $v$ on
$\widetilde{(Q')^{*}}$. It holds that $\widetilde{(Q')^{*}}$ is
a quadric $4$-fold of rank $5$ with $\widetilde{v}$ being the vertex.
We denote $\widetilde{(Q')^{*}}$ by 
\[
\{\mathsf{s}^{2}=4(\mathsf{z}_{1}\mathsf{z}_{2}+\mathsf{z}_{3}\mathsf{z}_{4})\}\subset\mP^{5},
\]
where $\widetilde{v}$ is the $\mathsf{z}_{5}$-point. The morphism
$\widetilde{(Q')^{*}}\setminus\{\widetilde{v}\}\to{\rm G}(2,W)$ is
defined by setting
\begin{equation}
\mathsf{m}_{12}=\mathsf{z}_{1},\mathsf{m}_{13}=\mathsf{z}_{4},\mathsf{m}_{14}=\mathsf{m}_{23}=\nicefrac{\mathsf{s}}{2},\mathsf{m}_{24}=\mathsf{z}_{3},\mathsf{m}_{34}=-\mathsf{z}_{2},\label{eq:mz}
\end{equation}
 which is constant on rulings of $\widetilde{(Q')^{*}}\setminus\{\widetilde{v}\}$,
and whose image coincides with
\[
\{\mathsf{m}_{12}\mathsf{m}_{34}-\mathsf{m}_{13}\mathsf{m}_{24}+\mathsf{m}_{23}^{2}=0\}={\rm G}(2,W)\cap\{\mathsf{m}_{14}=\mathsf{m}_{23}\}.
\]

\subsubsection{\textbf{The double cover $\widetilde{{\rm dP}}_{2}\to\mP(\wedge^{2}V)$\label{subsec:The-double-cover dp2} }}

We define the double cover $\widetilde{{\rm dP}}_{2}\to\mP(\wedge^{2}V)$,
which is compatible through the map $\mu\colon\mP(\wedge^{2}V)\dashrightarrow\mP(V^{*})$
with the double cover of $\mP(V^{*})$ defined as in Subsection \ref{subsec:The-double-cover G(2,W)}. 

\vspace{3pt}

\noindent\textbf{General Case.} We set 
\begin{equation}
\begin{cases}
m_{12}(r)={\rm Pf}_{2345}(r),m_{13}(r)=-{\rm Pf}_{1235}(r),m_{14}(r,s)=-\nicefrac{1}{2}({\rm Pf}_{1234}(r)-s),\\
m_{23}(r,s)=\nicefrac{1}{2}({\rm Pf}_{1234}(r)+s),m_{24}(r)={\rm Pf}_{1245}(r),m_{34}(r)={\rm Pf}_{1345}(r)
\end{cases}\label{eq:mijGen}
\end{equation}
(we sometimes denote $m_{ij}(r)$ by $m_{ij}(r,s)$ even when $(i,j)\not=(1,4),$
(2,3)).
\begin{defn}
We consider the weighted projective space $\mP(1^{10},2)$ with weight
one coordinates $\widetilde{r}_{ij}\,(1\leq i<j\leq5)$ and a weight
two coordinate $\widetilde{s}$. We define the double cover $\widetilde{{\rm dP}}_{2}$
of $\mP(\wedge^{2}V)$ as follows:
\[
\widetilde{{\rm dP}}_{2}:=\{m_{12}(\widetilde{r})m_{34}(\widetilde{r})-m_{13}(\widetilde{r})m_{24}(\widetilde{r})+m_{14}(\widetilde{r},\widetilde{s})m_{23}(\widetilde{r},\widetilde{s})=0\}\subset\mP(1^{10},2).
\]
The branched locus of $\widetilde{{\rm dP}}_{2}\to\mP(\wedge^{2}V)$
is the quartic hypersurface
\[
B_{{\rm G}}:=\left\{ 4\left({\rm Pf}_{2345}(\widetilde{r}){\rm Pf}_{1345}(\widetilde{r})+{\rm Pf}_{1235}(\widetilde{r}){\rm Pf}_{1245}(\widetilde{r})\right)-{\rm Pf}_{1234}(\widetilde{r})^{2}=0\right\} \subset\mP(\wedge^{2}V).
\]
\end{defn}

Note that ${\rm G}(2,V)=\{\mathrm{Pf}_{\bullet}(\widetilde{r})=0\}$
is contained in $B_{{\rm G}},$ so we also consider ${\rm G}(2,V)\subset\widetilde{{\rm dP}}_{2}$.

The above construction shows that $\mu\colon\mP(\wedge^{2}V)\dashrightarrow\mP(V^{*})$
described as in (\ref{eq:zPf}) induces the rational map $\widetilde{{\rm dP}}_{2}\dashrightarrow{\rm G}(2,W)$
defined by (\ref{eq:doublecoverG2W}) and (\ref{eq:mijGen}). 

\vspace{10pt}

\noindent\textbf{Special Case. }We set
\begin{equation}
\begin{cases}
m_{12}(r)={\rm Pf}_{2345}(r),m_{13}(r)=-{\rm Pf}_{1235}(r),m_{23}(s)=\nicefrac{s}{2},\\
m_{24}(r)={\rm Pf}_{1245}(r),m_{34}(r)={\rm Pf}_{1345}(r),\mathsf{z}_{5}(r)={\rm Pf}_{1234}(r)
\end{cases}\label{eq:mijSp}
\end{equation}
(we sometimes denote $m_{ij}(r)\,((i,j)\not=(2,3))$ by $m_{ij}(r,s)$,
$m_{23}(s)$ by $m_{23}(r,s)$ and $\mathsf{z}_{5}(r)$ by $\mathsf{z}_{5}(r,s)$).
\begin{defn}
We consider the weighted projective space $\mP(1^{10},2)$ with weight
one coordinates $\widetilde{r}_{ij}\,(1\leq i<j\leq5)$ and a weight
two coordinate $\widetilde{s}$. We define the double cover $\widetilde{{\rm dP}}_{2}$
of $\mP(\wedge^{2}V)$ as follows using (\ref{eq:mijSp}):
\begin{equation}
\widetilde{{\rm dP}}_{2}:=\{m_{12}(\widetilde{r})m_{34}(\widetilde{r})-m_{13}(\widetilde{r})m_{24}(\widetilde{r})+m_{23}(\widetilde{s})^{2}=0\}\subset\mP(1^{10},2).\label{eq:tildedp2 Sp}
\end{equation}
The branched locus of $\widetilde{{\rm dP}}_{2}\to\mP(\wedge^{2}V)$
is the quartic hypersurface
\[
B_{{\rm G}}:=\left\{ {\rm Pf}_{2345}(\widetilde{r}){\rm Pf}_{1345}(\widetilde{r})+{\rm Pf}_{1235}(\widetilde{r}){\rm Pf}_{1245}(\widetilde{r})=0\right\} .
\]
\end{defn}

Note that ${\rm G}(2,V)=\{\mathrm{Pf}_{\bullet}(\widetilde{r})=0\}$
is contained in $B_{{\rm G}},$ so we also consider ${\rm G}(2,V)\subset\widetilde{{\rm dP}}_{2}$.

The above construction shows that $\mu\colon\mP(\wedge^{2}V)\dashrightarrow\mP(V^{*})$
induces $\widetilde{{\rm dP}}_{2}\dashrightarrow\widetilde{(Q')^{*}}$
defined by (\ref{eq:mz}) and (\ref{eq:tildedp2 Sp}). 

\subsubsection{\textbf{The rational map from $\overline{{\rm G}}_{\mathscr{I}\mathscr{R}}$
to the double cover of $\mP(V^{*})$\label{subsec:The-rational-mapGoverline to P(V*)}}}

$\empty$

\vspace{3pt}

\noindent\textbf{General Case.} We set
\begin{align}
s(x,r): & =(x_{1}x_{4}+x_{2}x_{3})r_{12}-x{}_{2}^{2}r_{13}+x{}_{1}^{2}r_{14}+2x_{1}x_{2}r_{15}\label{eq:s(x,r)Gen}\\
 & +x_{4}^{2}r_{23}-x{}_{3}^{2}r_{24}+2x_{3}x_{4}r_{25}+(x_{2}x_{3}-x_{1}x_{4})r_{34}\nonumber \\
 & -2x_{2}x_{4}r_{35}+2x_{1}x_{3}r_{45},\nonumber 
\end{align}
and recall that $m_{ij}(r,s)$ are as in (\ref{eq:mijGen}). The rational
map $\overline{{\rm G}}_{\mathscr{I}\mathscr{R}}\dashrightarrow{\rm G}(2,W)$
is defined by
\begin{equation}
[\overline{x},\overline{r}]\mapsto[m_{12}(\overline{r}),m_{13}(\overline{r}),m_{14}(\overline{r},s(\overline{x},\overline{r})),m_{23}(\overline{r},s(\overline{x},\overline{r})),m_{24}(\overline{r}),m_{34}(\overline{r})].\label{eq:GenGoverline G2W}
\end{equation}


\vspace{10pt}

\noindent\textbf{Special Case.} We set
\begin{align}
s & (x,r):=2(x_{1}x_{2}r_{15}+x_{3}x_{4}r_{25}-x_{2}x_{4}r_{35}+x_{1}x_{3}r_{45})\label{eq:s(x,r)Sp}
\end{align}
and recall that $m_{ij}(r,s)$ are as in (\ref{eq:mijSp}). The rational
map $\overline{{\rm G}}_{\mathscr{IR}}\dashrightarrow\widetilde{(Q')^{*}}$
is defined by 
\begin{equation}
[\overline{x},\overline{r}]\mapsto[m_{12}(\overline{r}),m_{13}(\overline{r}),m_{23}(\overline{r},s(\overline{x},\overline{r})),m_{24}(\overline{r}),m_{34}(\overline{r}),{\rm Pf}_{1234}(\overline{r})].\label{eq:SpGoverline G2W}
\end{equation}


\subsubsection{\textbf{The rational map $\overline{{\rm G}}_{\mathscr{I}\mathscr{R}}\protect\dashrightarrow\mP(\sU)\,\,\text{or}\,\,\mP(\sU')$\label{subsec:Rational-mapGoverline to P(U)}}}

By an explicit calculation, we see that $\mP(\sU)$ is defined in
$\mP(W)\times\mP(\wedge^{2}W)$ by 
\begin{equation}
\begin{cases}
\mathsf{x}_{1}\mathsf{m}_{23}-\mathsf{x}_{2}\mathsf{m}_{13}+\mathsf{x}_{3}\mathsf{m}_{12}=0\\
\mathsf{x}_{1}\mathsf{m}_{24}-\mathsf{x}_{2}\mathsf{m}_{14}+\mathsf{x}_{4}\mathsf{m}_{12}=0\\
\mathsf{x}_{1}\mathsf{m}_{34}-\mathsf{x}_{3}\mathsf{m}_{14}+\mathsf{x}_{4}\mathsf{m}_{13}=0\\
\mathsf{x}_{2}\mathsf{m}_{34}-\mathsf{x}_{3}\mathsf{m}_{24}+\mathsf{x}_{4}\mathsf{m}_{23}=0\\
\mathsf{m}_{12}\mathsf{m}_{34}-\mathsf{m}_{13}\mathsf{m}_{24}+\mathsf{m}_{14}\mathsf{m}_{23}=0,
\end{cases}\label{eq:Univ}
\end{equation}
where $\mathsf{x}_{1},\dots,\mathsf{x}_{4}$ are the coordinates of
$\mP(W)$. 

\vspace{3pt}

\noindent\textbf{General Case.} By a direct computation, we can check
the following:
\begin{lem}
\label{lem:IR-RelationsGen} The points of $\overline{{\rm G}}_{\mathscr{IR}}$
satisfy the following equations with\textup{ (\ref{eq:mijGen})} and
(\ref{eq:s(x,r)Gen}):\textup{
\begin{equation}
\begin{cases}
\overline{x}_{1}m_{23}(\overline{r},s(\overline{x},\overline{r}))-\overline{x}_{2}m_{13}(\overline{r})+\overline{x}_{3}m_{12}(\overline{r})=0\\
\overline{x}_{1}m_{24}(\overline{r})-\overline{x}_{2}m_{14}(\overline{r},s(\overline{x},\overline{r}))+\overline{x}_{4}m_{12}(\overline{r})=0\\
\overline{x}_{1}m_{34}(\overline{r})-\overline{x}_{3}m_{14}(\overline{r},s(\overline{x},\overline{r}))+\overline{x}_{4}m_{13}(\overline{r})=0\\
\overline{x}_{2}m_{34}(\overline{r})-\overline{x}_{3}m_{24}(\overline{r})+\overline{x}_{4}m_{23}(\overline{r},s(\overline{x},\overline{r}))=0\\
m_{12}(\overline{r})m_{34}(\overline{r})-m_{13}(\overline{r})m_{24}(\overline{r})+m_{14}(\overline{r},s(\overline{x},\overline{r}))m_{23}(\overline{r},s(\overline{x},\overline{r}))=0.
\end{cases}\label{eq:Univ-3}
\end{equation}
}
\end{lem}

By (\ref{eq:GenGoverline G2W}), (\ref{eq:Univ}) and Lemma \ref{lem:IR-RelationsGen},
a rational map $\overline{{\rm G}}_{\mathscr{I}\mathscr{R}}\dashrightarrow{\rm \mP(\sU)}$
is defined by
\begin{equation}
[\overline{x},\overline{r}]\mapsto[\overline{x}]\times[m_{12}(\overline{r}),m_{13}(\overline{r}),m_{14}(\overline{r},s(\overline{x},\overline{r})),m_{23}(\overline{r},s(\overline{x},\overline{r})),m_{24}(\overline{r}),m_{34}(\overline{r})].\label{eq:genGoverline PU}
\end{equation}

\vspace{10pt}

\noindent\textbf{Special Case.} We denote by $\mP_{\widetilde{(Q')^{*}}}\simeq\mP^{5}$
the ambient projective space of $\widetilde{(Q')^{*}}$. Noting (\ref{eq:mz}),
we choose the coordinates of $\mP_{\widetilde{(Q')^{*}}}$ as $\mathsf{m}_{12}$,
$\mathsf{m}_{13}$, $\mathsf{m}_{23}$, $\mathsf{m}_{24}$, $\mathsf{m}_{34}$
and $\mathsf{z}_{5}$. By an explicit calculation, we see that the
closure of $\mP(\sU')$ in $\mP(W)\times\mP_{\widetilde{(Q')^{*}}}$
is defined by $\mathsf{m}_{14}=\mathsf{m}_{23}$ and the equations
(\ref{eq:Univ}). 

By a direct computation, we can check the following:
\begin{lem}
\label{lem:IR-RelationsSP} With $m_{ij}(r,s)$\textup{ as in (\ref{eq:mijSp})}
and $s(x,r)$ as in (\ref{eq:s(x,r)Sp}), the points of $\overline{{\rm G}}_{\mathscr{IR}}$
satisfy the equations (\ref{eq:Univ-3}) by setting $m_{14}=m_{23}$.
\end{lem}

By (\ref{eq:SpGoverline G2W}), (\ref{eq:Univ}) with $\mathsf{m}_{14}=\mathsf{m}_{23}$
and Lemma \ref{lem:IR-RelationsSP}, a rational map $\overline{{\rm G}}_{\mathscr{I}\mathscr{R}}\dashrightarrow{\rm \mP(\sU')}$
is defined by
\begin{equation}
[\overline{x},\overline{r}]\mapsto[\overline{x}]\times[m_{12}(\overline{r}),m_{13}(\overline{r}),m_{23}(\overline{r},s(\overline{x},\overline{r})),m_{24}(\overline{r}),m_{34}(\overline{r}),{\rm Pf}_{1234}(\overline{r})].\label{eq:SpOverline PU}
\end{equation}

\subsubsection{\textbf{The quasi-$\mP^{1}$-bundle }$\varphi\colon\widetilde{{\rm G}}_{\mathscr{\mathscr{I}R}}\to\widetilde{{\rm dP}}_{2}$\label{subsec:The-quasi--P1bundleG gen}}

By (\ref{eq:neweq}) and Lemmas \ref{lem:IR-RelationsGen} and \ref{lem:IR-RelationsSP},
we are led to define a variety $\widetilde{{\rm G}}_{\mathscr{IR}}$
over $\widetilde{{\rm dP}}_{2}$, which plays a role connecting $\overline{{\rm G}}_{\mathscr{IR}}$
and $\mP(\sU)$ or $\mP(\sU')$. 

First we define the toric variety $\mathbb{E}_{\widetilde{\mathscr{IR}}}$
as the quotient of $\left((\wedge^{2}V\oplus\mC)\setminus\{\bm{o}\}\right)\times\left((W\oplus\mC)\setminus\{\bm{o}\}\right)$
by the action of $(\mC^{*})^{2}$ with the weights given as in Table
\ref{Table:FIR}, where the entries of the second line are the coordinates
of $\left((\wedge^{2}V\oplus\mC)\setminus\{\bm{o}\}\right)\times\left((W\oplus\mC)\setminus\{\bm{o}\}\right)$,
those of the third and the fourth lines are the weights of the coordinates
by the action of the first and the second factors of $(\mC^{*})^{2}$
respectively.
\begin{table}[h]
\begin{tabular}{|c|c|c|c|c|}
\hline 
spaces & \multicolumn{2}{c|}{$(\wedge^{2}V\oplus\mC)\setminus\{\bm{o}\}$} & \multicolumn{2}{c|}{$(W\oplus\mC)\setminus\{\bm{o}\}$}\tabularnewline
\hline 
coordinates & $\widetilde{r}_{ij}\,(1\leq i<j\leq5)$ & $\widetilde{s}$ & $\widetilde{x}$ & $w$\tabularnewline
\hline 
1st weights & $1$ & 2 & $0$ & $-1$\tabularnewline
\hline 
2nd weights & $0$ & 0 & $1$ & $2$\tabularnewline
\hline 
\end{tabular}

\caption{Toric variety $\mathbb{\mathbb{E}_{\widetilde{\mathscr{IR}}}}$}
\label{Table:FIR}
\end{table}
 The quotient of $(\wedge^{2}V\oplus\mC)\setminus\{\bm{o}\}$ by the
action of the first factor of $(\mC^{*})^{2}$ is nothing but $\mP(1^{10},2)$.
Hence we have the natural projection $\mathbb{E_{\widetilde{\mathscr{IR}}}}\to\mP(1^{10},2)$,
which is seen to be a $\mP(1^{4},2)$-bundle. We denote by $\mF_{\widetilde{\mathscr{IR}}}\to\widetilde{{\rm dP}}_{2}$
the restriction of $\mathbb{E_{\widetilde{\mathscr{IR}}}}\to\mP(1^{10},2)$
over $\widetilde{{\rm dP}}_{2}$.
\begin{defn}
\label{def:GtildeIR}Let $m_{ij}$ be as in (\ref{eq:mijGen}) in
General Case, or as in (\ref{eq:mijSp}) considering $m_{14}=m_{23}$
in Special Case, and $\widetilde{q}_{ij}:=q_{ij}(\widetilde{x})\,(1\leq i<j\leq5)$.
We define $\widetilde{{\rm G}}_{\mathscr{\mathscr{I}R}}\subset\mF_{\widetilde{\mathscr{\mathscr{I}R}}}$
by 
\begin{equation}
\begin{cases}
\widetilde{x}_{1}m_{23}(\widetilde{r},\widetilde{s})-\widetilde{x}_{2}m_{13}(\widetilde{r})+\widetilde{x}_{3}m_{12}(\widetilde{r})=0,\\
\widetilde{x}_{1}m_{24}(\widetilde{r})-\widetilde{x}_{2}m_{14}(\widetilde{r},\widetilde{s})+\widetilde{x}_{4}m_{12}(\widetilde{r})=0,\\
\widetilde{x}_{1}m_{34}(\widetilde{r})-\widetilde{x}_{3}m_{14}(\widetilde{r},\widetilde{s})+\widetilde{x}_{4}m_{13}(\widetilde{r})=0,\\
\widetilde{x}_{2}m_{34}(\widetilde{r})-\widetilde{x}_{3}m_{24}(\widetilde{r})+\widetilde{x}_{4}m_{23}(\widetilde{r},\widetilde{s})=0,\\
m_{12}(\widetilde{r})m_{34}(\widetilde{r})-m_{13}(\widetilde{r})m_{24}(\widetilde{r})+m_{14}(\widetilde{r},\widetilde{s})m_{23}(\widetilde{r},\widetilde{s})=0,
\end{cases}\label{eq:Univ-1}
\end{equation}
and
\begin{align}
 & (\widetilde{r}_{ij}\widetilde{q}_{kl}+\widetilde{r}_{kl}\widetilde{q}_{ij})-(\widetilde{r}_{ik}\widetilde{q}_{jl}+\widetilde{r}_{jl}\widetilde{q}_{ik})+(\widetilde{r}_{il}\widetilde{q}_{jk}+\widetilde{r}_{jk}\widetilde{q}_{il})+w{\rm Pf}_{ijkl}(\widetilde{r})=0\label{eq:Rq}
\end{align}
for $(i,j,k,l)=(1,2,3,4),(1,2,3,5),(1,2,4,5),(1,3,4,5),(2,3,4,5)$,
where we consider $m_{14}=m_{23}$ in Special Case. Note that the
equations in (\ref{eq:Univ-1}) (resp.~(\ref{eq:Rq})) originate
from Lemmas \ref{lem:IR-RelationsGen} and \ref{lem:IR-RelationsSP}
(resp.~(\ref{eq:neweq})). 
\end{defn}

We denote by $\varphi\colon{\rm \widetilde{G}}_{\mathscr{IR}}\to\widetilde{{\rm dP}}_{2}$
the natural projection. We investigate basic properties of $\varphi$. 

\vspace{3pt}

\noindent \textbf{General Case.} It is convenient to regard the quadric
$Q$ as a homogeneous space of the group ${\rm Sp}(W)$. The following
fact is a standard one in projective geometry:
\begin{lem}
\label{lem:Sp-orbit} We consider ${\rm G}(2,V)$ as a subvariety
of $\mP(\wedge^{2}V)$ with coordinates $\widetilde{r}_{ij}\,(1\leq i<j\leq5)$.
The ${\rm {\rm Sp}}(W)$-orbits on ${\rm G}(2,V)$ are classified
as follows: \begin{enumerate}[$(a)$]

\item The orbit of a line intersecting simply $Q$ at two points,
which is the orbit of the $\widetilde{r}_{12}$-point and the unique
open orbit.

\item The orbit of a tangent line of $Q$, which is the orbit of
the $\widetilde{r}_{45}$-point.

\item The orbit of a line contained in $Q$, which is the orbit of
the $\widetilde{r}_{23}$-point.

\end{enumerate}
\end{lem}

\begin{prop}
\label{prop:rho-fiber-1}It holds that 

\begin{enumerate}[$(1)$]

\item the morphisms $\varphi$ is a $\mP^{1}$-bundle outside ${\rm G}(2,V)=\{\widetilde{s}=0,{\rm Pf}_{\bullet}(\widetilde{r})=0\}$,
and

\item the $\varphi$-fiber over a point in ${\rm G}(2,V)$ is described
as follows, where we use the notation as in Lemma \ref{lem:Sp-orbit}:

\begin{itemize}

\item The fiber over a point in the orbit $(a)$ is a pair of two
copies of $\mP(1^{2},2)$ intersecting at the $w$-point. 

\item The fiber over a point in the orbit $(b)$ is a $\mP(1^{2},2)$. 

\item The fiber over a point in the orbit $(c)$ is a $\mP(1^{3},2)$.

\end{itemize} 

\end{enumerate}
\end{prop}

\begin{proof}
(1). Note that the open subset $\widetilde{{\rm dP}}_{2}\setminus{\rm G}(2,V)$
is the union of the open subsets $\{m_{ij}\not=0\}\,(1\leq i<j\leq4)$.
Similarly to Proposition \ref{prop:f'}, we can check that $\varphi$
is a $\mP^{1}$-bundle over each of them. 

(2). By the ${\rm Sp}(W)$-action on ${\rm G}(2,V)$ as in Lemma \ref{lem:Sp-orbit},
we have only to determine the following fibers:

\begin{itemize}

\item The fiber over the $\widetilde{r}_{12}$-point: $\{\widetilde{q}_{34}=\widetilde{q}_{35}=\widetilde{q}_{45}=0\}=\{\widetilde{x}_{1}=\widetilde{x}_{2}=0\}\cup\{\widetilde{x}_{3}=\widetilde{x}_{4}=0\}$,
which is a pair of two copies of $\mP(1^{2},2)$ intersecting at the
$w$-point.

\item The fiber over the $\widetilde{r}_{45}$-point: $\{\widetilde{q}_{12}=\widetilde{q}_{13}=\widetilde{q}_{23}=0\}=2\{\widetilde{x}_{1}=\widetilde{x}_{3}=0\}$,
which is a $\mP(1^{2},2)$.

\item The fiber over the $\widetilde{r}_{23}$-point: $\{\widetilde{q}_{14}=\widetilde{q}_{15}=\widetilde{q}_{45}=0\}=\{\widetilde{x}_{4}=0\}$,
which is a $\mP(1^{3},2)$.

\end{itemize}
\end{proof}
\vspace{10pt}

\noindent \textbf{Special Case. }Let ${\rm SO}(V')$ be the special
orthogonal group associated to the quadric $Q'$, and $H$ the subgroup
of ${\rm GL}(V)$ consisting elements of the form $\left(\begin{array}{cc}
A & \bm{b}\\
\bm{o} & c
\end{array}\right)$ with $A\in{\rm SO}(V')$, $\bm{b}\in\mC^{4}$ and $c\in\mC^{\times}$.

By linear algebra, the following fact can be easily verified:
\begin{lem}
\label{lem:SO(4)}We consider ${\rm G}(2,V)$ as a subvariety of $\mP(\wedge^{2}V)$
with coordinates $\widetilde{r}_{ij}\,(1\leq i<j\leq5)$. The $H$-orbits
on ${\rm G}(2,V)$ are classified as follows: \begin{enumerate}[$(a)$]

\item The orbit of a line disjoint from $Q'$, which is the orbit
of the point with $\widetilde{r}_{15}=-\widetilde{r}_{25}=1$ and
the remaining coordinates being zero.

\item The orbit of a line intersecting simply $Q'$ at one point,
which is the orbit of the $\widetilde{r}_{15}$-point.

\item The orbit of a line intersecting simply $Q'$ at two points,
which is the orbit of the point with $\widetilde{r}_{13}=-\widetilde{r}_{14}=\widetilde{r}_{23}=-\widetilde{r}_{24}=1$
and the remaining coordinates being zero.

\item The orbit of a tangent line of $Q'$, which is the orbit of
the point with $\widetilde{r}_{13}=-\widetilde{r}_{23}=1$ and the
remaining coordinates being zero.

\item The orbit of a line in $Q'$, which is the orbit of the $\widetilde{r}_{13}$-point.

\end{enumerate}
\end{lem}

We set 
\begin{align*}
\Delta & :=\{{\rm Pf}_{1235}(\widetilde{r})={\rm Pf}_{1245}(\widetilde{r})={\rm Pf}_{1345}(\widetilde{r})={\rm Pf}_{2345}(\widetilde{r})=0\}\cap\widetilde{{\rm dP}}_{2}\\
 & =\{\widetilde{s}=0,{\rm Pf}_{\bullet}(\widetilde{r})=0\}\cup\{\widetilde{s}=\widetilde{r}_{15}=\widetilde{r}_{25}=\widetilde{r}_{35}=\widetilde{r}_{45}=0\},
\end{align*}
where we note that ${\rm G}(2,V)=\{\widetilde{s}=0,{\rm Pf}_{\bullet}(\widetilde{r})=0\}=\{{\rm Pf}_{\bullet}(\widetilde{r})=0\}\cap\widetilde{{\rm dP}}_{2}$. 
\begin{prop}
It holds that 

\begin{enumerate}[$(1)$]

\item the morphisms $\varphi$ is a $\mP^{1}$-bundle outside $\Delta$,
and

\item the $\varphi$-fiber over a point in $\Delta$ is described
as follows, where we use the notation as in Lemma \ref{lem:SO(4)}:

\begin{itemize}

\item The fiber over a point in the orbit $(a)$ is a pair of two
copies of $\mP(1^{2},2)$ intersecting at the $w$-point. 

\item The fiber over a point in the orbit $(b)$ is a $\mP(1^{2},2)$. 

\item The fiber over a point in the orbit $(c)$ is a quadric cone
in $\mP(1^{4},2)$.

\item The fiber over a point in the orbit (d) is a pair of two copies
of $\mP(1^{3},2)$.

\item The fiber over a point in the orbit (e) is a $\mP(1^{3},2).$

\item The fiber over a point in $\{\widetilde{s}=\widetilde{r}_{15}=\widetilde{r}_{25}=\widetilde{r}_{35}=\widetilde{r}_{45}=0\}\setminus{\rm G}(2,V)$
is a $\mP^{3}$.

\end{itemize} 

\end{enumerate}
\end{prop}

\begin{proof}
As in Proposition \ref{prop:rho-fiber-1}, we can directly check that
$\varphi$ is a $\mP^{1}$-bundle over each of the open subsets ${\rm Pf}_{ijkl}(\widetilde{r})\not=0$
with $(i,j,k,l)=(1,2,3,5)$, $(1,2,4,5)$, $(1,3,4,5)$ and $(2,3,4,5)$.

We can also directly check that the $\varphi$-fiber over a point
of $\{\widetilde{s}=\widetilde{r}_{15}=\widetilde{r}_{25}=\widetilde{r}_{35}=\widetilde{r}_{45}=0\}\setminus{\rm G}(2,V)$
is a copy of $\mP^{3}$.

Let us check the $\varphi$-fiber over a point of ${\rm G}(2,V)$.
Since ${\rm Pf}_{\bullet}=0$ on ${\rm G}(2,V)$, points of $\widetilde{{\rm G}}_{\mathscr{IR}}$
over ${\rm G}(2,V)$ satisfy the equations (\ref{eq:Rq}) with $w=0$.
We recall the situation of Proposition \ref{prop:bundleE}. Since
$\mP(\sE)$ is contained in $\mP(\wedge^{2}V)\times{\rm G}(2,V)$,
$\mP(q^{*}\sE)$ is contained in $\mP(\wedge^{2}V)\times\mP(W)$.
Let $\gamma$ be the fiber of $\mP(q^{*}\sE)\to\mP(\wedge^{2}V)$
over a point $\ell\in{\rm G}(2,V)$. Then the $\varphi$-fiber over
$\ell$ is the cone over $\gamma$ with the vertex being the $w$-point.
By (2) of the same proposition and Proposition \ref{prop:q(P(W))}
(1), $\gamma$ is the inverse image of the locus in $q(\mP(W))$ parameterizing
tangent lines or lines in $Q'$ intersecting $\ell$. Therefore the
description of the $\varphi$-fiber over $\ell$ follows by Proposition
\ref{prop:q(P(W))} (2) and (3).
\end{proof}

\subsubsection{\textbf{The birational morphism $\psi\colon\widetilde{{\rm G}}_{\mathscr{IR}}\to\overline{{\rm G}}_{\mathscr{IR}}$ }}

Let $\widetilde{{\rm G}}_{\mathscr{\mathscr{I}R}}\to\mP(1^{4},2^{10})$
be the morphism defined by

\begin{align}
\widetilde{{\rm G}}_{\mathscr{\mathscr{I}R}} & \ni[\widetilde{r},\widetilde{s};\widetilde{x},w]\mapsto[\widetilde{x},w\widetilde{r}]\in\mP(1^{4},2^{10}).\label{eq:psi}
\end{align}

This induces a morphism $\psi\colon\widetilde{{\rm G}}_{\mathscr{\mathscr{I}R}}\to\overline{{\rm G}}_{\mathscr{\mathscr{I}R}}$
by (\ref{eq:neweq}) and (\ref{eq:Rq}). Note that 
\begin{equation}
\psi^{-1}(\Pi)=\widetilde{{\rm G}}_{\mathscr{\mathscr{I}R}}\cap\{w=0\}.\label{eq:Psi-1Pi}
\end{equation}

\begin{prop}
\label{prop:The-birational-morphism psi}The following assertions
hold:

\begin{enumerate}[$(1)$]

\item The morphism $\psi$ is an isomorphism outside $\psi^{-1}(\Pi)$.

\item The induced morphism $\psi^{-1}(\Pi)\to\Pi$ is the $\mP^{6}$-bundle
$\mP(q^{*}\sE)$ over $\Pi\simeq\mP(W)$.

\end{enumerate}
\end{prop}

\begin{proof}
(1). Let $[\overline{x},\overline{r}]\in\overline{{\rm G}}_{\mathscr{\mathscr{I}R}}$
be a point of $\overline{{\rm G}}_{\mathscr{\mathscr{I}R}}$. The
morphism $\psi$ is surjective since it holds that $\psi([\overline{r},s(\overline{x},\overline{r});\overline{x},1])=[\overline{x},\overline{r}]$,
where $s(x,r)$ is defined as in (\ref{eq:s(x,r)Gen}) or (\ref{eq:s(x,r)Sp}),
and $[\overline{r},s(\overline{x},\overline{r});\overline{x},1]\in\widetilde{{\rm G}}_{\mathscr{\mathscr{I}R}}$
by (\ref{eq:neweq}), Lemma \ref{lem:IR-RelationsGen} or \ref{lem:IR-RelationsSP},
and Definition \ref{def:GtildeIR}. Now we show that $\psi$ is injective
over $\widetilde{{\rm G}}_{\mathscr{\mathscr{I}R}}\setminus\psi^{-1}(\Pi)$.
Assume that $\psi([\widetilde{r},\widetilde{s};\widetilde{x},w])=[\overline{x},\overline{r}]$
for a point $[\widetilde{r},\widetilde{s};\widetilde{x},w]\in\widetilde{{\rm G}}_{\mathscr{\mathscr{I}R}}\setminus\psi^{-1}(\Pi)$.
Since $\overline{r}\not=0$, we have $w\not=0$. By the $(\mC^{*})^{2}$-action,
we may assume that $w=1$. Therefore there exists $\alpha\in\mC^{*}$
such that $\widetilde{x}=\alpha\overline{x}$ and $\widetilde{r}=\alpha^{2}\overline{r}$.
By Table \ref{Table:FIR}, we see that
\[
[\widetilde{r},\widetilde{s};\widetilde{x},1]=[\alpha^{-2}\widetilde{r},\alpha^{-4}\widetilde{s};\alpha^{-1}\widetilde{x},1]=[\overline{r},\alpha^{-4}\widetilde{s};\overline{x},1].
\]
 Therefore we may assume that $\widetilde{r}=\overline{r}$ and $\widetilde{x}=\overline{x}.$
It remain to show the uniqueness of $\widetilde{s}$. If $\overline{x}\not=\bm{o}$,
then we may recover $\widetilde{s}$ by one of the first 4 equations
of (\ref{eq:Univ-1}). If $\overline{x}=\bm{o}$, then, by (\ref{eq:neweq}),
we have ${\rm Pf}_{ijkl}(\overline{r})=0\,(1\leq i<j<k<l\leq5)$.
Thus the 5th equation of (\ref{eq:Univ-1}) implies that $\widetilde{s}=0.$

(2). Similarly to Lemma \ref{lem:PE}, we see that $\psi^{-1}(\Pi)$
is isomorphic to the subvariety in $\mP(\wedge^{2}V)\times\mP(W)$
defined by the equation in (\ref{eq:Rq}) with $w=0$. Thus the second
assertion follows from Proposition \ref{prop:bundleE}.
\end{proof}

\subsection{Constructing the key varieties $\mathscr{\mathscr{I}R}$\label{subsec:IR}}

Let us construct the following Sarkisov link \begin{equation}\label{eq:SarkisovIR} \xymatrix{& \widehat{\mathscr{IR}}\ar@{-->}[rr]\ar[dl]_{\widehat{f}_{\mathscr{IR}}}\ar[dr]^{\widehat{g}_{\mathscr{IR}}} &  &\widetilde{\mathscr{IR}}\ar[dr]^{\widetilde{f}_{\mathscr{IR}}}\ar[dl]_{\widetilde{g}_{\mathscr{IR}}}\\
\mathscr{IR}& &\overline{\mathscr{IR}}&  & {\rm dP}_{2}.}
\end{equation}By taking a linear section of $\overline{{\rm G}}_{\mathscr{IR}}$,
we define $\overline{\mathscr{IR}}$ and further restrict ${\rm \widetilde{G}}_{\mathscr{IR}}\to\overline{{\rm G}}_{\mathscr{IR}}$
over $\overline{\mathscr{IR}}$ to obtain a small resolution $\widetilde{\mathscr{IR}}\to\overline{\mathscr{IR}}$.
The subsequent process is similar to that for Type $R$ though, for
Type $I\!R$, we do not provide explicit equations for $\text{\ensuremath{\mathscr{IR}}}$;
since $\text{\ensuremath{\mathscr{IR}}}$ depends on ${\rm dP}_{2}$,
its defining equations are more complicated than those of $\mathscr{R}$,
which we prefer to avoid in this paper. Fortunately, the flop $\widetilde{\mathscr{IR}}\dashrightarrow\widehat{\mathscr{IR}}$
is simple and can be handled without explicit equations. Moreover,
since the key variety $\mathscr{IR}$ is four-dimensional, the construction
can proceed in a reasonably concrete way by referring to the classification
of four-dimensional extremal contractions, even in the absence of
explicit equations. 
\begin{rem}
It seems that $\mathscr{IR}$ is not related with any $\mP^{1}\times\mP^{1}\times\mP^{1}$-fibration.
\end{rem}

\subsubsection{\textbf{The $3$-space $L_{\mathscr{IR}}\subset\mP(\wedge^{2}V)$
and its double cover\label{subsec:The-3-space IR}}}

Now we cut $\widetilde{{\rm dP}}_{2}$ by six linearly independent
linear forms $L_{1}(\widetilde{r}),\dots,L_{6}(\widetilde{r})$ with
respect to $\widetilde{r}=(\widetilde{r}_{ij}\,(1\leq i<j\leq5))$
and set 
\begin{equation}
{\rm dP}_{2}:=\widetilde{{\rm dP}}_{2}\cap\{L_{1}(\widetilde{r})=\cdots=L_{6}(\widetilde{r})=0\}.\label{eq:dP2}
\end{equation}
We may also consider that $L_{1}(\widetilde{r}),\dots,L_{6}(\widetilde{r})$
are linear forms on $\mP(\wedge^{2}V)$. We set 
\begin{align*}
L_{\mathscr{IR}} & :=\{L_{1}(\widetilde{r})=\cdots=L_{6}(\widetilde{r})=0\}\subset\mP(\wedge^{2}V),\\
B & :=B_{{\rm G}}\cap L_{\mathscr{IR}}.
\end{align*}

Note that $B$ is the branch locus of the double cover ${\rm dP}_{2}\to L_{\mathscr{IR}}$.
We assume the following:

\vspace{3pt}

\noindent \textbf{Conditions on $L_{\mathscr{IR}}$.}
\begin{enumerate}
\item ${\rm G}(2,V)\cap L_{\mathscr{IR}}$ consists of $5$ points, which
we denote by $\mathsf{q}_{1},\dots,\mathsf{q}_{5}$. This implies
that the inverse image of ${\rm G}(2,V)\cap L_{\mathscr{IR}}$ on
$\widetilde{{\rm dP}}_{2}$ also consists of $5$ points, which we
denote by $\widetilde{\mathsf{q}}_{1},\dots,\mathsf{\widetilde{q}}_{5}$. 
\item $B$ is a quartic surface with only $A_{\leq2}$-singularities at
$\mathsf{q}_{1},\dots,\mathsf{q}_{5}$. 
\end{enumerate}
\begin{lem}
\label{lem:Important}Let us take any ${\rm dP}_{2}$ as in Definition
\ref{def:dep2}. We may identify the target of the double cover ${\rm dP}_{2}\to\mP^{3}$
as $L_{\mathscr{IR}}\subset\mP(\wedge^{2}V)$ and we may assume that
the branch divisor $B$ of ${\rm dP}_{2}\to L_{\mathscr{IR}}$ satisfies
Conditions on $L_{\mathscr{IR}}$, namely, $B$ and its singular points
$\mathsf{q}_{1},\dots,\mathsf{q}_{5}$ are cut out by $L_{\mathscr{IR}}$
from $B_{G}$ and ${\rm G}(2,V)$, respectively. 
\end{lem}

\begin{proof}
Certainly, we may identify the target of the double cover ${\rm dP}_{2}\to\mP^{3}$
as $L_{\mathscr{IR}}\subset\mP(\wedge^{2}V)$ and we may assume that
the singular points $\mathsf{q}_{1},\dots,\mathsf{q}_{5}$ are cut
out by $L_{\mathscr{IR}}$ from ${\rm G}(2,V)$. Then the map $\mu\colon L_{\mathscr{IR}}\dashrightarrow{\rm S}_{3}$
is the restriction of $\mP(\wedge^{2}V)\dashrightarrow\mP(V^{*})$
described as (\ref{eq:zPf}). By the ${\rm GL}(V)$-action on $\mP(V^{*})$,
we transform $Q_{B}$ to $Q^{*}$. Let $L'_{\mathscr{IR}}\subset\mP(\wedge^{2}V)$
and ${\rm S}'_{3}\subset\mP(V^{*})$ be the images of $L_{\mathscr{IR}}$
and ${\rm S}{}_{3}$ by this transformation. Since $\mP(\wedge^{2}V)\dashrightarrow\mP(V^{*})$
is ${\rm G}(V)$-equivariant, the induced map $L'_{\mathscr{IR}}\dashrightarrow{\rm S}'_{3}$
is also the restriction of $\mP(\wedge^{2}V)\dashrightarrow\mP(V^{*})$.
Therefore we may assume $Q_{B}=Q^{*}$, $L_{\mathscr{IR}}=L'_{\mathscr{IR}}$
and ${\rm S}_{3}={\rm S}'_{3}$ from the beginning. Then, since the
branch locus $B$ of ${\rm dP}_{2}\to L_{\mathscr{IR}}$ is the birational
image of $Q^{*}\cap{\rm S}_{3}$, $B$ coincides with $B_{G}\cap L_{\mathscr{IR}}$.
\end{proof}
\begin{lem}
\label{lem:A1A2} Both in General and Special Cases, the five points
$\mathsf{\widetilde{q}}_{1},\dots,\widetilde{\mathsf{q}}_{5}$ are
contained in either the orbit $(a)$ or $(b)$. In Special Case, $L_{\mathscr{IR}}$
is disjoint from $\{\widetilde{r}_{15}=\widetilde{r}_{25}=\widetilde{r}_{35}=\widetilde{r}_{45}=0\}$. 
\end{lem}

\begin{proof}
\noindent \textbf{General Case.} By a direct computation, we see
that, at a point of the orbit of a line contained in $Q$, $B_{G}$
has a hypersurface singularity whose quadratic term has rank $1$
(we have only to check this at the $p_{23}$-point). Therefore such
a point cannot be an $A_{\leq2}$-singularity after cutting $B_{G}$
by six linearly independent linear forms.

\vspace{5pt}

\noindent \textbf{Special Case.} The same argument works in this
case for a point in ${\rm G}(2,V)$. Hence the five points $\mathsf{\widetilde{q}}_{1},\dots,\widetilde{\mathsf{q}}_{5}$
are contained in either the orbit $(a)$ or $(b)$ and the first assertion
follows. We can directly check that ${\rm G(2,V)\cap}\{\widetilde{r}_{15}=\widetilde{r}_{25}=\widetilde{r}_{35}=\widetilde{r}_{45}=0\}$
is the union of the orbits (c), (d), and (e), and hence $L_{\mathscr{IR}}$
is disjoint from $G(2,V)\cap\{\widetilde{r}_{15}=\widetilde{r}_{25}=\widetilde{r}_{35}=\widetilde{r}_{45}=0\}$
by the first assertion. Therefore, by the assumptions (1) and (2)
above, $\Sing({\rm G}(2,V)\cap L_{\mathscr{IR}})$ consists of $5$
points in $G(2,V)\setminus\{\widetilde{r}_{15}=\widetilde{r}_{25}=\widetilde{r}_{35}=\widetilde{r}_{45}=0\}$.
Since we can easily check that $B_{{\rm G}}$ is singular along $\{\widetilde{r}_{15}=\widetilde{r}_{25}=\widetilde{r}_{35}=\widetilde{r}_{45}=0\}$,
$L_{\mathscr{IR}}$ is disjoint from $\{\widetilde{r}_{15}=\widetilde{r}_{25}=\widetilde{r}_{35}=\widetilde{r}_{45}=0\}\setminus G(2,V)$
. Thus we have shown the second assertion.
\end{proof}
\begin{defn}
We define $\widetilde{\mathscr{\mathscr{I}R}}$ to be the restriction
of $\widetilde{{\rm G}}_{\mathscr{\mathscr{I}R}}$ over ${\rm dP}_{2}$,
and set $\widetilde{f}{}_{\mathscr{IR}}:=\varphi|_{\widetilde{\mathscr{\mathscr{I}R}}}$.
\end{defn}

\subsubsection{\textbf{The quasi-$\mP^{1}$-bundle }$\widetilde{{\rm \mathscr{\mathscr{I}R}}}$\textbf{
over ${\rm dP}_{2}$\label{subsec:The-quasi--bundle-IR}}}

We immediately derive the following proposition from Proposition \ref{prop:rho-fiber-1}
and Lemma \ref{lem:A1A2}:
\begin{prop}
\label{prop:five1/2}The morphism $\widetilde{f}{}_{\mathscr{IR}}$
is a $\mP^{1}$-bundle outside the $5$ points $\mathsf{\widetilde{q}}_{1},\dots,\widetilde{\mathsf{q}}_{5}$.
The $\widetilde{f}{}_{\mathscr{IR}}$- fiber over $\mathsf{\widetilde{q}}_{i}\,(1\le i\leq5)$
is the union of two copies of $\mP(1^{2},2)$ in $\mP(1^{4},2)$ intersecting
at the $w$-point when $\widetilde{\mathsf{q}}_{i}$ is contained
in the orbit $(a)$, and a $\mP(1^{2},2)$ when $\widetilde{\mathsf{q}}_{i}$
is contained in the orbit $(b)$. 
\end{prop}

\subsubsection{\textbf{The birational map $\widetilde{\mathscr{I}{\rm \mathscr{R}}}\protect\dashrightarrow\mathsf{IR}$\label{subsec:The-birational-mapIR PU}}}

$\empty$

\vspace{3pt}

\noindent\textbf{General Case.} By Table \ref{Table:FIR} and the
construction as in Subsection \ref{subsec:Rational-mapGoverline to P(U)},
we obtain the rational map $\mF_{\widetilde{\mathscr{IR}}}\dashrightarrow{\rm G}(2,W)\times\mP(V)$
defined by 
\begin{equation}
[\widetilde{r},\widetilde{s};\widetilde{x},w]\mapsto[m_{ij}(\widetilde{r},\widetilde{s})\,(1\leq i<j\leq4)]\times[\widetilde{x}],\label{eq:IRtilde to sfIR}
\end{equation}
where we recall that $m_{ij}(\widetilde{r},\widetilde{s})$ are defined
as in (\ref{eq:mijGen}). By (\ref{eq:Univ}) and (\ref{eq:Univ-1}),
this induces the rational map $\widetilde{G}_{\mathscr{IR}}\dashrightarrow\mP(\sU)$.
We can check the following square is commutative:

\begin{equation}\label{eq:sqG2}
\xymatrix{\widetilde{{\rm G}}_{\mathscr{IR}}\ar@{-->}[r]\ar[d]&\mP(\sU)\ar[d]\\
\widetilde{\rm{dP}}_2 \ar@{-->}[r]& \rm{G}(2,W).}
\end{equation}

\vspace{5pt}

\noindent\textbf{Special Case.} By Table \ref{Table:FIR} and the
construction as in Subsection \ref{subsec:Rational-mapGoverline to P(U)},
we obtain the rational map $\mF_{\widetilde{\mathscr{IR}}}\dashrightarrow\widetilde{(Q')^{*}}\times\mP(V)$
defined by 
\begin{equation}
[\widetilde{r},\widetilde{s};\widetilde{x},w]\mapsto[m_{ij}(\widetilde{r},\widetilde{s})\,(1\leq i<j\leq4)]\times[\widetilde{x}],\label{eq:IRtilde to sfIR2}
\end{equation}
where we recall that $m_{ij}(\widetilde{r},\widetilde{s})$ are defined
as in (\ref{eq:mijSp}) and we consider $m_{14}(\widetilde{r},\widetilde{s})=m_{23}(\widetilde{r},\widetilde{s})$.
By (\ref{eq:Univ}) with $\mathsf{m}_{14}=\mathsf{m}_{23}$ and (\ref{eq:Univ-1}),
this induces the rational map $\widetilde{G}_{\mathscr{IR}}\dashrightarrow\mP(\sU')$.
We can check the following square is commutative:

\begin{equation}\label{eq:sqG3}
\xymatrix{\widetilde{{\rm G}}_{\mathscr{IR}}\ar@{-->}[r]\ar[d]&\mP(\sU')\ar[d]\\
\widetilde{\rm{dP}}_2 \ar@{-->}[r]& \widetilde{(Q')^{*}} .}
\end{equation}

Let $\mathsf{f}_{\mathscr{IR}}\colon\mathsf{IR}\to\mathcal{S}$ be
the natural projection, where $\mathsf{IR}$ is defined as in (\ref{eq:mathsfIR}).
Restricting the diagrams (\ref{eq:sqG2}) and (\ref{eq:sqG3}), we
obtain the following commutative square:

\begin{equation}
\xymatrix{\widetilde{\mathscr{IR}}\ar@{-->}[r]\ar[d]_{\widetilde{f}_{\mathscr{IR}}}&\mathsf{IR}\ar[d]^{\mathsf{f}_{\mathscr{IR}}}\\
\rm{dP}_2 \ar@{-->}[r]& \mathcal{S}.}
\end{equation}

We can show the following proposition in the same way as Proposition
\ref{prop:Rtilde and P(Omega2)}.
\begin{prop}
\label{prop:IRtilde and P(U)S}The rational map $\widetilde{{\rm \mathscr{IR}}}\dashrightarrow\mathsf{IR}$
is an isomorphism outside $(\rho_{\mP}\circ\widetilde{f}_{\mathscr{IR}})^{-1}(\ell_{ij})\,(1\leq i<j\leq5)$.
The $\widetilde{f}_{\mathscr{IR}}$-fiber over a point of $\rho_{\mP}^{-1}(\ell_{ij}\setminus\{\mathsf{q}_{i},\mathsf{q}_{j})\}$
is mapped isomorphically to the $\mathsf{f}_{\mathscr{IR}}$-fiber
over the corresponding point in $\rho_{{\rm S}}^{-1}(\mathsf{q}_{ij})$
on $\mathcal{S}$.
\end{prop}

\begin{prop}
\label{prop:The-Picard-numbertwo} It holds that $\widetilde{\mathscr{IR}}$
has Picard number $2$ and only five $\nicefrac{1}{2}$-singularities
at the $\nicefrac{1}{2}$-singularities of a $\mP(1^{2},2)$ in the
fiber over $\widetilde{\mathsf{q}}_{i}\,(1\leq i\leq5)$. 
\end{prop}

\begin{proof}
Let $U$ be the base change of $\mathsf{IR}$ by the morphism $\nu_{\mathcal{S}}'\colon{\rm BldP}_{2}\to\mathcal{S}$
as in the diagram (\ref{eq:doublecover}). The variety $U$ is smooth
since so is ${\rm BldP}_{2}$. 

We consider Special Case in this paragraph. Note that $\mathsf{Q}_{i}$
does not contain the vertex $\widetilde{v}$ of $\widetilde{(Q')^{*}}$
by Proposition \ref{prop:23CI}. Therefore $\mathsf{Q}_{i}$ is isomorphic
to the image of $\mathsf{Q}_{i}$ by $\widetilde{(Q')^{*}}\setminus\{\widetilde{v}\}\to{\rm G}(2,W)$,
and the restriction of $\mathsf{IR}$ over $\mathsf{Q}_{i}$ is isomorphic
to that of $\mP(\sU)$ over the image of $\mathsf{Q}_{i}$. So, abusing
notation, we denote by $\mP(\sU|_{\mathsf{Q}_{i}})$ the restriction
of $\mathsf{IR}$ over $\mathsf{Q}_{i}$.

Hereafter we do not distinguish General and Special Cases until the
end of the proof. Let $U_{\mathsf{Q}_{i}'}$ be the strict transform
of $\mP(\sU|_{\mathsf{Q}_{i}})$ on $U$. Since the strict transform
$\mathsf{Q}_{i}'$ of $\mathsf{Q}_{i}$ on ${\rm BldP}_{2}$ is isomorphic
to $\mathsf{Q}_{i}$, $U_{\mathsf{Q}_{i}'}$ is isomorphic to $\mP(\sU|_{\mathsf{Q}_{i}})$.
Note that ${\rm bl}_{{\rm dP}_{2}}\colon{\rm BldP}_{2}\to{\rm dP}_{2}$
is the contraction of the quadric surface $\mathsf{Q}'_{i}$ near
each $\mathsf{Q}_{i}'$. Since ${\rm dP}_{2}$ is factorial by Definition
\ref{def:dep2}, the relative Picard number of ${\rm bl}_{{\rm dP}_{2}}$
near each $\mathsf{Q}_{i}'$ is $1$. This implies that the relative
Picard number of $U\to{\rm dP}_{2}$ is $2$ near each $U_{\mathsf{Q}_{i}'}$
since $U_{\mathsf{Q}_{i}'}\to\mathsf{Q}'_{i}$ is a $\mP^{1}$-bundle. 

Now we play the 2-ray game over ${\rm dP}_{2}$ near $\widetilde{\mathsf{q}}_{i}$
starting from $U$. Since $\mP(\sU|_{\mathsf{Q}_{i}})$ is contained
in $\mathsf{Q}_{i}\times\mP(W)$, the natural projection $p_{i}\colon\mP(\sU|_{\mathsf{Q}_{i}})\to\mP(W)\simeq\mP^{3}$
is induced. We will see below that $p_{i}\,(1\leq i\leq5)$ are birational
morphisms and these extend to a numerically $K$-trivial small contraction
$U\to U'$ over ${\rm dP}_{2}$ such that the exceptional locus coincides
with the union of the $p_{i}$-exceptional loci, and it is a flopping
contraction near each $p_{i}$-exceptional locus. By \cite{bchm},
there exists the birational map $U\dashrightarrow U^{+}$ which coincides
with the flop near each $p_{i}$-exceptional locus. Constructing explicitly
$U\dashrightarrow U^{+}$ near each $p_{i}$-exceptional locus, let
us see that $U\dashrightarrow U^{+}$ induces $p_{i}\colon\mP(\sU|_{\mathsf{Q}_{i}})\to\mP(W)\simeq\mP^{3}$
(hence the strict transforms of $\mP(\sU|_{\mathsf{Q}_{i}})$ on $U^{+}$
are isomorphic to $\mP^{3}$), and their normal bundles are $\sO_{\mP^{3}}(-2)$.

\vspace{3pt}

\noindent \textit{Case $\mathsf{Q}_{i}\simeq\mP^{1}\times\mP^{1}$}:
By \cite[Def.1.3, Thm.1.4 and Ex.1.5]{O}, we have $\sU|_{\mathsf{Q}_{i}}\simeq\sO_{\mP^{1}\times\mP^{1}}(-1,0)\oplus\sO_{\mP^{1}\times\mP^{1}}(0,-1)$.
This implies that the natural projection $U_{\mathsf{Q}_{i}'}\simeq\mP(\sU|_{\mathsf{Q}_{i}})\to\mP(W)$
is the blow-up of $\mP^{3}$ along the two mutually disjoint lines
$\delta_{1}$ and $\delta_{2}$ which are the images of the divisors
$\mP(\sO_{\mP^{1}\times\mP^{1}}(-1,0)\oplus0)$ and $\mP(0\oplus\sO_{\mP^{1}\times\mP^{1}}(0,-1))$.
By a numerical calculation, we see that, inside the 4-fold $U$, these
divisors are connected components of the exceptional locus of a flopping
contraction near $U_{\mathsf{Q}_{i}'}$. Computing the normal bundle
of flopping curves, we see the flopping contraction is of Atiyah's
type along these exceptional loci. Constructing the flop locally by
composing a blow-up and a blow-down in a standard way, we see the
strict transform $U_{\mathsf{Q}_{i}'}^{+}$ of $U_{\mathsf{Q}_{i}'}$
on the flopped $4$-fold $U^{+}$ is $\mP^{3}\simeq\mP(W)$ with normal
bundle $\sO_{\mP^{3}}(-2)$. 

\vspace{3pt}

\noindent  \textit{Case $\mathsf{Q}_{i}\simeq\mP(1^{2},2)$}: Note
that the quadric surface $\mathsf{Q}_{i}$ is contained in a smooth
hyperplane section $L$ of ${\rm G}(2,W)$. Since ${\rm G}(2,W)$
and the bundle $\sU$ are homogeneous under the action of ${\rm GL}(W)$,
the bundle $\sU|_{L}$ is uniquely determined up to isomorphism, which
is nothing but the Spinor bundle on the smooth quadric 3-fold $L$.
Moreover, since $L$ and $\sU|_{L}$ are homogeneous under the action
of the Spinor group by \cite{O}, $\sU|_{\mathsf{Q}_{i}}$ is also
uniquely determined up to isomorphism. Therefore, we may assume 
\begin{equation}
\mathsf{Q}_{i}=\{\mathsf{m}_{12}=0,\mathsf{m}_{13}=\mathsf{m}_{24}\}\cap{\rm G}(2,W).\label{eq:Qi sing}
\end{equation}
Then, using (\ref{eq:Univ}) and (\ref{eq:Qi sing}), we see that
the natural projection $\mP(\sU|_{\mathsf{Q}_{i}})\to\mP(W)$ is an
isomorphism outside $\Delta:=\mP(\sU|_{\mathsf{Q}_{i}})\cap\{\mathsf{x}_{1}=\mathsf{x}_{2}=0\}$,
which is mapped by the projection onto the line $\delta:=\{\mathsf{x}_{1}=\mathsf{x}_{2}=0\}\subset\mP(W)$,
and $\Delta\to\delta$ is a $\mP^{1}$-bundle. We may directly check
that, $\Delta$ is the exceptional locus of a flopping contraction
near $U_{\mathsf{Q}_{i}'}$ inside $U$. Now let us construct the
flop of $\Delta$ (the following is a 4-fold version of the construction
in Step 3. Case $\mathsf{Q}_{i}'\simeq\mP(1^{2},2)$ of the proof
of Proposition \ref{prop:GenExOf f'-IR} (B). The discussion below
proceeds in parallel with that one). Let $\beta\colon U'\to U$ be
the blow-up of $U$ along $\Delta$, and $U'_{\mathsf{Q}_{i}'}$ the
strict transform of $U_{\mathsf{Q}_{i}'}$ on $U'$. Let $\Delta'\subset U'_{\mathsf{Q}_{i}'}$
be the strict transform of $\Delta$ by $\beta|_{U'_{\mathsf{Q}_{i}'}}$,
and $\Gamma$ the $\beta|_{U'_{\mathsf{Q}_{i}'}}$-exceptional divisor,
which is mapped to the fiber of the $\mP^{1}$-bundle $U_{\mathsf{Q}_{i}'}\to\mathsf{Q}_{i}'$
over the vertex of $\mathsf{Q}_{i}'$. By a local computation, we
see that $(\beta|_{U'_{\mathsf{Q}_{i}'}})^{*}\Delta=\Delta'+\Gamma$
and the intersection between the $\beta$-exceptional divisor and
$U'_{\mathsf{Q}_{i}'}$ coincides with $(\beta|_{U'_{\mathsf{Q}_{i}'}})^{*}\Delta$.
From this, we see that the normal bundle of a fiber of $\Delta'\to\delta$
in $U'$ is $\sO_{\mP^{1}}(-1)^{\oplus2}\oplus\sO_{\mP^{1}}$, and
hence $\Delta'$ is a flopping surface of Atiyah type. Therefore we
can flop $\Delta'$ by composing a blow-up and blow-down in a standard
way. The flop induces the contraction $U'_{\mathsf{Q}_{i}'}\to U''_{\mathsf{Q}_{i}'}$
of $\Delta'$ to $\delta$ and the isomorphism between the $\beta$-exceptional
divisor and its strict transform on the flopped $4$-fold for $\Delta'$.
The image of $\Gamma$ on the flopped 4-fold for $\Delta'$ is the
intersection between the strict transform of the $\beta$-exceptional
divisor and $U''_{\mathsf{Q}_{i}}$. Finally we can contract the strict
transform of the $\beta$-exceptional divisor and the flop of $\Delta$
is completed. Then we can check that the strict transform of $U_{\mathsf{Q}_{i}'}$
on the flopped 4-fold $U^{+}$ is isomorphic to $\mP^{3}\simeq\mP(W)$
with normal bundle $\sO_{\mP^{3}}(-2)$. 
\vspace{3pt}

Let us denote by $U^{+}\to\widetilde{\mathscr{IR}}'$ the contraction
of the strict transforms of $\mP(\sU|_{\mathsf{Q}_{i}})$ on $U^{+}$.
By Proposition \ref{prop:IRtilde and P(U)S}, $\widetilde{\mathscr{IR}}$
and $\widetilde{\mathscr{IR}}'$ are isomorphic outside the fibers
over $\widetilde{\mathsf{q}}_{1},\dots,\widetilde{\mathsf{q}}_{5}$.
By Proposition \ref{prop:five1/2}, the fibers of $\widetilde{\mathscr{IR}}$
over each $\widetilde{\mathsf{q}}_{i}$ is of codimension 2. By the
above construction, we can determine the fiber of $\widetilde{\mathscr{IR}}'\to{\rm dP}_{2}$
over each $\widetilde{\mathsf{q}}_{i}$ as follows: Assume that $\mathsf{Q}_{i}\simeq\mP^{1}\times\mP^{1}$.
Then, since the flopped surfaces corresponding to $\mP(\sO_{\mP^{1}\times\mP^{1}}(-1,0)\oplus0)$
and $\mP(0\oplus\sO_{\mP^{1}\times\mP^{1}}(0,-1))$ intersect $U_{\mathsf{Q}_{i}'}^{+}\simeq\mP^{3}$
along the lines $\delta_{1}$ and $\delta_{2}$, we see that the self-intersection
numbers of the lines in the flopped surfaces are $-2$. Therefore
the images of the flopped surfaces by the contraction of $U_{\mathsf{Q}_{i}'}^{+}$
are $\mP(1^{2},2)$, which intersect at the $\nicefrac{1}{2}$-singularity
and form the fiber of $\widetilde{\mathscr{IR}}'\to{\rm dP}_{2}$
over $\widetilde{\mathsf{q}}_{i}$. By the same reason, if $\mathsf{Q}_{i}\simeq\mP(1^{2},2)$,
the image of the flopped surface corresponding to $\Delta$ by the
contraction of $U_{\mathsf{Q}_{i}}^{+}$ is $\mP(1^{2},2)$, which
is the fiber of $\widetilde{\mathscr{IR}}'\to{\rm dP}_{2}$ over $\widetilde{\mathsf{q}}_{i}$.
In particular, the fiber of $\widetilde{\mathscr{IR}}'\to{\rm dP}_{2}$
over each $\widetilde{\mathsf{q}}_{i}$ is of codimension $2$. Therefore,
by \cite[Lem.5.5]{Taka3}, $\widetilde{\mathscr{IR}}$ and $\widetilde{\mathscr{IR}}'$
are actually isomorphic. Therefore, the proposition is proved.
\end{proof}
We denote by $\sO_{\widetilde{\mathscr{IR}}}(1)$ the divisor linearly
equivalent to $\{\widetilde{x}_{i}=0\}\cap\widetilde{\mathscr{IR}}$,
which we call the\textit{ tautological divisor }of $\widetilde{\mathscr{IR}}$. 
\begin{prop}
\label{prop:IRtildecanonicaldiv}The anticanonical divisor $-K_{\widetilde{\mathscr{IR}}}$
is equal to $\sO_{\widetilde{\mathscr{IR}}}(2)$. 
\end{prop}

\begin{proof}
Note that, by the proof of Proposition \ref{prop:The-Picard-numbertwo},
the birational map $\mathsf{IR}\dashrightarrow\widetilde{\mathscr{IR}}$
is contracting (any divisor on $\widetilde{\mathscr{IR}}$ corresponds
to a divisor on $\mathsf{IR}$). Therefore $-K_{\widetilde{\mathscr{IR}}}$
is the push-forward of $-K_{\mathsf{IR}}$. Let $\sO_{\mathsf{IR}}(1)$
be the tautological divisor associated to $\sU$ or $\sU'$. Since
$\det\sU=\det\sU'=-K_{\mathcal{S}}=\sO_{\mP^{5}}(1)|_{\mathcal{S}}$,
it holds that $-K_{\mathsf{IR}}=\sO_{\mathsf{IR}}(2)$. By the construction
of $\widetilde{\mathscr{IR}}\dashrightarrow\mathsf{IR}$ as in Subsection
\ref{subsec:The-birational-mapIR PU}, the push-forward of $\sO_{\mathsf{IR}}(1)$
coincides with $\sO_{\widetilde{\mathscr{IR}}}(1)$. Therefore $-K_{\widetilde{\mathscr{IR}}}=\sO_{\widetilde{\mathscr{IR}}}(2).$
\end{proof}

\subsubsection{\textbf{The midpoint $\overline{\mathscr{I}\mathscr{R}}$ and the
flopping contraction $\widetilde{g}_{\mathscr{IR}}\colon\widetilde{\mathscr{IR}}\to\overline{{\rm \mathscr{IR}}}$}}

Now we set 
\begin{equation}
\overline{\mathscr{I}\mathscr{R}}:=\overline{{\rm G}}_{\mathscr{\mathscr{I}R}}\cap\{L_{1}(\overline{r})=\cdots=L_{6}(\overline{r})=0\},\label{eq:overlineIR}
\end{equation}
which is compatible with (\ref{eq:dP2}). Note that $\overline{\mathscr{I}\mathscr{R}}$
contains $\Pi$.

We can show the following in the same way as Proposition \ref{prop:RvarGor}:
\begin{prop}
\label{prop:IRvarGor}The scheme $\overline{\mathscr{IR}}$ is a $\mQ$-Gorenstein
variety of codimension $3$ in $\mP(1^{4},2^{4})$.
\end{prop}

By (\ref{eq:KG}) with $n=4$, we have $-K_{\overline{{\rm G}}_{\mathscr{IR}}}=\sO_{\overline{{\rm G}}_{\mathscr{IR}}}(14)$.
Therefore, by (\ref{eq:overlineIR}) and Proposition \ref{prop:IRvarGor},
we obtain
\begin{equation}
-K_{\mathscr{\overline{IR}}}=\sO_{\overline{\mathscr{IR}}}(2).\label{eq:canonical-divisor-overlineIR}
\end{equation}

The birational morphism \textbf{$\psi\colon\widetilde{{\rm G}}_{\mathscr{IR}}\to\overline{{\rm G}}_{\mathscr{IR}}$
}restricts to a birational morphism $\widetilde{\mathscr{IR}}\to\overline{{\rm \mathscr{IR}}}$,
which we denote by $\widetilde{g}_{\mathscr{IR}}$.

We set $\widetilde{\Pi}:=\widetilde{\mathscr{IR}}\cap\{w=0\}$. By
the weights on the coordinates of $\widetilde{{\rm G}}_{\mathscr{IR}}$
as in Table \ref{Table:FIR}, we have the following:
\begin{lem}
\label{lem:2H-L IR}The divisor $\widetilde{\Pi}$ on $\widetilde{\mathscr{IR}}$
is linearly equivalent to $\sO_{\widetilde{\mathscr{IR}}}(2)-\widetilde{f}_{\mathscr{IR}}^{*}\sO_{{\rm dP}_{2}}(1).$
\end{lem}

\begin{prop}
\label{prop:IRflop}The birational morphism \textbf{$\widetilde{g}_{\mathscr{IR}}\colon\widetilde{\mathscr{IR}}\to\overline{{\rm \mathscr{IR}}}$}
is a flopping contraction whose nontrivial fibers are $\mP^{1}$ and
the image of the $\widetilde{g}_{\mathscr{IR}}$-exceptional locus
is a curve in $\Pi$.
\end{prop}

\begin{proof}
Let us show that any nontrivial $\widetilde{g}_{\mathscr{IR}}$-fiber
is $\mP^{1}$. Since $\psi^{-1}(\Pi)\to\Pi$ is the $\mP^{6}$-bundle
$\mP(q^{*}\sE)$ over $\Pi$ by Proposition \ref{prop:The-birational-morphism psi}
(2), a nontrivial $\widetilde{g}_{\mathscr{IR}}$-fiber is a linear
subspace of $\mP^{6}$. Since such a nontrivial $\widetilde{g}_{\mathscr{IR}}$-fiber
is mapped isomorphically onto the image in $L_{\mathscr{IR}}$, it
is also mapped isomorphically onto the image in the ${\rm dP}_{2}$.
This implies that it must be $\mP^{1}$ since the ${\rm dP}_{2}$
is factorial by Definition \ref{def:dep2} and then it does not contain
$\mP^{2}$. 

By the construction of $\psi$ as in (\ref{eq:psi}), the pull-back
of $\sO_{\overline{\mathscr{IR}}}(1)$ coincides with $\sO_{\widetilde{\mathscr{IR}}}(1).$
Therefore 
\begin{equation}
-K_{\mathscr{\widetilde{IR}}}=\widetilde{g}_{\mathscr{IR}}^{*}(-K_{\overline{\mathscr{IR}}})\label{eq:-Ktilde=00003D-Koverline}
\end{equation}
 by Proposition \ref{prop:IRtildecanonicaldiv} and (\ref{eq:canonical-divisor-overlineIR}). 

Now we show that $\widetilde{g}_{\mathscr{IR}}$ is a small contraction.
Assume for contradiction that $\widetilde{g}_{\mathscr{IR}}$ is not
a small contraction. Then, since the Picard number of $\widetilde{\mathscr{IR}}$
is $2$ by Proposition \ref{prop:The-Picard-numbertwo}, $\widetilde{g}_{\mathscr{IR}}$
is a divisorial contraction, which is crepant by (\ref{eq:-Ktilde=00003D-Koverline}).
Let $G$ be $\widetilde{g}_{\mathscr{IR}}$-exceptional divisor. Since
the Picard number of $\widetilde{\mathscr{IR}}$ is $2$, we may write
$G\equiv\sO_{\widetilde{\mathscr{IR}}}(a)+\sO_{{\rm dP}_{2}}(b)$
with some $a,b\in\mQ$. Since $G\cdot\gamma=-2$ for a $\widetilde{g}_{\mathscr{IR}}$-exceptional
curve $\gamma\simeq\mP^{1}$, we have $b=-2$. Since $G\cdot\ell$
is a positive integer for a general $\widetilde{f}_{\mathscr{IR}}$-fiber
$\ell$, so is $a$. Since $G$ is an irreducible component of $\widetilde{\Pi}$
and $G\not=\widetilde{\Pi}$, we have $a=1$ by Lemma \ref{lem:2H-L IR}.
Since $\widetilde{g}_{\mathscr{IR}}(G)$ is a surface, it holds that
$\sO_{\widetilde{\mathscr{IR}}}(1)^{3}\cdot G=0,$ thus we have $\sO_{\widetilde{\mathscr{IR}}}(1)^{4}=2\sO_{\widetilde{\mathscr{IR}}}(1)^{3}\cdot\sO_{{\rm dP}_{2}}(1)$.
Since $\sO_{\widetilde{\mathscr{IR}}}(1)=\widetilde{g}_{\mathscr{IR}}^{*}\sO_{\mathscr{\overline{IR}}}(1)$,
we have $\sO_{\widetilde{\mathscr{IR}}}(1)^{4}=\sO_{\overline{\mathscr{IR}}}(1)^{4}$.
By (\ref{eq:degG}) and Proposition \ref{prop:IRvarGor}, we have
$\sO_{\overline{\mathscr{IR}}}(1)^{4}=\nicefrac{5}{2}$. Therefore
we obtain $\sO_{\widetilde{\mathscr{IR}}}(1)^{3}\cdot\sO_{{\rm dP}_{2}}(1)=\nicefrac{5}{4}$.
This derive a contradiction; since $\Bs|\sO_{{\rm dP}_{2}}(1)|=\emptyset$,
we may assume that $\sO_{{\rm dP}_{2}}(1)$ is contained in the smooth
locus of $\widetilde{\mathscr{IR}}$ by Proposition \ref{prop:The-Picard-numbertwo}.
Then $\sO_{\widetilde{\mathscr{IR}}}(1)^{3}\cdot\sO_{{\rm dP}_{2}}(1)$
must be an integer.

Since the Picard number of $\widetilde{\mathscr{IR}}$ is 2 and $\widetilde{\mathscr{IR}}$
is $\mQ$-factorial and has only terminal singularities by Proposition
\ref{prop:The-Picard-numbertwo}, $\widetilde{g}_{\mathscr{IR}}$
is a flopping contraction. 

The argument as above implies that the image of the $\widetilde{g}_{\mathscr{IR}}$-exceptional
locus is at most a curve in $\Pi$. By changing the linear forms $L_{1}(r),\dots,L_{6}(r)$
with their linear combinations if necessary, we may assume that $\mP(q^{*}\sE)\cap\{L_{1}(\widetilde{r})=\cdots=L_{5}(\widetilde{r})=0\}$
is a $\mP^{1}$-bundle near a fiber $\gamma$ of the $\widetilde{g}_{\mathscr{IR}}$-exceptional
locus. In this situation, $L_{6}(\widetilde{r})$ can be regarded as
a global section of the rank 2 bundle associated to this $\mP^{1}$-bundle,
and the image of $\gamma$ is contained in the zero locus of this
global section. Since the zero locus is of codimension at most two,
the fiber $\gamma$ is not isolated. Therefore the image of the $\widetilde{g}_{\mathscr{IR}}$-exceptional
locus is a curve in $\Pi$. 
\end{proof}
\begin{prop}
\label{prop:IRmult2} It holds that $\overline{\mathscr{IR}}$ has
hypersurface singularities with multiplicity $2$ along the image
of the $\widetilde{g}_{\mathscr{IR}}$-exceptional locus.
\end{prop}

\begin{proof}
We can show this proposition in the same way as Proposition \ref{prop:hyp double pt}.
\end{proof}
\begin{rem}
\label{rem:general }If we choose $L_{\mathscr{IR}}$ generally, then,
by the description of the bundle $\sE$ in Proposition \ref{prop:The-birational-morphism psi},
we can show that the image of $\widetilde{g}_{\mathscr{IR}}$-exceptional
locus is a smooth curve $\Gamma_{\mathscr{IR}}$ of degree $10$ and
genus $11$, and the normal bundle of the nontrivial $\widetilde{g}_{\mathscr{IR}}$-fiber
over a point of $\Gamma_{\mathscr{IR}}$ is $\sO_{\mP^{1}}\oplus\sO_{\mP^{1}}(-1)^{\oplus2}$.
Therefore, $\widetilde{g}_{\mathscr{IR}}$ is a flopping contraction
of Atiyah type and $\overline{\mathscr{IR}}$ has ordinary double
points along $\Gamma_{\mathscr{IR}}$.
\end{rem}

\subsubsection{\textbf{The flop $\widetilde{\mathscr{I}{\rm \mathscr{R}}}\protect\dashrightarrow\widehat{\mathscr{\mathscr{I}R}}$} }

By \cite{bchm}, the flop for the flopping contraction $\widetilde{g}_{\mathscr{IR}}\colon\widetilde{\mathscr{IR}}\to\overline{{\rm \mathscr{IR}}}$
exists, which we write $\widetilde{\mathscr{IR}}\dashrightarrow\overline{\mathscr{IR}}$.
We denote by $\widehat{g}_{\mathscr{IR}}\colon\widehat{\mathscr{IR}}\to\overline{\mathscr{IR}}$
the flopped contraction.
\begin{prop}
\label{prop:flopped IR}The variety $\widehat{\mathscr{IR}}$ has
only five $\nicefrac{1}{2}$-singularities and smooth near $\widehat{g}_{\mathscr{IR}}^{-1}(\Pi)$.
Any non-trivial $\widehat{g}_{\mathscr{IR}}$-fiber is isomorphic
to $\mP^{1}$.
\end{prop}

\begin{proof}
The assertion follows from Proposition \ref{prop:IRmult2} and \cite[Prop.2.2 and Ex.2.3]{Ko}.
\end{proof}

\subsubsection{\textbf{The divisorial contraction $\widehat{f}_{\mathscr{IR}}\colon\widehat{\mathscr{IR}}\to\mathscr{IR}$}}

We denote by $\sO_{\widehat{\mathscr{IR}}}(1)$ the pull-back of the
divisor $\sO_{\overline{\mathscr{IR}}}(1)$. By (\ref{eq:canonical-divisor-overlineIR}),
we have 

\begin{equation}
-K_{\mathscr{\widehat{IR}}}=\sO_{\widehat{\mathscr{IR}}}(2).\label{eq:hatIR canonicaldiv}
\end{equation}

\begin{prop}
\label{prop:fhat IR}The strict transform $\widehat{\Pi}$ of $\Pi$
on $\widehat{\mathscr{IR}}$ is isomorphic to $\Pi\simeq\mP^{3}$
and the normal bundle of $\widehat{\Pi}$ is $\sO_{\mP^{3}}(-2)$.
There is a divisorial contraction $\widehat{f}_{\mathscr{IR}}\colon\widehat{\mathscr{IR}}\to\mathscr{IR}$
of $\widehat{\Pi}$ which contracts $\widehat{\Pi}$ to a $\nicefrac{1}{2}$-singularity.
\end{prop}

\begin{proof}
By Lemma \ref{lem:2H-L IR} and Proposition \ref{prop:IRflop}, it
holds that $\widetilde{\Pi}$ is negative for a $\widetilde{g}_{\mathscr{IR}}$-exceptional
curve. Therefore, by a property of flop, $\widehat{\Pi}$ is positive
for a $\widehat{g}_{\mathscr{IR}}$-exceptional curve. Let $\gamma$
be a line in $\Pi\simeq\mP^{3}$ which is disjoint from the image
of the $\widetilde{g}_{\mathscr{IR}}$-exceptional locus. Then $K_{\overline{\mathscr{IR}}}+\Pi|_{\Pi}=K_{\Pi}$
holds near $\gamma$. Since $K_{\overline{\mathscr{IR}}}\cdot\gamma=-2$
by (\ref{eq:canonical-divisor-overlineIR}) and $K_{\Pi}\cdot\gamma=-4$,
we have $\Pi\cdot\gamma<0$. Therefore, we have $\widehat{\Pi}\cdot\widehat{\gamma}<0$,
where $\widehat{\gamma}$ is the strict transform of $\gamma$. Since
the Picard number of $\widehat{\mathscr{IR}}$ is two, and $K_{\widehat{\mathscr{IR}}}$
is not nef, the only non-trivial morphisms from $\widehat{\mathscr{IR}}$
are $\widehat{g}_{\mathscr{IR}}$ and a $K$-negative contraction,
say, $\widehat{f}_{\mathscr{IR}}\colon\widehat{\mathscr{IR}}\to\mathscr{IR}$.
Since there exists a curve negative for $\widehat{\Pi}$ and $\widehat{\Pi}$
is $\widehat{g}_{\mathscr{IR}}$-positive, $\widehat{\Pi}$ is $\widehat{f}_{\mathscr{IR}}$-negative.
Note that $\widehat{\Pi}$ is contained in the smooth locus of $\widehat{\mathscr{IR}}$.
By (\ref{eq:hatIR canonicaldiv}), $\widehat{f}_{\mathscr{IR}}$ cannot
be a flipping contraction by \cite[(2.1) Thm.]{KA}. Thus $\widehat{f}_{\mathscr{IR}}$
is the divisorial contraction of $\widehat{\Pi}$. If the image of
the $\widehat{f}_{\mathscr{IR}}$-exceptional divisor were a surface,
then a general nontrivial $\widehat{f}_{\mathscr{IR}}$-fiber $\gamma$
would satisfy $-K_{\widehat{\mathscr{IR}}}\cdot\gamma=1$ by \cite{AW2},
a contradiction to (\ref{eq:hatIR canonicaldiv}). If the image of
the $\widehat{f}_{\mathscr{IR}}$-exceptional divisor were a curve,
then, by \cite{Taka1} and (\ref{eq:hatIR canonicaldiv}), $\widehat{\Pi}$
would be a $\mP^{2}$-bundle over a smooth curve, say, $\Delta$.
Since $\widehat{\Pi}$ is rational, $\Delta\simeq\mP^{1}$. Since
$\widehat{\Pi}$ has a contraction to $\Pi\simeq\mP^{3}$, we must
have $\widehat{\Pi}\simeq\mP(\sO(1)\oplus\sO^{\oplus2}$). Then, for
an exceptional curve $\gamma$ for $\widehat{\Pi}\to\Pi$, we have
$-K_{\widehat{\Pi}}\cdot\gamma=1$. This implies that $\widehat{\Pi}\cdot\gamma=-1$
since $-K_{\widehat{\mathscr{IR}}}\cdot\gamma=0.$ This is, however,
impossible since $\widehat{\Pi}\cdot\gamma$ is positive by the property
of $\widehat{g}_{\mathscr{IR}}$. Therefore, the image of the $\widehat{f}_{\mathscr{IR}}$-exceptional
divisor is a point. Then, $\widehat{\Pi}\to\Pi\simeq\mP^{3}$ is a
finite and birational morphism, and hence is an isomorphism by the
Zariski main theorem. Let $\gamma$ be a line in $\widehat{\Pi}$.
Then, by (\ref{eq:hatIR canonicaldiv}), we have 

\[
0>\widehat{\Pi}\cdot\gamma=K_{\widehat{\Pi}}\cdot\gamma+(-K_{\widehat{\mathscr{IR}}})\cdot\gamma=-4+\sO_{\widehat{\mathscr{IR}}}(2)\cdot\gamma.
\]
Thus we have $\sO_{\widehat{\mathscr{IR}}}(1)\cdot\gamma=1$ and $\widehat{\Pi}\cdot\gamma=-2$.
Thus, $\widehat{f}_{\mathscr{IR}}$ contracts $\widehat{\Pi}$ to
a $\nicefrac{1}{2}$-singularity.
\end{proof}
\begin{rem}
If we choose $L_{\mathscr{IR}}$ generally as in Remark \ref{rem:general },
then Proposition \ref{prop:fhat IR} immediately follows from the
construction of the flop $\widetilde{\mathscr{IR}}\dashrightarrow\widehat{\mathscr{IR}}$
of Atiyah type. 
\end{rem}

\subsubsection{\textbf{The variety $\mathscr{IR}$ and its linear section }}
\begin{prop}
\label{prop:CoseqIR}The variety $\mathscr{IR}$ is a terminal $\mQ$-factorial
$\mQ$-Fano variety with Picard number $1$ and with only six $\nicefrac{1}{2}$-singularities.
It holds that $-K_{\mathscr{R}}=\sO_{\mathscr{R}}(2)$, $\sO_{\mathscr{R}}(1)^{10}=3$
and $\sO_{\mathscr{R}}(1)$ is numerically primitive.
\end{prop}

\begin{proof}
The Picard number of the variety $\widehat{\mathscr{IR}}$ is $2$
by Proposition \ref{prop:The-Picard-numbertwo} since the flop \textbf{$\widetilde{\mathscr{I}{\rm \mathscr{R}}}\dashrightarrow\widehat{\mathscr{\mathscr{I}R}}$}
does not change the Picard numbers. Therefore the Picard number of
the variety $\mathscr{IR}$ is $1$. It holds that $\mathscr{IR}$
has only six $\nicefrac{1}{2}$-singularities by Propositions \ref{prop:flopped IR}
and \ref{prop:fhat IR}. The equality $-K_{\mathscr{IR}}=\sO_{\mathscr{R}}(2)$
follows from (\ref{eq:hatIR canonicaldiv}).\textbf{ }The remaining
facts follow in the same way as the proof of Proposition \ref{prop:CoseqR}.
\end{proof}
Now we have finished the constructions of $\mathscr{IR}$ and the Sarkisov diagram (\ref{eq:SarkisovIR}).

We can show the following corollary in the same way as Corollary \ref{cor:numX}: 
\begin{cor}
\label{cor: NumX IR} Let $X$ be a linear section of $\mathscr{IR}$
with only six $\nicefrac{1}{2}$-singularities. Then $X$ is a $\mQ$-Fano
$3$-fold with $(-K_{X})^{3}=3$ and $g(X)=1$ such that $-K_{X}$
is numerically primitive. 
\end{cor}

\section{\textbf{Proof of Theorem \ref{thm:main} (2)\label{sec:Proof-of-Theorem(2)}}}

\subsection{Preliminary results}
\begin{lem}
\label{lem:unique-IR} For any ${\rm dP}_{2}$, we adopt the identifications
as in Lemma \ref{lem:Important}. Let $\Pi_{C}$ be a plane in $\mathcal{Q}_{B}$
such that $C_{\mathcal{S}}:=\mathcal{S}\cap\Pi_{C}$ is a cubic elliptic
curve. Let $C$ be the birational image of $C_{\mathcal{S}}$ on the
${\rm dP}_{2}$ by the birational map ${\rm dP}_{2}\dashrightarrow\mathcal{S}$,
which is an quartic elliptic curve through $\mathsf{\widetilde{q}}_{1},\dots,\mathsf{\widetilde{q}}_{5}$
by Proposition \ref{prop:CS cubic} (B). We take the divisor\textbf{
$\mathsf{X}=\{\sum_{i=1}^{4}a_{i}\mathsf{x}_{i}=0\}\subset\mathsf{IR}$}
for $\Pi_{C}$ as in Proposition \ref{prop:YS} and\textbf{ }let $\widetilde{X}=\{\sum_{i=1}^{4}a_{i}\widetilde{x}_{i}=0\}\subset\widetilde{\mathscr{IR}}$
be the divisor with the same equation as $\mathsf{X}$.

\begin{enumerate}[$(A)$]

\item The following hold:

\begin{enumerate}[$(1)$]

\item $\widetilde{X}$ is irreducible and normal,

\item $-K_{\widetilde{X}}$ is $\mQ$-Cartier and relatively ample
over ${\rm dP}_{2}$, and 

\item the induced morphism $\widetilde{X}\to{\rm dP}_{2}$ has only
$1$-dimensional non-trivial fibers and coincides with the blow-up
along $C$ outside the fibers over $\mathsf{\widetilde{q}}_{1},\dots,\mathsf{\widetilde{q}}_{5}$. 

\end{enumerate} 

\item If $\widetilde{X}'$ and the morphism $\widetilde{X}'\to{\rm dP}_{2}$
satisfies the same properties as $\widetilde{X}$ and $\widetilde{X}\to{\rm dP}_{2}$
stated in $(A)$ $(1)$--$(3),$ then $\widetilde{X}$ and $\widetilde{X}'$are
isomorphic over ${\rm dP}_{2}$.

\end{enumerate}
\end{lem}

\begin{proof}
(A) (3). By Proposition \ref{prop:YS}, the naturally induced morphism
$\mathsf{X}\to\mathcal{S}$ is the blow-up along $C_{\mathcal{S}}$.
Therefore $\widetilde{X}\to{\rm dP}_{2}$ is the blow-up along $C$
possibly outside the inverse image of $\cup_{1\leq i<j\leq5}\rho_{\mP}^{-1}(\ell_{ij})$
by Proposition \ref{prop:IRtilde and P(U)S}. Let us check that $\widetilde{X}\to{\rm dP}_{2}$
is an isomorphism near the fiber over any point of $\cup_{1\leq i<j\leq5}\rho_{\mP}^{-1}(\ell_{ij}\setminus\{\mathsf{q}_{i},\mathsf{q}_{j})\}$.
If $\widetilde{X}\to{\rm dP}_{2}$ is not an isomorphism near the
fiber over a point of $\cup_{1\leq i<j\leq5}\rho_{\mP}^{-1}(\ell_{ij}\setminus\{\mathsf{q}_{i},\mathsf{q}_{j})\}$,
then $\mathsf{X}\to\mathcal{S}$ is not an isomorphism at a point
in $\rho_{{\rm S}}^{-1}(\mathsf{q}_{ij})$ by Proposition \ref{prop:IRtilde and P(U)S}.
This implies that $C_{\mathcal{S}}$ contains a point in $\rho_{{\rm S}}^{-1}(\mathsf{q}_{ij})$.
This is a contradiction since $\rho_{{\rm S}}^{-1}(\mathsf{q}_{ij})\subset\Sing\mathcal{S}$
by Proposition \ref{prop:Sing calS} but we can show that $C_{\mathcal{S}}$
does not pass through a singular point of $\mathcal{S}$ in the same
way as Proposition \ref{prop:CS cubic} (A). It remains to show that
the fibers of $\widetilde{X}\to{\rm dP}_{2}$ over $\widetilde{\mathsf{q}}_{1},\dots,\widetilde{\mathsf{q}}_{5}$
are $1$-dimensional. Assume for contradiction that the fiber of $\widetilde{X}\to{\rm dP}_{2}$
over a point $\widetilde{\mathsf{q}}_{i}$ is not $1$-dimensional.
Then, by Proposition \ref{prop:five1/2}, the fiber is $2$-dimensional,
and $\widetilde{X}$ contains an irreducible component of the fiber
of $\widetilde{\mathscr{IR}}\to{\rm dP}_{2}$ over $\widetilde{\mathsf{q}}_{i}$.
Then, by the construction of the birational map $\mP(\sU|_{\mathcal{S}})\dashrightarrow\widetilde{\mathscr{IR}}$
as in the proof of Proposition \ref{prop:The-Picard-numbertwo}, we
see that $\mP(\sU|_{\mathsf{Q}_{i}})\cap\mathsf{X}$ contains $\Delta:=\mP(\sO_{\mP^{1}\times\mP^{1}}(-1,0)\oplus0)$
or $\mP(0\oplus\sO_{\mP^{1}\times\mP^{1}}(0,-1))$ if $\mathsf{Q}_{i}\simeq\mP^{1}\times\mP^{1}$
(resp.$\,$$\Delta$ as in the proof if $\mathsf{Q}_{i}\simeq\mP(1^{2},2)$).
Since both of $\Delta$ and $\mP(\sU|_{\mathsf{Q}_{i}})\cap\mathsf{X}$
are generically a section of $\mP(\sU|_{\mathsf{Q}_{i}})$ and $\Delta\not=\mP(\sU|_{\mathsf{Q}_{i}})\cap\mathsf{X}$,
we see that $\mP(\sU|_{\mathsf{Q}_{i}})\cap\mathsf{X}$ contains an
irreducible component which is mapped to a curve in $\mathsf{Q}_{i}$.
This is a contradiction since then $\mathsf{X}\to\mathcal{S}$ cannot
be the blow-up along $C_{\mathcal{S}}$. 

(1) and (2). By the same argument as the last paragraph of the proof
of Lemma \ref{lem:unique}, we see that $\widetilde{X}$ is irreducible
and normal.\textbf{ }Since $\widetilde{X}$ is a tautological divisor,
so is $-K_{\widetilde{X}}$ by Proposition \ref{prop:IRtildecanonicaldiv}.\textbf{
}Then this is relatively ample for $\widetilde{X}\to{\rm {\rm dP}_{2}}$. 

(B) The last assertion follows from \cite[Lem.5.5]{Taka3} since $\widetilde{X}$
and $\widetilde{X}'$ are isomorphic in codimension $1$.
\end{proof}
\begin{lem}
\textbf{\label{lem:RestIRToX}} Let $\widetilde{X}$ be a tautological
divisors of $\widetilde{\mathscr{IR}}$. Assume that

\begin{enumerate}[$(i)$]

\item $\widetilde{X}$ is an irreducible normal $3$-fold with only
five $\nicefrac{1}{2}$-singularities, 

\item $\widetilde{X}$ has Picard number $2$, and

\item the birational morphism $\widetilde{g}\colon\widetilde{X}\to\overline{X}$
induced from $\widetilde{g}_{\mathscr{IR}}$ is a flopping contraction.

\end{enumerate}

Then the strict transform $X\subset\mathscr{IR}$ of $\overline{X}$
is a Type $I\!R$ $\mQ$-Fano $3$-fold with respect to the point
$\widehat{f}_{\mathscr{IR}}(\widehat{\Pi})$, and the Sarkisov diagram
(\ref{eq:SarkisovIR}) restricts to (\ref{eq:3-foldSarkisov}).
\end{lem}

\begin{proof}
We can show the assertion in the same way as the proof of Lemma \ref{lem:RestRToX}
except showing that the restriction $\widetilde{f}\colon\widetilde{X}\to{\rm dP}_{2}$
of $\widetilde{f}_{\mathscr{IR}}$ has the desired properties. Since
$\dim\widetilde{X}=\dim Y$, $\widetilde{f}$ is birational. Let $\widetilde{E}$
be the $\widetilde{f}$-exceptional divisor. If $\widetilde{f}$ contracts
$\widetilde{E}$ to a point, then, $\widetilde{X}$ cannot have five
$\nicefrac{1}{2}$-singularities by the classification of $\widetilde{f}$
as in \cite{Ka1,Ka2} since ${\rm dP}_{2}$ has only Gorenstein singularities,
a contradiction. Thus $\widetilde{f}$ contracts $\widetilde{E}$
to a curve, which we denote by $C$. Let $\widetilde{L}$ be the $\widetilde{f}$-pull-back
of the ample generator of ${\rm Pic}\,{\rm dP}_{2}$ and $\widehat{E}_{\widetilde{X}}:=\widetilde{\Pi}\cap\widetilde{X}$.
Note that\textbf{ }$\widehat{E}_{\widetilde{X}}$ is the strict transform
of $\widehat{E}$. By Lemma \ref{lem:2H-L IR}, it holds that 
\begin{equation}
\widetilde{L}=2(-K_{\widetilde{X}})-\widehat{E}_{\widetilde{X}}.\label{eq:L=00003D2-K-E IR}
\end{equation}
Note that the equalities (\ref{eq:KE})\textbf{ }hold in this situation.
Therefore, since $\widetilde{L}^{3}=2$, we have
\[
2=\left(2(-K_{\widetilde{X}})-\widehat{E}_{\widetilde{X}}\right)^{3}=-4-\widehat{E}_{\widetilde{X}}^{3},\text{i.e.,\,\,}\widehat{E}_{\widetilde{X}}^{3}=-6.
\]
By (\ref{eq:L=00003D2-K-E IR}) and $-K_{\widetilde{X}}=2\widetilde{L}-\widetilde{E}$,
we have $\widetilde{E}=3(-K_{\widetilde{X}})-2\widehat{E}_{\widetilde{X}}$,
and then
\begin{align*}
 & -\widetilde{E}^{2}L=-\left(3(-K_{\widetilde{X}})-2\widehat{E}_{\widetilde{X}}\right){}^{2}\left(2(-K_{\widetilde{X}})-\widehat{E}_{\widetilde{X}}\right)=4,
\end{align*}
which implies that $\deg C=4$. We also have 
\end{proof}
\[
(-K_{\widetilde{E}})^{2}=(-K_{\widetilde{X}}-\widetilde{E})^{2}\widetilde{E}=\left(2(-K_{\widetilde{X}})-2\widehat{E}_{\widetilde{X}}\right){}^{2}\left(3(-K_{\widetilde{X}})-2\widehat{E}_{\widetilde{X}}\right)=-10,
\]
which is equal to $8\left(1-p_{g}(C)\right)-2m-18m'$ as in \cite[Prop.7.1]{Taka3},
where $m$ and $m'$ are the numbers of fibers of type (A) and (B)
respectively as in there. Since $Y={\rm dP}_{2}$ has five $A_{\leq2}$-singularities,
all the five $\nicefrac{1}{2}$-singularities of $\widetilde{X}$
are contained in $\widetilde{E}$, and it holds that $m=5$, $m'=0$
and $p_{g}(C)=1$. Since $m'=0$, $C$ is a smooth curve, hence $C$
is an elliptic curve. Therefore $X$ is of Type $I\!R$ with respect
to the $\nicefrac{1}{2}$-singularity $\widehat{f}_{\mathscr{IR}}(\widehat{\Pi})$.

\subsection{Proof of Theorem \ref{thm:main} (2-1)}

We can show this in the same way as in the case of Type $R$ using
Lemmas \ref{lem:unique-IR} and \ref{lem:RestIRToX}.$\hfill\square$

\subsection{Proof of Theorem \ref{thm:main} (2-2)}
\begin{lem}
\label{lem:Construct-example-IR}There exists a linear section of
$\mathscr{IR}$ of codimension $1$ that is a Type $I\!R$ $\mQ$-Fano
$3$-fold.
\end{lem}

\begin{proof}
We construct such an example starting from the plane $\Pi_{C}$ in
$\mathcal{Q}_{B}$ as in the proof of Proposition \ref{prop:GenExOf f'-IR}.
We use the notation as in there. Since $\Pi_{C}$ as in the proposition
satisfies the assumptions as in Lemma \ref{lem:unique-IR}, and $\widetilde{X}$
and the morphism $\widetilde{f}\colon\widetilde{X}\to{\rm dP}_{2}$
constructed in the proposition satisfies the assumption of $\widetilde{X}'$
and the morphism $\widetilde{X}'\to{\rm dP}_{2}$ in the lemma (B),
$\widetilde{X}$ is a tautological divisor of $\widetilde{\mathscr{IR}}$.
Let $\overline{X}$ be the $\widetilde{g}_{\mathscr{IR}}$-image of
$\widetilde{X}\subset\overline{\mathscr{IR}}$, and we denote by \textit{$\widetilde{g}\colon\widetilde{X}\to\overline{X}$
}the restriction of $\widetilde{g}_{\mathscr{IR}}$. An important
point to note is that the construction in the proof of the proposition
shows that the Picard number of $\widetilde{X}$ is $2$. By the explicit
construction as in the proposition, we are ready to apply Lemma \ref{lem:RestIRToX}
once we check the following:

\vspace{5pt}

\noindent Claim. \textit{It holds that $\widetilde{g}\colon\widetilde{X}\to\overline{X}$
is a flopping contraction.}

\vspace{5pt}

\noindent \textit{Proof.} Let $\widetilde{L}$ be the $\widetilde{f}$-pull-back
of the ample generator of ${\rm Pic}\,{\rm dP}_{2}$. Suppose for
contradiction that $\widetilde{g}$ is not a flopping contraction.
Then it is a crepant divisorial contraction. Let $\widetilde{G}$
be the $\widetilde{g}$-exceptional divisor. As in the same way for
Type $R$, we may write $\widetilde{G}\equiv p(-K_{\widetilde{X}})-q\widetilde{L}$
with some $p,q\in\mQ$, and we have 
\begin{equation}
0=5/2\,p-q(-K_{\widetilde{X}})^{2}\widetilde{L}.\label{eq:pq-2}
\end{equation}
We note that $\widetilde{E}\widetilde{L}^{2}=0$ and $\widetilde{E}^{2}\widetilde{L}=-\deg C=-4$.
Then, since $-K_{\widetilde{X}}=2\widetilde{L}-\widetilde{E}$, we
obtain

\begin{align}
 & (-K_{\widetilde{X}})^{2}\widetilde{E}+2(-K_{\widetilde{X}})\widetilde{E}^{2}+\widetilde{E}^{3}=0,\label{eq:eqA}\\
 & (-K_{\widetilde{X}})\widetilde{E}^{2}+\widetilde{E}^{3}=-8.\label{eq:eqB}
\end{align}
Since $(-K_{\widetilde{E}})^{2}=(-K_{\widetilde{X}}-\widetilde{E})^{2}\widetilde{E}=-10$
by Proposition \ref{prop:GenExOf f'}, we have 
\begin{equation}
(-K_{\widetilde{X}})^{2}\widetilde{E}-2(-K_{\widetilde{E}})\widetilde{E}^{2}+\widetilde{E}^{3}=-10.\label{eq:eqC}
\end{equation}
Therefore, by (\ref{eq:eqA}), (\ref{eq:eqB}), and (\ref{eq:eqC}),
we obtain $(-K_{\widetilde{X}})^{2}\widetilde{E}=\nicefrac{11}{2}$.
Since $-K_{\widetilde{Y}}=2\widetilde{L}-\widetilde{E}$, we have
\[
(-K_{\widetilde{X}})^{2}\widetilde{L}=\nicefrac{1}{2}\left((-K_{\widetilde{X}})^{3}+(-K_{\widetilde{X}})^{2}\widetilde{E}\right)=8.
\]
From this and (\ref{eq:pq-2}), we have 
\begin{equation}
\widetilde{G}\equiv p\left((-K_{\widetilde{X}})-\nicefrac{5}{16}\widetilde{L}\right).\label{eq:Gpq-3}
\end{equation}
 Pushing this down to ${\rm dP}_{2}$, we have $\widetilde{f}{}_{*}\widetilde{G}\equiv\nicefrac{27p}{16}\widetilde{L}.$
Thus $m:=\nicefrac{27p}{16}$ is a natural number. By Proposition
\ref{prop:flopped IR}, $\widetilde{g}(\widetilde{G})$ is a curve,
hence we have $\widetilde{G}\cdot\gamma=-2$ for a general $\widetilde{g}|_{\widetilde{G}}$-fiber
$\gamma$. Thus, by (\ref{eq:Gpq-3}), it holds that 
\[
-2=-\nicefrac{5p}{16}(\widetilde{L}\cdot\gamma)=-\nicefrac{5m}{27}(\widetilde{L}\cdot\gamma),
\]
which has no integral solutions for $m$ and $\widetilde{L}\cdot\gamma$,
a contradiction.$\hfill\square$

Therefore, by Lemma \ref{lem:RestIRToX}, we obtain a Type $I\!R$
$\mQ$-Fano $3$-fold $X$ with respect to the $\nicefrac{1}{2}$-singularity
$\widehat{f}_{\mathscr{IR}}(\widehat{\Pi})$ as desired. 
\end{proof}
\noindent \textit{Proof of Theorem} \ref{thm:main} \textit{(2-2).}
We assume that $X$ is a linear section of $\mathscr{IR}$ of codimension
$1$ with only six $\nicefrac{1}{2}$-singularities. By \cite[Thm.4.1]{To}
and Lemma \ref{lem:Construct-example-IR}, $X$ has Picard number
$1$. Therefore, $X$ is a prime $\mQ$-Fano $3$-fold of $g(X)=1$
by Corollary \ref{cor: NumX IR}. 

Let $\widehat{X}\subset\widehat{\mathscr{IR}}$, $\overline{X}\subset\overline{\mathscr{IR}}$,
and $\widetilde{X}\subset\widetilde{\mathscr{IR}}$ be the strict
transforms of $X$. In the same way as the proof of Theorem \ref{thm:main}
(1-2), we see that $\widetilde{X}$ is a tautological divisor of $\widetilde{\mathscr{IR}}$,
the restriction $\widehat{f}\colon\widehat{X}\to X$ of $\widehat{f}_{\mathscr{IR}}$
is the blow-up at $\widehat{f}_{\mathscr{IR}}(\widehat{\Pi}),$ the
restriction $\widehat{g}\colon\widehat{X}\to\overline{X}$ of $\widehat{g}_{\mathscr{IR}}$
is a flopping contraction, and $\widehat{X}\dashrightarrow\widetilde{X}$
is the flop for $\widehat{g}$. Therefore the assumptions (i)--(iii)
of Lemma \ref{lem:RestIRToX} are satisfied, and then $X$ is of Type
$IR$ with respect to the $\nicefrac{1}{2}$-singularity $\widehat{f}_{\mathscr{IR}}(\widehat{\Pi})$. 

$\hfill\square$

\end{document}